\documentclass[12pt, eqno]{article}
\usepackage{amssymb,latexsym}
\usepackage{amsmath,latexsym}

\newcommand{\K}{{\mathbb K}}
\newcommand{\F}{{\mathbb F}}
\newcommand{\N}{{\mathbb N}}

\newcommand{\R}{{\mathbb R}}

\newcommand{\Z}{{\mathbb Z}}

\newtheorem{theorem}{Theorem}[section]

\newtheorem{corollary}[theorem]{Corollary}

\newtheorem{remark}[theorem]{Remark}

\newtheorem{lemma}[theorem]{Lemma}
\newtheorem{question}[theorem]{Question}

\newtheorem{proposition}[theorem]{Proposition}
\newtheorem{claim}[theorem]{Claim}

\begin{document}


\title{The Conley conjecture for
Hamiltonian systems\\
 on the cotangent bundle\\
  and its analogue for
Lagrangian systems }

\author{\\Guangcun Lu
\thanks{Partially supported by the NNSF  10671017 of
China and the Program for New Century Excellent Talents of the
Education Ministry of China.}\\
{\normalsize School of Mathematical Sciences, Beijing Normal University},\\
{\normalsize Laboratory of Mathematics
 and Complex Systems,  Ministry of
  Education},\\
  {\normalsize Beijing 100875, The People's Republic
 of China}\\
{\normalsize (gclu@bnu.edu.cn)}}
\date{Revised, 29 June,  2008}
\maketitle \vspace{-0.1in}

\abstract{In this paper, the Conley conjecture,  which were recently
proved by
 Franks and Handel \cite{FrHa} (for
surfaces of positive genus), Hingston \cite{Hi} (for tori) and
Ginzburg \cite{Gi} (for closed symplectically aspherical manifolds),
is proved for $C^1$-Hamiltonian systems on the cotangent bundle of a
$C^3$-smooth compact manifold $M$ without boundary, of a time
$1$-periodic $C^2$-smooth Hamiltonian $H:\R\times T^\ast M\to\R$
which is strongly convex and has quadratic growth on the fibers.
Namely,  we show that such a Hamiltonian system  has an infinite
sequence of contractible integral periodic solutions such that any
one of them cannot be obtained from others by iterations.
 If $H$ also satisfies $H(-t,q, -p)=
H(t,q, p)$ for any $(t,q, p)\in\R\times T^\ast M$, it is shown that
the time-one map of the
 Hamiltonian system  (if exists) has infinitely many periodic points siting in
the zero section of $T^\ast M$. If $M$ is $C^5$-smooth and $\dim
M>1$, $H$ is of $C^4$ class and independent of time $t$, then for
any $\tau>0$ the corresponding system has an infinite sequence of
contractible  periodic solutions of periods of integral multiple of
$\tau$ such that any one of them cannot be obtained from others by
iterations or rotations. These results are obtained by proving
similar results for the Lagrangian system of the Fenchel transform
of $H$, $L:\R\times TM\to\R$, which is proved to be strongly convex
and to have quadratic growth in the velocities yet.} \vspace{-0.1in}
\medskip
\tableofcontents

\section{Introduction and main results}
\setcounter{equation}{0}

Recently, a remarkable progress in Symplectic geometry and
Hamiltonian dynamics is that the Conley conjecture \cite{Co, SaZe}
were proved by  Franks and Handel \cite{FrHa} (for surfaces of
positive genus,  also see \cite{Le} for  generalizations to
Hamiltonian homeomorphisms), Hingston \cite{Hi} (for tori) and
Ginzburg \cite{Gi} (for closed symplectically aspherical manifolds).
See \cite{FrHa, Le, Gi} and references therein for a detailed
history and related studies.

In this paper we always assume that $M$ is a $n$-dimensional,
connected \textsf{$C^3$-smooth compact manifold without boundary}
without special statements.  For a time $1$-periodic $C^2$-smooth
Hamiltonian $H:\R\times T^\ast M\to\R$, let $X_H$ be the Hamiltonian
vector field of $H$ with respect to the standard symplectic
structure on $T^\ast M$, $\omega_{\rm can}:=-dq\wedge dp$ in local
coordinates $(q,p)$ of $T^\ast M$, that is, $\omega(X_H(t,q,p),\xi)=
- dH(t,q,p)(\xi)\;\forall\xi\in T_{(q,p)}T^\ast M$. Unlike the case
of compact symplectic manifolds we only consider subharmonic
solutions of the Hamiltonian equations
\begin{equation}\label{e:1.1}
\dot{x}(t)=X_H(t,x(t))
\end{equation}
for $C^2$-smooth Hamiltonians $H:\R\times T^\ast M\to\R$ satisfying
the following conditions (H1)-(H3):
\begin{description}
\item[(H1)] $H(t+1, q, p)= H(t, q, p)$ for all $(t,q, p)\in\R\times
T^\ast M$.
\end{description}
 In any local coordinates $(q_1,\cdots, q_n)$, there exist
constants $0<C_1<C_2$, depending on the local coordinates, such that
\begin{description}
\item[(H2)] $C_1|{\bf u}|^2\le \sum_{ij}\frac{\partial^2 H}{\partial p_i\partial p_j}(t,q,
p)u_iu_j\le C_2|{\bf u}|^2\quad\forall {\bf
u}=(u_1,\cdots,u_n)\in\R^n$,

\item[(H3)] $\Bigl|\frac{\partial^2 H}{\partial
q_i\partial p_j}\bigl(t, q, p\bigr)\Bigr| \le C_2 (1+ |p|)$,\;
$\Bigl|\frac{\partial^2 H}{\partial q_i\partial q_j}\bigl(t, q,
p\bigr)\Bigr| \le C_2\Bigl(1+ |p|^2\Bigr)$.
 \end{description}
A class of important examples of such Hamiltonians are Physical
Hamiltonian (including $1$-periodic potential
 and electromagnetic forces in time ) of the form
\begin{equation}\label{e:1.2}
H(t,q,p)=\frac{1}{2}\|p-A(t,q)\|^2+ V(t,q)
\end{equation}

For $C^r$-smooth Hamiltonians $H:\R\times T^\ast M\to\R$ satisfying
the conditions (H1)-(H3), $r\ge 2$, by the inequality in the left
side of the condition (H2),   we can use the inverse Legendre
transform to get a
 fiber-preserving $C^{r-1}$-diffeomorphism
\begin{equation}\label{e:1.3}
\frak{L}_H: \R/\Z\times T^\ast M\to \R/\Z\times TM, \quad (t, q,
p)\mapsto \left(t, q, D_pH(t,q, p)\right),
\end{equation}
and a $C^r$-smooth function $L:\R\times TM\to\R$:
\begin{eqnarray}\label{e:1.4}
L(t,q, v)&=& \max_{p\in T_qM}\{\langle p, v\rangle-
H(t,q, p)\}\nonumber\\
&=& \langle p(t,q, v), v\rangle- H(t,q, p(t,q, v)),
\end{eqnarray}
where $p=p(t,q, v)$ is a unique point determined by the equality
$v=D_pH(t,q, p)$. (See (\cite[Prop.2.1.6]{Fa})). By (\ref{e:1.4}) we
have
\begin{description}
\item[(L1)] $L(t+1, q, v)=L(t, q, v)$ for all $(t,q,v)\in\R\times
TM$.
\end{description}
It is easily checked that the corresponding $L$ with the physical
Hamiltonian in (\ref{e:1.2}) is given by
$$
L(t,q,v) = \frac{1}{2}\|v\|^2+ \langle A(t,q), v\rangle - V(t,q).
$$
In Appendix we shall prove\vspace{2mm}

\noindent{\bf Proposition A.}\quad{\it Under the condition {\rm
(H1)}, {\rm (H2)} is equivalent to the following {\rm (L2)} plus the
third inequality in {\rm (L3)}, and  ${\rm (H2)+(H3)}\Leftrightarrow
{\rm (L2)+(L3)}$. }\vspace{2mm}\\
In any local coordinates $(q_1,\cdots, q_n)$, there exist constants
$0<c<C$, depending on the local coordinates, such that
\begin{description}
\item[(L2)] $\sum_{ij}\frac{\partial^2L}{\partial v_i\partial
v_j}(t,q,v)u_iu_j\ge c|{\bf u}|^2\quad\forall {\bf
u}=(u_1,\cdots,u_n)\in\R^n$,

\item[(L3)] $\Bigl| \frac{\partial^2L}{\partial q_i\partial
q_j}(t,q,v)\Bigr|\le C(1+ |v|^2),\quad \Bigl|
\frac{\partial^2L}{\partial q_i\partial v_j}(t,q,v)\Bigr|\le C(1+
|v|),\quad\hbox{and}\\
 \Bigl| \frac{\partial^2L}{\partial v_i\partial
v_j}(t,q,v)\Bigr|\le C$.
\end{description}
(One can also write these two conditions in the free coordinates,
see \cite[\S 2]{AbSc}.) So Proposition A shows that the conditions
(L2)-(L3) have the same properties as (H2)-(H3). ({\it Note}: we do
not claim that the condition (H2) (resp. (H3)) is equivalent to (L2)
(resp.(L3)).) By (L2), the Legendre transform produces the inverse
of $\frak{L}_H$,
 \begin{equation}\label{e:1.5}
\frak{L}_L: \R/\Z\times TM\to \R/\Z\times T^\ast M, \quad (t, q,
v)\mapsto \left(t, q, D_vL(t,q, v)\right),
\end{equation}
and $H$
and $L$ are related by :
$$
H(t,q,p)= \langle p, v(t,q,p)\rangle- L(t,q, v(t,q,p)),
$$
where $v=v(t,q,p)$ is a unique point determined by the equality
$p=D_vL(t,q,v)$. In this case, it is well-known that a curve $\R\to
T^\ast M,\;t\mapsto x(t)=(\gamma(t), \gamma^\ast(t))$ is a solution
of (\ref{e:1.1}) if and only if
$\gamma^\ast(t)=D_vL(t,\gamma(t),\dot\gamma(t))\;\forall t\in\R$ and
$\gamma$ is a solution of  the Lagrangian system on $M$:
\begin{equation}\label{e:1.6}
\frac{d}{dt}\Big(\frac{\partial
L}{\partial\dot{q_i}}\Big)-\frac{\partial L}{\partial q_i}=0
\end{equation}
in any local coordinates $(q_1,\cdots, q_n)$.

Hence we only need to study the existence of infinitely many
distinct integer periodic solutions of the system (\ref{e:1.6})
under the assumptions (L1)-(L3). To describe our results we
introduce the following \textsf{notations and notions}.

For any $T>0$, each map in $C(\R/T\Z, M)$ represent a
 homotopy class of free loops in $M$. As topological
 spaces $C(\R/T\Z, M)$ and $C(\R/\Z, M)$ are always homeomorphic.
For a homotopy class $\alpha$ of free loops in $M$, denote by
$C(\R/T\Z, M;\alpha)$  the subset of maps in $C(\R/T\Z, M)$
representing $\alpha$. For $k\in\N$, if we view $\gamma\in C(\R/T\Z,
M;\alpha)$ as a $T$-periodic map $\gamma:\R\to M$, it is also viewed
as a $kT$-periodic map from $\R$ to $M$ and thus yields an element
of $C(\R/kT\Z, M)$, \textsf{called the $k$-th iteration of} $\gamma$
and \textsf{denoted by} $\gamma^k$. This $\gamma^k\in C(\R/kT\Z, M)$
 represents a free homotopy class in $M$,
\textsf{denoted by} $\alpha^k$. So $\gamma^k\in C(\R/kT\Z,
M;\alpha^k)$. Note also that topological
 spaces $C(\R/T\Z, M;\alpha)$ and $C(\R/\Z, M;\alpha)$ are always
 homeomorphic yet. For $m\in\N$ let $C^m(\R/T\Z, M)$ denote the subset
 of all $C^m$-loops $\gamma:\R/T\Z\to M$.

A periodic map $\gamma:\R\to M$ is called {\bf reversible} (or {\bf
even}) if $\gamma(-t)=\gamma(t)$ for any $t\in\R$. \textsf{Note that
such a map is always contractible}! For $\gamma\in C(\R/T\Z, M)$  we
define \textsf{rotations} of $\gamma$ via $s\in\R$ as maps
$s\cdot\gamma:\R\to M$  defined by $s\cdot\gamma(t)=\gamma(t+s)$ for
$t\in\R$. Then $s\cdot\gamma\in C(\R/T\Z, M)$ and
$(s\cdot\gamma)^m=s\cdot\gamma^m$ for any $s\in\R$ and $m\in\N$. We
 call the set
 $$
 \{\gamma^m\}_{m\in\N}\quad \left({\rm resp}.\;  \{s\cdot\gamma^m\}^{s\in\R}_{m\in\N}\;\right)
 $$
 a $T$-{\bf periodic map tower} (resp.  $T$-{\bf periodic orbit tower})
 based on $\gamma$ (a $T$-periodic map from $\R$ to $M$).
A $T_1$-periodic map tower  $\{\gamma_1^m\}_{m\in\N}$ (resp.
$T_1$-periodic orbit tower $\{s\cdot\gamma_1^m\}^{s\in\R}_{m\in\N}$)
based on a $T_1$-periodic map $\gamma_1:\R\to M$ is called {\bf
distinct} with $\{\gamma^m\}_{m\in\N}$ (resp.
$\{s\cdot\gamma^m\}_{m\in\N}$) if there is no $\tau$-periodic map
 $\beta:\R\to M$ such that $\gamma=\beta^p$ and $\gamma_1=\beta^q$
for some $p,q\in\N$ (resp. $\gamma=s\cdot\beta^p$ and
$\gamma_1=s'\cdot\beta^q$ for some $p,q\in\N$ and $s, s'\in\R$).
When $\gamma$ is contractible as a map from $\R/T\Z$ to $M$, we call
the $T$-{periodic map tower} $\{\gamma^m\}_{m\in\N}$ (resp.
$T$-{periodic orbit tower} $\{s\cdot\gamma^m\}^{s\in\R}_{m\in\N}$)
{\bf contractible}.

 For $\tau\in\N$, if $\gamma:\R\to M$ is  a
$\tau$-periodic solution of (\ref{e:1.6}), we call the set
$\{\gamma^m\}_{m\in\N}$
 a $\tau$-{\bf periodic solution tower} of (\ref{e:1.6}) based on $\gamma$.
Two periodic solution towers of (\ref{e:1.6}) are said to be {\bf
distinct} if they are distinct as periodic map towers. Furthermore,
if $s\cdot\gamma$ is also a $\tau$-periodic solution of
(\ref{e:1.6}) for any $s\in\R$, (for example, in the case $L$ is
independent of $t$), we call $\{s\cdot\gamma^m\}^{s\in\R}_{m\in\N}$
a $\tau$-{\bf periodic solution orbit tower} of (\ref{e:1.6}). When
two periodic solution orbit towers are distinct as periodic orbit
towers we call them {\bf distinct} {\bf periodic solution orbit
towers} of (\ref{e:1.6}) based on $\gamma$. Clearly, the existence
of infinitely many distinct integer periodic solution towers (resp.
solution orbit towers) of (\ref{e:1.6}) implies that there exist an
infinite sequence of integer periodic solutions of (\ref{e:1.6})
such that each of them cannot be obtained from others by iterations
(resp. iterations or rotations). \textsf{The following is the first
main result of this paper}.

\begin{theorem}\label{th:1.1}
Let $M$ be a $C^3$-smooth compact n-dimensional manifold without
boundary, and $C^2$-smooth map $L:\R\times TM\to\R$ satisfy the
conditions (L1)-(L3). Then:
\begin{description}
\item[(i)]  Suppose that for a homotopy class $\alpha$ of free
loops in $M$ and an abelian group $\K$ the singular homology groups
$H_r(C(\R/\Z, M; \alpha^k); \K)$ have nonzero ranks for some integer
$r\ge n$ and all $k\in\N$. Then \verb"either" for some $l\in\N$
there exist infinitely many distinct $l$-periodic solutions of
(\ref{e:1.6}) representing $\alpha^l$, \verb"or" there exist
infinitely many positive integers $l_1<l_2<\cdots$, such that for
each $i\in\N$ the system (\ref{e:1.6}) has a periodic solution with
minimal period $l_i$ and representing $\alpha^{l_i}$.

\item[(ii)]  Suppose that for some abelian group $\K$ and integer
$r\ge n$ the singular homology groups $H_r(C(\R/\Z, M); \K)$ have
nonzero ranks. Then \verb"either" for some $l\in\N$ there exist
infinitely many distinct $l$-periodic solutions of (\ref{e:1.6}),
\verb"or" there exist infinitely many  positive integers
$l_1<l_2<\cdots$, such that for each $i\in\N$ the system
(\ref{e:1.6}) has a periodic solution with minimal period $l_i$.
 \end{description}
\end{theorem}

Let $0$ denote the free homotopy class of contractible loops in $M$,
i.e., $C(\R/\Z, M; 0)$ consists of all contractible loops
$\gamma:\R/\Z\to M$. The obvious inclusion $\imath:M\to C(\R/\Z, M;
0)$ and the evaluation
$$
{\rm EV}: C(\R/\Z, M; 0)\to M,\;\gamma\mapsto \gamma(0)
$$
satisfy ${\rm EV}\circ\imath=id_M$. It easily follows that
$$
\imath_\ast: H_k(M; \Z_2)\to H_k\bigl( C(\R/\Z, M; 0); \Z_2\bigr)
$$
is injective for any $k\in\N$. Since $H_n(M,\Z_2)=\Z_2$ for $n=\dim
M$, we get
\begin{equation}\label{e:1.8}
{\rm rank}H_n\bigl( C(\R/\Z, M; 0);\Z_2\bigr)\ne 0.
\end{equation}

\begin{corollary}\label{cor:1.2}
Let $M$ be a $C^3$-smooth compact n-dimensional manifold without
boundary, and $C^2$-smooth map $L:\R\times TM\to\R$ satisfy the
conditions (L1)-(L3). Then the system (\ref{e:1.6}) possesses
infinitely many distinct {\rm contractible} integer periodic
solution towers.
\end{corollary}

\begin{remark}\label{rem:1.3}
{\rm  $1^\circ$ \textsf{When $M$ has finite fundamental group},
Benci \cite{Be} first proved that the system (\ref{e:1.6}) has
infinitely many distinct contractible $1$-periodic solutions for
$C^2$-smooth Lagrangian $L$ satisfying the conditions (L1)-(L3) and
$$
\Bigl| \frac{\partial L}{\partial q_i}(t,q,v)\Bigr|\le C(1+
|v|^2),\quad \Bigl| \frac{\partial L}{\partial v_i}L(t,q,v)\Bigr|\le
C(1+ |v|)
$$
in some local coordinates $(q_1,\cdots, q_n)$ for some constant
$C>0$. Recently, under weaker assumptions than (L1)-(L3), i.e.
Tonelli conditions and (L5) below,  Abbondandolo and Figalli
\cite[Cor.3.2]{AbF} showed that the system (\ref{e:1.6}) has an
infinite sequence of $1$-periodic contractible solutions with
diverging action and diverging Morse index. The key in \cite{Be,
AbF} is the fact  that the space of free loops in a compact simply
connected manifold has infinitely many nonzero (co)homology groups
with real coefficients \cite{Su}.  A new technique in \cite{AbF} is
to modify their Tonelli Lagrangian $L$ to one satisfying (L1)-(L3).
\\
$2^\circ$ \textsf{On $n$-dimensional torus $T^n$}, for the
Lagrangian of the form
\begin{equation}\label{e:1.9}
L(t,q,v)=\frac{1}{2}g_q(v,v)+ U(t,q)
\end{equation}
for all $(t, q, v)\in\R\times TT^n=\R\times T^n\times\R^n$, where
$g$ is a $C^3$-smooth Riemannian metric on $T^n$ and $U\in
C^3(\R/\Z\times T^n, \R)$, (such a $L$ satisfies the conditions
(L1)-(L3)),  Yiming Long \cite{Lo2} proved that
  the system (\ref{e:1.6}) possesses
infinitely many distinct contractible integer periodic solution
towers.

We refer the reader to \cite{Lo2} and the references given there for
the detailed history on the integer periodic solutions of the
Lagrangian system. }
\end{remark}

 If $L:\R\times TM\to\R$
   also satisfies
\begin{description}
 \item[(L4)] $L(-t,q, -v)=L(t,q,v)$ for any $(t,q,v)\in\R\times TM$,
\end{description}
we can  improve Corollary~\ref{cor:1.2} as follows.

\begin{theorem}\label{th:1.4}
Let $M$ be a $C^3$-smooth compact n-dimensional manifold without
boundary, and $C^2$-smooth map $L:\R\times TM\to\R$ satisfy the
conditions (L1)-(L4).  Then the system (\ref{e:1.6}) possesses
infinitely many distinct contractible integer periodic solution
towers based on \verb"reversible" periodic solutions.
\end{theorem}

This result was proved by the author and Mingyan Wang \cite{LuW2} in
the case that $M=T^n$ and that $L$ has the form (\ref{e:1.9}) and
satisfies (L4), i.e. $U(-t, q)=U(t, q)$ for any $(t,q)\in\R\times
T^n$.  In particular, we have a generalization of
\cite[Th.1.6]{LuW2}.

\begin{corollary}\label{cor:1.5}
 If $L\in C^2(TM, \R)$ satisfies (L2)-(L4), then for any real number
$\tau>0$,  the following three claims have at least one to be
true:\\
$\bullet$ $L$ has infinitely many critical points sitting in
$M=0_{TM}$ and thus the system (\ref{e:1.6}) possesses infinitely
many different constant
solutions in $M$;\\
$\bullet$ there exists some positive integer $k$ such that the
system (\ref{e:1.6}) possesses infinitely many different nonconstant
$k\tau$-periodic solution orbit towers based on   \verb"reversible" periodic solutions of (\ref{e:1.6});\\
$\bullet$  there exist  infinitely many positive integers
$k_1<k_2<\cdots$, such that for each $k_m$ the system (\ref{e:1.6})
possesses a  \verb"reversible" periodic solution with minimal period
$k_m\tau$, $m=1,2,\cdots$.
\end{corollary}

When $M=T^n$ and $L$ has the form (\ref{e:1.9}) with real analytic
 $g$ and nonconstant, autonomous and real analytic $U$,  the author and
Mingyan Wang \cite{LuW1} observed that suitably improving the
arguments in \cite{CaTa} can give a simple proof of
Corollary~\ref{cor:1.5}. It should also be noted that even if $M$ is
simply connected the methods in \cite{Be, AbF} cannot produce
infinitely many reversible integer periodic solutions because the
space of reversible loops in $M$ can contract to the zero section of
$TM$ and therefore has no infinitely many nonzero Betti numbers.

If $L\in C^2(TM, \R)$ only satisfies (L2)-(L3), it is possible that
two distinct solutions $\gamma_1$ and $\gamma_2$ obtained by
Theorem~\ref{th:1.1} only differ a rotation, i.e.,
$\gamma_1(t)=\gamma_2(s+ t)$ for some $s\in\R$ and any $t\in\R$.
However, we can combine the proof of Theorem~\ref{th:1.1} with the
method in \cite{LoLu} to improve the results in Theorem~\ref{th:1.1}
as follows:

\begin{theorem}\label{th:1.6}
Let $M$ be a $C^5$-smooth compact n-dimensional manifold without
boundary, and $C^4$-smooth map $L:TM\to\R$ satisfy the conditions
(L2)-(L3). Then for any $\tau>0$ the following results hold:
\begin{description}
\item[(i)]  Suppose that for a homotopy class $\alpha$ of free
loops in $M$ and an abelian group $\K$ the singular homology groups
$H_r(C(\R/\Z, M; \alpha^k); \K)$ have nonzero ranks for some integer
$r\ge n$ and all $k\in\N$. If either $r\ge n+ 1$ or $r=n>1$, then
\verb"either" for some $l\in\N$ there exist infinitely many distinct
 periodic solution orbit towers based on $l\tau$-periodic solutions of
 (\ref{e:1.6}) representing $\alpha^l$,
\verb"or" there exist infinitely many positive integers
$l_1<l_2<\cdots$, such that for each $i\in\N$ the system
(\ref{e:1.6}) has a periodic solution orbit tower based on a
periodic solution with minimal period $l_i\tau$ and representing
$\alpha^{l_i}$.

\item[(ii)]  Suppose that the singular homology groups $H_r(C(\R/\Z, M); \K)$
have nonzero ranks for some integer $r\ge n$ and some abelian group
$\K$. If either $r\ge n+ 1$ or $r=n>1$, then \verb"either" for some
$l\in\N$ there exist infinitely many distinct periodic solution
orbit towers based on $l\tau$-periodic solutions of (\ref{e:1.6}),
\verb"or" there exist infinitely many positive integers
$l_1<l_2<\cdots$, such that for each $i\in\N$ the system
(\ref{e:1.6}) has a periodic solution orbit tower based on a
periodic solution with minimal period $l_i\tau$.
 \end{description}
\end{theorem}

By (\ref{e:1.8}) we immediately get:

\begin{corollary}\label{cor:1.7}
Let $M$ be a $C^5$-smooth compact  manifold of dimension $n>1$ and
without boundary, and $C^4$-smooth map $L:TM\to\R$ satisfy the
conditions (L2)-(L3). Then for any $\tau>0$  the system
(\ref{e:1.6}) possesses infinitely many distinct periodic solution
orbit towers based on \verb"contractible" periodic solutions of
integer multiple periods of $\tau$.
\end{corollary}

Clearly, when (L4) is satisfied Corollary~\ref{cor:1.5} seems to be
stronger than Corollary~\ref{cor:1.7}. If $n=1$ and (L4) is not
satisfied, we do not know whether Corollary~\ref{cor:1.7} is still
true. Moreover, the reason that we require higher smoothness in
Theorem~\ref{th:1.6} and Corollary~\ref{cor:1.7} is to assure that
the normal bundle of a nonconstant periodic orbit is $C^2$-smooth.

When $M=T^n$ and $L$ has the form (\ref{e:1.9}) with \textsf{flat
$g$ and autonomous $U$}, Yiming Long and the author \cite{LoLu}
developed the equivariant version of the arguments in \cite{Lo2} to
prove Corollary~\ref{cor:1.7}. Even if $g$ is not flat, the author
and Mingyan Wang \cite[Th.1.6]{LuW2} also derived a stronger result
than  Corollary~\ref{cor:1.7} in the case that $M=T^n$.
 Campos and Tarallo \cite{CaTa} obtained a
similar result provided
 that the metric $g$ is real analytic, and that the  potential $U$ is autonomous, real
 analytic and nonconstant.

Even if $L=\frac{1}{2}g$ for a $C^4$-Riemannian metric $g$ on $M$,
it seems that Theorem~\ref{th:1.6} or Corollary~\ref{cor:1.7} cannot
yield infinitely many geometrically distinct closed geodesics.

Assume that $L$ also satisfies
\begin{description}
 \item[(L5)] For any $(q,v)\in TM$ there exists
an unique solution of (\ref{e:1.6}), $\gamma:\R\to M$, such that
$(\gamma(0), \dot\gamma(0))=(q,v)$.
\end{description}
By \cite[\S2]{AbF}, this assumption can be satisfied if
\begin{equation}\label{e:1.11}
-\partial_t L(t, q, v)\le c\left(1+
D_vL(t,q,v)[v]-L(t,q,v)\right)\quad\forall (t, q, v)\in\R\times TM.
\end{equation}
(Clearly, the left side may be replaced by $const- \partial_t L(t,
q, v)$ since (L5) is also satisfied up to adding a constant to $L$.
Moreover, that $L$ satisfies (L1)-(L3) is equivalent to that the
Fenchel transform $H$ of $L$ given by (\ref{e:1.4}) satisfies the
assumptions (H1)-(H3) below. In this case (\ref{e:1.11}) is
equivalent to (\ref{e:1.13}) below. Hence (\ref{e:1.11}) holds if
$L$ is independent of $t$ as noted below (\ref{e:1.13}).)  Under the
assumption (L5), we have an one-parameter family of
$C^1$-diffeomorphisms $\Phi_L^t\in {\rm Diff}(TM)$ satisfying
$\Phi_L^t(\gamma(0), \dot\gamma(0))=(\gamma(t), \dot\gamma(t))$.
(See \cite[Th.2.6.5]{Fa}). Following \cite{Lo2}, the time-1-map
$\Phi_L=\Phi^1_L$ is called the {\bf Poincar\'{e} map} of the system
(\ref{e:1.6}) corresponding to the Lagrangian function $L$. Every
integer periodic solution $\gamma$ of (\ref{e:1.6}) gives a periodic
point $(\gamma(0), \dot\gamma(0))$ of $\Phi_L$. If $\gamma$ is even,
then the periodic point $(\gamma(0), \dot\gamma(0))$ sits in the
zero section $0_{TM}$ of $TM$. So Corollary~\ref{cor:1.2} and
Theorem~\ref{th:1.4} yield the following

\begin{corollary}\label{cor:1.8}
Let $M$ be a $C^3$-smooth compact n-dimensional manifold without
boundary, and $C^2$-smooth map $L:\R\times TM\to\R$ satisfy the
conditions (L1)-(L3) and (L5). Then the Poincar\'{e} map $\Phi_L$
has infinitely many distinct periodic points. Furthermore, if (L4)
is also satisfied then the Poincar\'{e} map $\Phi_L$ has infinitely
many distinct periodic points sitting in the zero section $0_{TM}$
of $TM$.
\end{corollary}

If $L$ is independent of $t$, for a periodic point $(\gamma(0),
\dot{\gamma}(0))$ of $\Phi_L$ generated by a $\tau$-periodic
solution $\gamma$, then all points of $\{(\gamma(s),
\dot{\gamma}(s))\,|\, s\in\R\}$ are periodic points of $\Phi_L$. We
call such period points {\bf orbitally same}. By remarks below
(\ref{e:1.11}), using Corollary~\ref{cor:1.7} we can improve
Corollary~\ref{cor:1.8} as follows:

\begin{corollary}\label{cor:1.9}
Let $M$ be a $C^5$-smooth compact  manifold of dimension $n>1$ and
without boundary, and $C^4$-smooth map $L:TM\to\R$ satisfy the
conditions (L2)-(L3). Then the Poincar\'{e} map $\Phi_L$ has
infinitely many \verb"orbitally distinct" periodic points.
\end{corollary}

It is easily checked that the assumption (L4) is equivalent to the
following:
\begin{description}
\item[(H4)] $H(-t,q, -p)=H(t,q, p)$ for any $(t,q,p)\in\R\times T^\ast M$.
 \end{description}
In this case,  $v=v(t,q,p)$ uniquely determined by the equality
$p=D_vL(t,q,v)$ satisfies
\begin{equation}\label{e:1.12}
v(-t,q,-p)=-v(t,q,p)\quad\forall (t,q,p)\in\R\times T^\ast M.
\end{equation}
So if a solution $\gamma:\R\to M$ of (\ref{e:1.6}) satisfies
 $\gamma(-t)=\gamma(t)\;\forall t\in\R$, then
$\gamma^\ast(-t)=-\gamma^\ast(t)$ for all $t\in\R$.

With the same way as the definition of solution towers and solution
orbit towers to (\ref{e:1.6}) we can define
  solution towers to (\ref{e:1.1}), and solution orbit towers to (\ref{e:1.1}) in the case
  $H$ is independent of $t$. Then the Hamiltonian versions from
  Theorem~\ref{th:1.1} to Corollary~\ref{cor:1.7} can be obtained
  immediately. For example, from Corollary~\ref{cor:1.2},
  Theorem~\ref{th:1.4} and Corollary~\ref{cor:1.7} we directly
  derive:

\begin{theorem}\label{th:1.10}
$1^\circ)$ Let $M$ be a $C^3$-smooth compact n-dimensional manifold
without boundary, and $C^2$-smooth map $H:\R\times T^\ast M\to\R$
satisfy the conditions (H1)-(H3).  Then the system (\ref{e:1.1})
possesses infinitely many distinct {\rm contractible} integer
periodic solution towers. Furthermore, if (H4)  is also satisfied
then  the system (\ref{e:1.1}) possesses infinitely many distinct
contractible integer periodic solution towers based on
 periodic solutions with \verb"reversible" projections to $M$.\\
$2^\circ)$ Let $M$ be a $C^5$-smooth compact  manifold of dimension
$n>1$ and without boundary, and $C^4$-smooth map $H:\R\times T^\ast
M\to\R$ satisfy the conditions (H2)-(H3). Then for any $\tau>0$  the
system (\ref{e:1.1}) has infinitely many distinct periodic solution
orbit towers based on \verb"contractible" periodic solutions of
integer multiple periods of $\tau$.
\end{theorem}

\begin{remark}\label{rm:1.11}
{\rm If $\pi_1(M)$ is finite, Cieliebak \cite{Ci} showed that the
system (\ref{e:1.1}) has infinitely many contractible 1-periodic
solutions (with unbounded actions) provided that $H\in
C^\infty(\R/\Z\times T^\ast M, \R)$ satisfies
\begin{description}
\item[(HC1)] $dH(t,q,p)\left[p \frac{\partial}{\partial p}\right] - H(t,q,p)\geq
h_0 \|p\|^2- h_1$,
\item[(HC2)] $\left|\frac{\partial^2 H}{\partial p_i\partial
p_j}(t,q,p)\right|\le d\quad\hbox{and}\quad \left|\frac{\partial^2
H}{\partial p_i\partial q_j}(t,q,p)\right|\le d$,
\end{description}
for all $(t,q,p)\in\R\times T^\ast M$, with respect to a suitable
metric on the bundle $T^\ast M\to M$ and constants $h_0>0, h_1$ and
$d$.  Here $q_1,\cdots, q_n, p_1,\cdots,p_n$ are coordinates on
$T^\ast M$ induced by geodesic normal coordinates $q_1,\cdots, q_n$
on $M$.\\
 Recently, Abbondandolo and Figalli  stated in \cite[Remark 7.4]{AbF} that the same
result can be derived from \cite[Th.7.3]{AbF} if the
 assumptions (HC1)-(HC2) are replaced by
\begin{description}
\item[(HAF1)] $dH(t,q,p)\left[p \frac{\partial}{\partial p}\right] - H(t,q,p)\geq
a(|p|_q)$  for some function $a:[0, \infty)\to\R$ with $\lim_{s\to
+\infty}a(s)=+\infty$,
\item[(HAF2)] $H(t, q,p)\ge h(|p|_q)$ for some function $h:[0, \infty)\to\R$
with $\lim_{s\to +\infty}\frac{h(s)}{s}=+\infty$ and all $(t,
q,p)\in\R\times T^\ast M$,
\end{description}
and (H5) below. Note that no convexity assumption on $H$ was made in
\cite{Ci, AbF} and therefore that their results cannot be obtained
from one on Lagrangian system via the Legendre transform. }
\end{remark}

It is easily seen that the assumption  (L5) is equivalent to the
following:
\begin{description}
 \item[(H5)] For any $(q, p)\in T^\ast M$ there exists
an unique solution of $\dot{x}(t)=X_H(t,x(t))$, $x:\R\to M$, such
that $x(0)=(q, p)$.
\end{description}
The assumption  can be satisfied under the following equivalent
condition of (\ref{e:1.11}):
\begin{equation}\label{e:1.13}
\partial_t H(t, q, p)\le c\left(1+
H(t,q, p)\right)\quad\forall (t, q, p)\in\R\times T^\ast M,
\end{equation}
see \cite[pp.629]{AbF}. Since (H2) implies that $H$ is superlinear
on the fibers of $T^\ast M$, (\ref{e:1.13}) holds clearly if $H$ is
independent of time $t$. The condition (H5) guarantees that the
global flow of $X_H$ exists on $T^\ast M$. Thus we have an
one-parameter family of Hamiltonian diffeomorphisms $\Psi^H_t\in
{\rm Ham}(T^\ast M, \omega_{\rm can})$ satisfying
$\Psi^H_t(\gamma(0), \dot\gamma^\ast(0))=(\gamma(t), \dot\gamma^\ast
(t))$. As usual, the time-1-map $\Psi^H=\Psi_1^H$ is called the {\bf
Poincar\'{e} map} of the system (\ref{e:1.1}) corresponding to the
Hamiltonian function $H$. For each $t\in\R$ recall that the Legendre
transform associated with $L_t(\cdot)=L(t, \cdot)$ is given by
$$
\frak{L}_{L_t}: TM\to  T^\ast M, \quad (q, v)\mapsto \left(q,
D_vL(t,q, v)\right).
$$
It is easy to check that
\begin{equation}\label{e:1.14}
\Psi^H_t\circ\frak{L}_{L_0}=\frak{L}_{L_t}\circ\Phi^t_L\quad\hbox{for
any}\; t\in\R.
\end{equation}
From this one immediately gets the following equivalent Hamiltonian
versions of Corollary~\ref{cor:1.8} and Corollary~\ref{cor:1.9}.

\begin{theorem}\label{th:1.12}
$1^\circ)$ Let $M$ be a $C^3$-smooth compact n-dimensional manifold
without boundary, and $C^2$-smooth map $H:\R\times T^\ast M\to\R$
satisfy the conditions (H1)-(H3) and (H5).  Then the Poincar\'{e}
map $\Psi^H$ has infinitely many distinct periodic points.
Furthermore, if (H4)  is also satisfied then the Poincar\'{e} map
$\Psi^H$ has infinitely many distinct periodic
points sitting in the zero section $0_{T^\ast M}$ of $T^\ast M$.\\
$2^\circ)$ Let $M$ be a $C^5$-smooth compact  manifold of dimension
$n>1$ and without boundary, and $C^4$-smooth map $H:\R\times T^\ast
M\to\R$ satisfy the conditions (H2)-(H3). Then the Poincar\'{e} map
$\Psi^H$ has infinitely many orbitally distinct periodic points.
(That is, any two do not sit the same Hamiltonian orbit.)
\end{theorem}

Theorems~\ref{th:1.10},~\ref{th:1.12}  may be viewed a solution for
the Conley conjecture  for Hamiltonian systems on cotangent bundles,
and Corollary~\ref{cor:1.8} and Corollary~\ref{cor:1.9} may be
viewed as confirm answers of Lagrangian systems analogue of the
Conley conjecture for Hamiltonian systems.

The main proof ideas come from \cite{Lo2}. We shall prove
Theorems~\ref{th:1.1}, \ref{th:1.6} in the case $r=n$, and
Theorem~\ref{th:1.4}  by generalizing the variational arguments in
\cite{Lo2}, \cite{LoLu} and \cite{LuW2} respectively. Some new ideas
are needed because we do not lift to the universal cover space of
$M$ as done in \cite{Lo2, LoLu, LuW2} for the tori case. We also
avoid using finite energy homologies used in \cite{Lo2, LoLu, LuW2}.
Let us outline the variational setup and new ideas as follows.  For
$\tau>0$, let
$$
S_\tau:=\R/\tau\Z=\{[s]_{\tau}\,|\,[s]_{\tau}=s+\tau\Z,\, s\in\R\},
\quad{\rm and}\quad E_{\tau}=W^{1,2}(S_{\tau}, M)
$$
denote the space of all loops $\gamma:S_\tau\to M$ of Sobolev class
$W^{1,2}$. For a homotopy class $\alpha$ of free loops in $M$, let
$$
H_\tau(\alpha),\quad H_\tau=H_\tau(0),\quad EH_\tau
$$
respectively denote the subset of loops of $E_\tau$ representing
$\alpha$, that of all contractible loops in $E_\tau$, and that of
all reversible loops in $E_\tau$. Then $EH_\tau\subset H_\tau$.

For integer $m\ge 2$,  if $M$ is $C^m$-smooth, all these spaces
$E_\tau$, $H_\tau(\alpha)$ and $EH_\tau$ have  $C^{m-1}$-smooth
Hilbert manifold structure \cite{Kl}, and the tangent space of
$E_\tau$ at $\gamma$ is $T_\gamma E_\tau=W^{1,2}(\gamma^\ast TM)$.
Moreover, any ($C^{m-1}$) Riemannian metric
$\langle\cdot,\cdot\rangle$ on $M$  induces a complete Riemannian
metric  on $E_\tau$:
\begin{eqnarray}\label{e:1.15}
&&\langle\!\langle \xi , \eta \rangle\!\rangle_\tau = \int_0^\tau
\left( \langle \xi(t), \eta(t)\rangle_{\gamma(t)}+ \langle\nabla_t
\xi(t),
  \nabla_t \eta(t)\rangle_{\gamma(t)}\right)\, dt\\
&&\qquad \forall\gamma\in E_\tau,\; \xi,\eta\in T_\gamma
  E_\tau=W^{1,2}(\gamma^\ast TM).\nonumber
\end{eqnarray}
 Here $\nabla_t$ denotes the covariant derivative in direction $\dot{\gamma}$ with respect
 to the Levi-Civita connection $\nabla$ of $\langle\cdot,\cdot\rangle$. Let
$\|\xi\|_\tau=\sqrt{\langle\!\langle \xi ,
\xi\rangle\!\rangle_\tau}\quad\forall\xi\in T_\gamma E_\tau$. Then
the distance on $E_\tau$ induced by $\|\cdot\|_\tau$ is {\it
complete} and also {\it compatible} with the manifold topology on
$E_\tau$. Consider the functional ${\cal L}_\tau:E_\tau\to\R$,
\begin{equation}\label{e:1.16}
{\cal L}_{\tau}(\gamma)=\int^{\tau}_0L(t, \gamma(t),\dot
\gamma(t))dt \quad\forall \gamma\in E_{\tau}.
\end{equation}
For integer $m\ge 3$,  if $M$ is $C^m$-smooth and $C^{m-1}$-smooth
$L:\R\times TM\to\R$ satisfies the assumptions (L1)-(L3), then the
functional ${\cal L}_\tau$ is $C^2$-smooth, bounded below, satisfies
the Palais-Smale condition, and all critical points of it have
finite Morse indexes and nullities, (see \cite[Prop.4.1, 4.2]{AbF}
and \cite{Be}). By \cite[Th.3.7.2]{Fa}, all critical points of
${\cal L}_\tau$ are all of class $C^{m-1}$ and therefore correspond
to all $\tau$-periodic solutions of (\ref{e:1.6}).

 Let ${\cal L}_\tau^E$
denote the restriction of ${\cal L}_\tau$ on $EH_{\tau}$. When $L$
satisfies (L4), it is not hard to prove that a map $\gamma:\R\to M$
is a $\tau$-periodic even solution to (\ref{e:1.6}) if and only if
$\gamma$ is a critical point of ${\cal L}^E_\tau$ on $EH_{\tau}$,
cf. \cite[Lem.1.7]{LuW2}.

When we attempt to prove Theorem~\ref{th:1.1} by the method of
\cite{Lo2}, we first need to know how to relate the Morse index and
nullity of a critical point $\gamma\in E_\tau$ of ${\cal L}_\tau$ to
those of the $k$-th iteration $\gamma^k\in E_{k\tau}$ as a critical
point of ${\cal L}_{k\tau}$ on $E_{k\tau}$. Since we do not assume
that $M$ is orientable or $\gamma$ is contractible, the bundle
$\gamma^\ast TM\to S_\tau$ might not be trivial. However, for the
$2$-th iteration $\gamma^2$, the pullback bundle $(\gamma^2)^\ast
TM\to S_{2\tau}$ is always trivial. Since our proof is indirect by
assuming that the conclusion does not hold, the arguments can be
reduced to the case that all $\tau$-periodic solutions have trivial
pullback bundles (as above Lemma~\ref{lem:5.2}).  For such periodic
solutions we can choose suitably coordinate charts around them on
$E_{k\tau}$ so that the question is reduced to the case $M=\R^n$ as
in Lemma~\ref{lem:3.2}. Hence we can get expected iteration
inequalities as in Theorem~\ref{th:3.1}. The \textsf{second new
idea} is that under the assumption each ${\cal L}_{k\tau}$ has only
isolated critical points we show in Lemma~\ref{lem:5.2} how to use
an elementary arguments as above Corollary~\ref{cor:1.2} and the
Morse theory to get a non-minimal saddle point with nonzero th-$n$
critical module with $\Z_2$-coefficient;  the original method  in
\cite[Lemma 4.1]{Lo2} is to use Lemma II.5.2 on the page 127 of
\cite{Ch} to arrive at this goal, which seems to be difficult for me
generalizing it to manifolds. It is \textsf{worth noting that we
avoid using finite energy homologies} used in \cite{Lo2, LoLu,
LuW2}. That is based on an observation, that is, the composition
$(j_{k\tau})_*\circ \psi^k_*$ in (\ref{e:5.13}) has a good
decomposition $({\bf J}_k)_\ast\circ(\psi^k)_\ast\circ ({\bf
I}_1)_\ast$ as in (\ref{e:5.15}) such that for each $\omega\in
C_n({\cal L}_{\tau}, \gamma;\K)$, $({\bf I}_1)_\ast(\omega)$ is a
singular homology class of a $C^1$-Hilbert manifold and hence has a
$C^1$-singular cycle representative. It is the final claim that
allows us to use the singular homology to complete the remained
arguments in Long'method of \cite{Lo2}. A merit of this improvement
is to reduce the smoothness of the Largangian $L$. That is, we only
need to assume that $L$ is of class $C^2$. However, a new problem
occurs, i.e. $\tilde\Theta_{k\tau}$ in (\ref{e:4.12}) is only a
homeomorphism. It is very fortunate that $\tilde\alpha_{k\tau}$ is
also of class $C^2$ as noted at the end of proof of Theorem~5.1 (the
generalized Morse lemma) on the page 44 of \cite{Ch}. Using the
image of Gromoll-Meyer of $\tilde\alpha_{k\tau}(\eta)+
\tilde\beta_{k\tau}(\xi)$ under $\tilde\Theta_{k\tau}$, called
topological Gromoll-Meyer, to replace a Gromoll-Meyer of
$\tilde{\cal L}_{k\tau}$ at $\tilde\gamma^k$, we construct
topological Gromoll-Meyer pairs of ${\cal L}_\tau$ at $\gamma\in
H_\tau(\alpha)$ and
 of ${\cal L}_{k\tau}$ at $\gamma^k\in H_{k\tau}(\alpha^k)$, to
 satisfy Theorem~\ref{th:4.4} which is enough to complete our proof
 of Theorem~\ref{th:1.1}. For the proof of Theorem~\ref{th:1.6} we
 need to complete more complex arguments as in \S4.3. But the ideas
 are similar.

The paper is organized as follows. Section 2 will review some basic
facts concerning the Maslov-type indices and relations between them
and Morse indexes. In Section 3 we give some iteration inequalities
of the Morse indexes. Section 4 studies changes of the critical
modules under iteration maps. In Sections 5, 6 and 7, we give the
proofs of Theorems~\ref{th:1.1}, ~\ref{th:1.4} and ~\ref{th:1.6}
respectively. Motivated by the second claim in Theorem~\ref{th:1.10}
$1^\circ$), a more general question than the Conley's conjecture and
a program in progress are proposed in Section 8. In Appendix of
Section 9 we prove Proposition A and a key Lemma A.4, which is a
generalization of \cite[Lemma 2.3]{Lo2}.\vspace{2mm}

{\bf Acknowledgements}: I am greatly indebted to Professor Yiming
Long for leading me this to question ten years ago. The author
sincerely thanks Professors  Le Calvez and C. Viterbo for organizing
a seminar of symplectic dynamics at Beijing International
Mathematics Center in May 2007, where my interest for this question
was aroused again. He also sincerely thanks Professor Alberto
Abbondandolo for some helps in understanding his paper. The results
and outlines of proofs in this paper were reported in the workshop
on Floer Theory and Symplectic Dynamics at CRM of University of
Montreal, May 19-23, 2008.  I would like to thank the organizers for
their invitation, and CRM for hospitality. Finally, I sincerely
thank Professor Kung-Ching Chang for his helps in correcting
mistakes in the first draft.

\section{Maslov-type indices and Morse index}
\setcounter{equation}{0}

\noindent{\bf 2.1. A review on Maslov-type indices}.\; Let ${\rm
Sp}(2n,\R)=\{M\in\R^{2n\times 2n}\,|\, M^TJ_0M=J_0\}$, where
$J_0={\scriptscriptstyle\left(\begin{array}{cc}
0& -I_n\\
I_n& 0\end{array}\right)}$. For $\tau>0$, denoted by
\begin{eqnarray*}
&&{\cal P}_\tau(2n)=\{\Psi\in C([0,\tau], {\rm
Sp}(2n,\R))\,|\,\Psi(0)=I_{2n}\},\\
&&{\cal P}^\ast_\tau(2n)=\{\Psi\in {\cal
P}_\tau(2n)\,|\,\det(\Psi(\tau)-I_{2n})\ne 0\}.
\end{eqnarray*}
The paths in ${\cal P}^\ast_\tau(2n)$ are called {\bf
nondegenerate}. The Maslov-type index (or Conley-Zehnder index)
theory for the paths in ${\cal P}^\ast_\tau(2n)$ was defined by
\cite{CoZe}, \cite{Lo1} and \cite{Vi2}. Yiming Long \cite{Lo4}
extended this theory to all paths in ${\cal P}_\tau(2n)$. The
Maslov-type index of a path $\Psi\in{\cal P}_\tau(2n)$ is a pair of
integers $({i}_\tau(\Psi),\nu_{\tau}(\Psi))$, where
\begin{eqnarray*}
&&\nu_{\tau}(\Psi)=\dim_{\R}{\rm
Ker}_{\R}(\Psi(\tau)-I_{2n})\quad{\rm
and}\\
&& {i}_\tau(\Psi)=\inf\{i_\tau(\beta)\,|\, \beta\in{\cal
P}^\ast_\tau(2n)\;\hbox{is sufficiently $C^0$ close to $\Psi$
in}\;{\cal P}_\tau(2n)\}
\end{eqnarray*}
with $i_\tau(\beta)$ defined as in \cite{CoZe}.
 Clearly, the map $i_\tau:{\cal P}_\tau(2n)\to\Z$ is lower
semi-continuous. For  any paths $\Psi_k\in{\cal P}_\tau(2n)$,
$k=0,1$,
$(i_\tau(\Psi_0),\nu_{\tau}(\Psi_0))=(i_\tau(\Psi_1),\nu_{\tau}(\Psi_1))$
if and only if there exists a homotopy $\Psi_s, 0\le s\le 1$ from
$\Psi_0$ to $\Psi_1$ in ${\cal P}_\tau(2n)$ such that
$\Psi_s(0)=I_{2n}$ and
$\nu_{\tau}(\Psi_s(\tau))\equiv\nu_\tau(\Psi_0)$ for any $s\in
[0,1]$.

 For $a<b$ and any path $\Psi\in C([a,b],{\rm Sp}(2n,\R))$, choose
$\beta\in{\cal P}_1(2n)$ with $\beta(1)=\Psi(a)$, and define
$\phi\in{\cal P}_1(2n)$ by $\phi(t)=\beta(2t)$ for $0\le t\le 1/2$,
and
$$
\phi(t)=\Psi(a+(2t-1)(b-a))\quad{\rm for}\quad 1/2\le t\le 1.
$$
 It was showed
in \cite{Lo4} that the difference $i_1(\phi)-i_1(\beta)$ only
depends on $\Psi$, and was called the {\bf Maslov-type index} of
$\Psi$ , denoted by
\begin{equation}\label{e:2.1}
i(\Psi, [a,b]):=i_1(\phi)-i_1(\beta).
\end{equation}
Clearly, $i(\Psi, [0,1])=i_1(\Psi)$ for any $\Psi\in{\cal P}_1(2n)$.

 Let
$(F,\{\cdot,\cdot\})$ be the symplectic space with
$F=\R^{2n}\oplus\R^{2n}$ and
$$
\{u,v\}=\langle{\cal J}u, v\rangle\quad\forall u,v\in F,\quad{\rm
where}\;{\cal J}={\scriptscriptstyle\left(\begin{array}{cc}
-J_0& 0\\
0& J_0\end{array}\right)}.
$$
{\bf All vectors are understand as column vectors in this paper
without special statements.} Let ${\rm Lag}(F)$ be the manifold of
Lagrangian Grassmannian of $(F,\{\cdot,\cdot\})$, and $\mu^{\rm
CLM}$ be the Cappell-Lee-Miller index characterized by properties
I-VI of \cite[pp. 127-128]{CLM}. There  exists the following
relation between $\mu^{\rm CLM}$ and the index defined by
(\ref{e:2.1}),
\begin{equation}\label{e:2.2}
i(\Psi, [a,b])=\mu^{\rm CLM}_F(W, {\rm Gr}\left(\Psi),
[a,b]\right)-n,
\end{equation}
where $W=\{(x^T,x^T)^T\in\R^{4n}\,|\,x\in\R^{2n}\}$.

 With $U_1=\{0\}\times\R^n$ and $U_2=\R^n\times\{0\}$, two new
Maslov-type indices for any path $\Psi\in C([a,b],{\rm Sp}(2n,\R))$
were defined in \cite{LoZZ} as follows:
\begin{equation}\label{e:2.3}
\mu_k(\Psi,[a,b])=\mu^{\rm CLM}_{\R^{2n}}(U_k,\Psi U_k,
[a,b]),\;k=1,2.
\end{equation}
 Let
$\Psi(b)={\scriptscriptstyle \left(\begin{array}{cc}
A& B\\
C& D\end{array}\right)}$, where $A, B, C, D\in\R^{n\times n}$. In
terms of \cite[(2.21)]{LoZZ}, define
\begin{equation}\label{e:2.4}
\nu_1(\Psi, [a,b])=\dim{\rm Ker}(B)\quad{\rm and}\quad \nu_2(\Psi,
[a,b])=\dim{\rm Ker}(C).
\end{equation}
In particular, for $\Psi\in{\cal P}_\tau(2n)$ and $k=1,2$ we denote
by
\begin{equation}\label{e:2.5}
\mu_{k, \tau}(\Psi)=\mu_k(\Psi, [0, \frac{\tau}{2}])\quad{\rm
and}\quad \nu_{k, \tau}(\Psi)=\nu_k(\Psi, [0, \frac{\tau}{2}]).
\end{equation}

\noindent{\bf Assumption B}. (B1) Let $B\in C(\R, \R^{2n\times 2n})$
be a path of symmetric matrix which is $\tau$-periodic in time $t$,
i.e., $B(t+\tau)=B(t)$ for any $t\in\R$.\\
(B2) Let $B(t)=\left(\begin{array}{cc}B_{11}(t) & B_{12}(t)\\
B_{21}(t)& B_{22}(t)\end{array}\right)$, where $B_{11},
B_{22},\;t\mapsto\R^{n\times n}$ are \verb"even" at $t=0$ and
$\tau/2$, and $B_{12}, B_{21},\;t\mapsto\R^{n\times n}$ are
\verb"odd" at $t=0$ and $\tau/2$. \vspace{2mm}

Under the assumption (B1), let $\Psi$ be  the fundamental solution
of the problem
\begin{equation}\label{e:2.6}
\dot{\Psi}(t)=J_0B(t)\Psi(t),\quad \Psi(0)=I_{2n}.
\end{equation}
By the classical Floquet theory, $\nu_\tau(\Psi)$ is the dimension
of the solution space of  the linear Hamiltonian system
$$
\dot{\bf u}(t)=J_0B(t){\bf u}(t)\quad{\rm and}\quad {\bf
u}(t+\tau)={\bf u}(t).
$$
Similarly, under the assumptions (B1) and (B2),
 it was also shown in
\cite[Prop.1.3]{LoZZ}) that  $\nu_{1,\tau}(\Psi)$ and
$\nu_{2,\tau}(\Psi)$ are the dimensions of the solution spaces of
the following two problems
 respectively,
\begin{eqnarray*}
\left\{\begin{array}{ll}
 &\dot{\bf u}(t)=J_0B(t){\bf u}(t),\\
&{\bf u}(t+\tau)={\bf u}(t),\quad {\bf u}(-t)=N{\bf u}(t),
\end{array}\right.\\
\left\{\begin{array}{ll}
 &\dot{\bf u}(t)=J_0B(t){\bf u}(t),\\
&{\bf u}(t+\tau)={\bf u}(t),\quad {\bf u}(-t)=-N{\bf u}(t),
\end{array}\right.
\end{eqnarray*}
where $N=\left(\begin{array}{cc}-I_n & 0\\0& I_n\end{array}\right)$.
Let $(x_1,\cdots, x_n, y_1,\cdots, y_n)$ denote the coordinates in
$\R^{2n}=\R^n\times\R^n$. Denote by $\omega_0=\sum^n_{k=1}dx_k\wedge
dy_k$ the standard symplectic structure  on $\R^{2n}$, i.e.
$\omega_0({\bf u}, {\bf v})=\langle J_0{\bf u}, {\bf
v}\rangle\;\forall {\bf u}, {\bf v}\in\R^{2n}$. Here
$\langle\cdot,\cdot\rangle$ is the standard inner product on
$\R^{2n}$. Define $H:\R\times\R^{2n}\to\R$ by $H(t, {\bf
u})=\frac{1}{2}\langle B(t){\bf u}, {\bf u}\rangle$. Let $X_H$ be
the corresponding Hamiltonian vector field defined by
\begin{equation}\label{e:2.7}
\omega_0(X_H(t,{\bf u}), {\bf v})=-d_{\bf u}H(t,{\bf u})({\bf v}).
\end{equation}
Then $X_H(t,{\bf u})=J_0B(t){\bf u}$ for any ${\bf u}\in\R^{2n}$.

For $\Psi\in{\cal P}_\tau(2n)$, extend the definition of $\Psi$ to
$[0, +\infty)$ by
\begin{equation}\label{e:2.8}
\Psi(t)=\Psi(t-j\tau)\Psi(\tau)^j,\quad\forall j\tau\le t\le (j+
1)\tau,\;j\in\N,
\end{equation}
and define the $m$-th iteration $\Psi^m$ of $\Psi$ by
\begin{equation}\label{e:2.9}
\Psi^m=\Psi|_{[0, m\tau]}.
\end{equation}
It was proved in \cite[pp. 177-178]{Lo3} that the mean index per
$\tau$ of $\Psi\in{\cal P}_\tau(2n)$,
\begin{equation}\label{e:2.10}
\hat i_\tau(\Psi):=\lim_{m\to +\infty}\frac{i_{m\tau}(\Psi^m)}{m}
\end{equation}
always exists.

\begin{lemma}\label{lem:2.1}
{\rm (i)} For any $\Psi\in{\cal P}_\tau(2n)$ it holds that
$$
{\rm max}\left\{0,m\hat i_{\tau}(\Psi)-n\right\}\le
i_{m\tau}(\Psi^m)\le m\hat i_{\tau}(\Psi)+n- \nu_{m\tau}(\Psi^m),
\quad\forall m\in\N.
$$
{\rm (ii)} $|\mu_1(\Psi)-\mu_2(\Psi)|\le n$ for any
$\Psi\in{\cal P}_\tau(2n)$ with $\tau>0$.\\
{\rm (iii)} Under Assumption B, let $\Psi:[0, +\infty)\to{\rm
Sp}(2n,\R)$ be the fundamental solution of the problem
(\ref{e:2.6}). (It must satisfy (\ref{e:2.8})). Then
\begin{equation}\label{e:2.11}
\mu_{1, m\tau}\left(\Psi|_{[0, \frac{m\tau}{2}]}\right) +  \mu_{2,
m\tau}\left(\Psi|_{[0, \frac{m\tau}{2}]}\right)=
i_{m\tau}\left(\Psi|_{[0, m\tau]}\right) + n\quad\forall m\in\N,
\end{equation}
(or equivaliently $\mu_1(\Psi, [0, m\tau]) + \mu_2(\Psi, [0,
m\tau])= i_{m\tau}\left(\Psi|_{[0, m\tau]}\right) + n\;\forall
m\in\N$).
 Moreover, for $k=1,2$ the mean indices of
$\Psi$ per $\tau$ defined by
\begin{equation}\label{e:2.12}
\hat \mu_{k,\tau}(\Psi):=\lim_{m\to +\infty}\frac{\mu_{k,
m\tau}\left(\Psi|_{[0, m\tau]}\right)}{m}
\end{equation}
always exist and equal to $\frac{1}{2}\hat i_\tau(\Psi)$.
\end{lemma}

(i) comes from \cite{LiLo} or \cite[p. 213, (17)]{Lo3}, (ii) is
\cite[Th.3.3]{LoZZ}, and (iii) is \cite[Prop.C, Cor.6.2]{LoZZ}
(precisely is derived from the proof of \cite[Prop.C,
Cor.6.2]{LoZZ}). It is easily checked that (i) implies
$|i_{m\tau}-mi_\tau|\le (m+1)n$ for any $m\in\N$. A similar
inequality to the latter was also derived in \cite[(12)]{DDP}
recently.\vspace{2mm}

\noindent{\bf 2.2. Relations between Maslov-type indices and Morse
indices}.\quad
\begin{lemma}\label{lem:2.2}\!\!{\rm (\cite{Vi1, LoAn})}.\quad
 Let the Lagrangian $L:\R\times\R^{2n}\to\R$ be given by
$$
 L(t, y, v)=\frac{1}{2}P(t)v\cdot v + Q(t)y\cdot v +
\frac{1}{2}R(t)y\cdot y,
$$
where $P, Q, R:\R\to\R^{n\times n}$ are $C^1$-smooth and
$\tau$-periodic,  $R(t)=R(t)^T$, and each $P(t)=P(t)^T$ is also
positive definite. The corresponding Lagrangian system is
\begin{equation}\label{e:2.13}
\frac{d}{dt}\Big(\frac{\partial  L}{\partial v}(t, y,\dot
y)\Big)-\frac{\partial  L}{\partial y}(t,y,\dot y)= (P\dot{y}+
Qy)^\cdot- Q^T\dot{y}- Ry=0.
\end{equation}
Let $\tilde y$ be a critical point of the functional
$$
f_\tau(y)=\int^\tau_0L(t,y(t),\dot y(t))dt
$$
on $W^{1,2}(S_\tau,\R^n)$, and  the second differential of $f_\tau$
at it be given by
$$
d^2f_\tau(\tilde y)(y,z)=\int^\tau_0\left[(P\dot{y}+
Qy)\cdot\dot{z}+ Q^T\dot{y}\cdot z+ Ry\cdot z\right]dt.
$$
The linearized system of (\ref{e:2.13}) at $\tilde y$ is the Sturm
system:
$$
-(P\dot{y}+ Qy)^\cdot + Q^T\dot{y}+ Ry=0.
$$
 Let
\begin{equation}\label{e:2.14}
S(t)=\left(\begin{array}{cc}
P(t)^{-1}& -P(t)^{-1}Q(t)\\
-Q(t)^TP(t)^{-1}& Q(t)^TP(t)^{-1}Q(t)-R(t)\end{array}\right),
\end{equation}
and $\Psi:[0, +\infty)\to{\rm Sp}(2n,\R)$ be the fundamental
solution of the problem
\begin{equation}\label{e:2.15}
\dot{\bf u}(t)=J_0S(t){\bf u}
\end{equation}
with $\Psi(0)=I_{2n}$.  Suppose that each $P(t)$ is symmetric
positive definite, and that each $R(t)$ is symmetric. Then $f_\tau$
at $\tilde y\in W^{1,2}(S_\tau,\R^n)$ has finite \verb"Morse index"
$m_\tau(f_\tau, \tilde y)$ and \verb"nullity" $m^0_\tau(f_\tau, \tilde y)$,
and
\begin{equation}\label{e:2.16}
m^-_\tau(f_\tau, \tilde y)=i_\tau(\Psi)\quad{\rm and}\quad
m^0_\tau(f_\tau, \tilde y)=\nu_\tau(\Psi).
\end{equation}
\end{lemma}\vspace{2mm}

\begin{remark}\label{rm:2.3}\quad{\rm
Since ${L}_{vv}(t, y, v)=P(t)$ is invertible for every $t$, $L$ has
the Legendre transform $H:\R\times\R^{2n}\to\R$:
$$
 H(t, x, y)=x\cdot v(t, x, y)-L(t, x,
v(t, x,y)),
$$
where $v(t,x,y)\in\R^n$ is determined by $L_v(t,y, v(t, x,y))=x$.
Precisely, $v(t,x,y)=P(t)^{-1}[x-Q(t)y]$ and
\begin{eqnarray*}
H(t, x, y)&=&\frac{1}{2}P(t)^{-1}x\cdot x- P(t)^{-1}x\cdot Q(t)y\\
& + &
\frac{1}{2}P(t)^{-1}Q(t)y\cdot Q(t)y-\frac{1}{2}R(t)y\cdot y.
\end{eqnarray*}
Then $X_H(t, x, y)=J_0S(t){\bf u}$ with ${\bf u}=(x^T, y^T)^T$, and
$\tilde{\bf u}=(\tilde x^T, \tilde y^T)^T$ is a $\tau$-periodic
solution of (\ref{e:2.15}). }
\end{remark}

Let
\begin{eqnarray*}
&&EW^{1,2}(S_\tau, \R^n)=\left\{\,y\in W^{1,2}(S_\tau,
\R^n)\,|\,y(-t)=y(t)\,\forall
t\in\R\right\},\\
&&OW^{1,2}(S_\tau, \R^n)=\left\{\,y\in
W^{1,2}(S_\tau,\R^n)\,|\,y(-t)=-y(t)\,\forall t\in\R\right\}.
\end{eqnarray*}

\begin{lemma}\label{lem:2.4}\!\!{\rm (\cite[Th.3.4]{LuW2})}\quad
Under the assumptions of Lemma~\ref{lem:2.2}, suppose furthermore
that
\begin{equation}\label{e:2.17}
\left\{\begin{array}{ll}
&P(t+ \tau)=P(t)=P(t)^T=P(-t)\;\forall t\in\R,\\
& R(t+ \tau)=R(t)=R(t)^T=R(-t)\;\forall t\in\R,\\
& Q(t+ \tau)=Q(t)=-Q(-t)\;\forall t\in\R,
 \end{array}\right.
\end{equation}
and thus $L$ in Lemma~\ref{lem:2.2} satisfies {\rm (L4)}. So the
present $S(t)$ in (\ref{e:2.14}) also satisfies the Assumption {\rm
B}. Let $\tilde y$ be a critical point of the \verb"restriction"
$f^E_\tau$ of the functional $f_\tau$
 to $EW^{1,2}(S_\tau,\R^n)$.
(It is also  a critical point of the functional $f_\tau$
 on $W^{1,2}(S_\tau,\R^n)$ because $f_\tau$ is even).
As in Lemma~\ref{lem:2.1}, let $\Psi$ denote the fundamental
solution of (\ref{e:2.15}). Let
\begin{eqnarray*}
&&EW^{1,2}(S_\tau, \R^n)=EW^{1,2}(S_\tau, \R^n)^+\oplus
EW^{1,2}(S_\tau, \R^n)^0\oplus EW^{1,2}(S_\tau, \R^n)^-,\\
&&OW^{1,2}(S_\tau, \R^n)=OW^{1,2}(S_\tau, \R^n)^+\oplus
OW^{1,2}(S_\tau, \R^n)^0\oplus OW^{1,2}(S_\tau, \R^n)^-
\end{eqnarray*}
be  respectively $d^2f_\tau(\tilde y)$-orthogonal decompositions
 according to $d^2f_\tau(\tilde y)$ being positive, null,
and negative definite. Then
\begin{eqnarray}
&&\dim EW^{1,2}(S_\tau, \R^n)^-= m^-_\tau(f^E_\tau,
 \tilde y)=\mu_{1,\tau}(\Psi),\label{e:2.18}\\
&&\dim EW^{1,2}(S_\tau, \R^n)^0=m^0_\tau(f^E_\tau,
\tilde y)=\nu_{1,\tau}(\Psi),\label{e:2.19}\\
&&\dim OW^{1,2}(S_\tau, \R^n)^-=\mu_{2,\tau}(\Psi)-n,\label{e:2.20}\\
&& \nu_\tau(\Psi)=\nu_{1,\tau}(\Psi)+
\nu_{2,\tau}(\Psi).\label{e:2.21}
\end{eqnarray}
\end{lemma}

For conveniences we denote by
\begin{eqnarray}
&&m^-_{2,\tau}(f_\tau,
 \tilde y):=\dim OW^{1,2}(S_\tau, \R^n)^-,\label{e:2.22}\\
&&m^0_{2,\tau}(f_\tau,
\tilde y):=\dim OW^{1,2}(S_\tau, \R^n)^0.\label{e:2.23}
\end{eqnarray}
Then under the assumptions of Lemma~\ref{lem:2.4},
Lemma~\ref{lem:2.1}(ii)(iii) and (\ref{e:2.21}) become
\begin{eqnarray}
&& \bigl|n+ m^-_{2,\tau}(f_\tau,
 \tilde y)-m^-_\tau(f^E_\tau,
 \tilde y)\bigr|\le n,\label{e:2.24}\\
&&m^-_{2,\tau}(f_\tau,
 \tilde y)+ m^-_\tau(f^E_\tau,
 \tilde y)=m^-_\tau(f_\tau,
 \tilde y), \label{e:2.25}\\
&& m^0_\tau(f_\tau, \tilde y)=m^0_\tau(f^E_\tau, \tilde y)+
m^0_{2,\tau}(f_\tau, \tilde y).\label{e:2.26}
\end{eqnarray}

\section{Iteration inequalities of the Morse index}
\setcounter{equation}{0}

\noindent{\bf 3.1. The case of general periodic solutions}.\quad In
this subsection we always assume: $M$ is $C^3$-smooth, $L$ is
$C^2$-smooth and satisfies (L1)-(L3). Let $\gamma\in E_\tau$ be a
critical point of the functional ${\cal L}_\tau$ on $E_\tau$. It is
a $\tau$-periodic map from $\R$ to $M$. For each $k\in\N$,
$\gamma:\R\to M$ is also $k\tau$-periodic map and therefore
determines an element in $E_{k\tau}$, denoted by $\gamma^k$ for the
sake of clearness. It is not difficult to see that $\gamma^k$ is a
critical point of ${\cal L}_{k\tau}$ on $E_{k\tau}$. Let
$$
m^-_{k\tau}(\gamma^k)\quad{\rm and}\quad m^0_{k\tau}(\gamma^k)
$$
denote the \verb"Morse index" and \verb"nullity" of ${\cal
L}_{k\tau}$ on $E_{k\tau}$ respectively. Note that
$$
0\le m^0_{k\tau}(\gamma^k)\le 2n\quad\forall k\in\N.
$$
(This can be derived from (\ref{e:2.16}) and Lemma~\ref{lem:3.2}
below). A natural question is how to estimate $m^-(\gamma^k)$ in
terms of $m^-_{\tau}(\gamma)$, $m^0_{\tau}(\gamma)$ and
$m^0_{k\tau}(\gamma^k)$. The following theorem gives an answer.

\begin{theorem}\label{th:3.1}
 For  a critical point $\gamma$ of ${\cal
L}_\tau$ on $E_\tau$, assume that $\gamma^\ast TM\to S_\tau$ is
trivial. Then the mean Morse index
\begin{equation}\label{e:3.1}
\hat
m^-_{\tau}(\gamma):=\lim_{k\to\infty}\frac{m^-_{k\tau}(\gamma^k)}{k}
\end{equation}
always exists, and it holds that
\begin{equation}\label{e:3.2}
{\rm max}\left\{0, k\hat m^-_{\tau}(\gamma)-n\right\}\le
m^-_{k\tau}(\gamma^k)\le k\hat m_{\tau}(\gamma)+n-
m^0_{k\tau}(\gamma^k) \quad\forall k\in\N.
\end{equation}
Consequently, for any critical point $\gamma$ of ${\cal L}_\tau$ on
$E_\tau$, $\hat m^-_{2\tau}(\gamma^2)$ exists and
\begin{equation}\label{e:3.3}
{\rm max}\left\{0, k\hat m^-_{2\tau}(\gamma^2)-n\right\}\le
m^-_{2k\tau}(\gamma^{2k})\le k\hat m_{2\tau}(\gamma)+n-
m^0_{2k\tau}(\gamma^{2k}) \quad\forall k\in\N
\end{equation}
because $(\gamma^2)^\ast TM\to S_{2\tau}$ is always trivial.
\end{theorem}

Before proving this result it should be noted that the following
special case is a direct consequence of Lemma~\ref{lem:2.1}(i) and
Lemma~\ref{lem:2.2}.

\begin{lemma}\label{lem:3.2}
Under the assumptions of Lemma~\ref{lem:2.2}, for each $k\in\N$,
$\tilde y$ is also a $k\tau$-periodic solution of (\ref{e:2.13}),
denoted by $\tilde y^k$. Then $\tilde y^k$ is a critical point of
the functional
$$
f_{k\tau}(y)=\int^{k\tau}_0L(t,y(t),\dot y(t))dt
$$
on $W^{1,2}(S_{k\tau},\R^n)$, and
\begin{eqnarray}
  &&\hat m^-_\tau(f_\tau, \tilde y):=\lim_{k\to+\infty}
  \frac{m^-_{k\tau}(f_{k\tau}, \tilde y^k)}{k}=\lim_{k\to +\infty}\frac{i_{k\tau}(\Psi^k)}{k}=\hat i_\tau(\Psi),
  \label{e:3.4}\\
&&{\rm max}\{0, k\hat m^-_\tau(f_\tau, \tilde y)-n\}\le
m^-_{k\tau}(f_{k\tau}, \tilde y^k)\nonumber\\
&&\hspace{44mm}\le k\hat m^-_\tau(f_\tau, \tilde y)+ n-
m^0_{k\tau}(f_{k\tau}, \tilde y^k)\label{e:3.5}
\end{eqnarray}
with $0\le m^0_\tau(f_{k\tau}, \tilde y^k)\le 2n$ for any $k\in\N$.
\end{lemma}
This result  was actually used in \cite{Lo2, LoLu, LuW2}. In
the following we shall show that Theorem~\ref{th:3.1} can be reduced
to the special case. \vspace{2mm}

\noindent{\bf Proof of Theorem~\ref{th:3.1}}.\quad {\it Step 1.
Reduce to the case $M=\R^n$}. Let $\gamma\in E_\tau$ be  a critical
point $\gamma$ of ${\cal L}_\tau$ on $E_\tau$ with {\bf trivial
pullback $\gamma^\ast TM\to S_\tau$}. Take a $C^2$-smooth loop
$\gamma_0:S_\tau\to M$ such that $\max_t d(\gamma(t),
\gamma_0(t))<\rho$, where $d$ and $\rho$ are the distance and
injectivity radius of $M$ with respect to some chosen Riemannian
metric on $M$ respectively. (Actually we can choose
$\gamma_0=\gamma$ because $\gamma_0$ is $C^2$-smooth under the
assumptions of this subsection). Clearly, $\gamma$ and $\gamma_0$
are homotopic, and thus $\gamma_0^\ast TM\to S_\tau$ is trivial too.
 Since  $\gamma_0$ is $C^2$-smooth, we can choose a $C^2$-smooth
orthogonal trivialization
\begin{equation}\label{e:3.6}
S_\tau\times\R^n\to \gamma_0^\ast TM,\;(t, q)\mapsto \Phi(t)q.
\end{equation}
It naturally leads to a  smooth orthogonal trivialization of
$(\gamma_0^k)^\ast TM$ for any $k\in\N$,
\begin{equation}\label{e:3.7}
S_{k\tau}\times\R^n\to (\gamma_0^k)^\ast TM,\;(t, q)\mapsto
\Phi(t)q.
\end{equation}
 Let $B^n_\rho(0)$
denote an open ball in $\R^n$ centered at $0$ with radius $\rho$.
 Then for each $k\in \N$, we
have a coordinate chart on $E_{k\tau}$ containing $\gamma^k$,
\begin{equation}\label{e:3.8}
\phi_{k\tau}: W^{1,2}(S_{k\tau}, B^n_\rho(0))\to E_{k\tau},\quad
\phi_{k\tau}(\tilde\alpha)(t)=\exp_{\gamma_0^k(t)}(\Phi(t)\tilde\alpha(t)).
\end{equation}
Clearly, $\phi_{k\tau}(\tilde\alpha)$ has a period $\tau$ if and
only if $\tilde\alpha$ is actually $\tau$-periodic. Thus we have a
unique $\tilde\gamma\in W^{1,2}(S_{\tau}, B^n_\rho(0))$ such that
$\phi_{k\tau}(\tilde\gamma^k)=\gamma^k$ for any $k\in\N$. Denote by
the iteration maps
\begin{eqnarray*}
&&\psi^k: E_\tau\to E_{k\tau},\; \alpha\mapsto \alpha^k,\\
&&\psi^k: T_\alpha E_\tau\to T_{\alpha^k}E_{k\tau},\;\; \xi\mapsto
\xi^k,\\
&&\tilde\psi^k: W^{1,2}(S_{\tau}, \R^n)\to W^{1,2}(S_{k\tau},
\R^n),\; \tilde\alpha\mapsto \tilde\alpha^k.
\end{eqnarray*}
It is easy to see that
\begin{equation}\label{e:3.9}
\phi_{k\tau}\circ{\tilde\psi}^k=\psi^k\circ\phi_\tau\quad\forall
k\in\N.
\end{equation}
For any $k\in \N$, set
\begin{equation}\label{e:3.10}
\tilde{\cal L}_{k\tau}:W^{1,2}(S_{k\tau},
B^n_\rho(0))\to\R,\;\tilde{\cal L}_{k\tau}={\cal
L}_{k\tau}\circ\phi_{k\tau}.
\end{equation}
 Then
$\tilde\gamma=\phi_\tau^{-1}(\gamma)$ is a  critical point of
$\tilde{\cal L}_\tau$, and therefore
$\tilde\gamma^k=\phi_{k\tau}^{-1}(\gamma^k)=\tilde\psi^k(\tilde\gamma)$
is a critical point of  $\tilde{\cal L}_{k\tau}$ for any $k\in\N$.
Moreover, the Morse indexes and nullities of these critical points
satisfy the relations:
\begin{equation}\label{e:3.11}
m^-_{k\tau}(\tilde
\gamma^k)=m^-_{k\tau}(\gamma^k)\quad\hbox{and}\quad
m^0_{k\tau}(\tilde\gamma^k)=m^0_{k\tau}(\gamma^k),\quad \forall
k\in\N.
\end{equation}

Viewing $\gamma_0$ a $\tau$-periodic map from $\R\to M$, consider
the $C^2$-smooth map
\begin{equation}\label{e:3.12}
\Xi:\R\times B^n_\rho(0)\to  M,\quad (t, \tilde q)\mapsto
\exp_{\gamma_0(t)}\bigl(\Phi(t)\tilde q\bigr).
\end{equation}
Then $\Xi(t+\tau,\tilde q)=\Xi(t, \tilde q)$ for any $(t,\tilde
q)\in\R\times M$. Clearly,
\begin{eqnarray}
 &&  \phi_{k\tau}(\tilde\alpha)(t)=\Xi(t,
\tilde\alpha(t))\quad{\rm and}\label{e:3.13}\\
&&\frac{d}{dt}(\phi_{k\tau}(\tilde\alpha))(t)=\frac{d}{dt}\Xi(t,\tilde
q)|_{\tilde q=\tilde\alpha(t)}+ d_{\tilde q}\Xi(t,\tilde\alpha(t)
)(\dot{\tilde\alpha}(t))\label{e:3.14}
\end{eqnarray}
for any $t\in\R$ and $\tilde\alpha\in W^{1,2}(S_{k\tau},
B^n_\rho(0))$. Define $\tilde L: \R\times
 B^n_\rho(0)\times\R^n\to\R$ by
\begin{equation}\label{e:3.15}
\tilde L(t, \tilde q, \tilde v)=L\left(t, \Xi(t,\tilde q),
\frac{d}{dt}\Xi(t,\tilde q)+ d_{\tilde q}\Xi(t,\tilde q)(\tilde
v)\right).
\end{equation}
Then $\tilde L(t+ \tau, \tilde q, \tilde v)=\tilde L(t, \tilde q,
\tilde v)\;\forall (t, \tilde q, \tilde v)\in \R\times
 B^n_\rho(0)\times\R^n$, and $\tilde L$ also
satisfies the conditions (L2')-(L3') (up to changing the constants).
For $\tilde\alpha\in W^{1,2}(S_{k\tau}, B^n_\rho(0))$, by
(\ref{e:3.10}) we have
\begin{eqnarray}\label{e:3.16}
\tilde{\cal L}_{k\tau}(\tilde\alpha)&=&{\cal
L}_{k\tau}((\phi^k(\tilde\alpha))\nonumber\\
 &=&\int^{k\tau}_0 L\left(t,
\phi^k(\tilde\alpha)(t),
\frac{d}{dt}(\phi^k(\tilde\alpha))(t)\right)dt\nonumber\\
&=&\int^{k\tau}_0 \tilde L\left(t, \tilde\alpha(t),
\dot{\tilde\alpha}(t)\right)dt.
\end{eqnarray}
\textbf{Therefore we may assume} $M=\R^n$. That is, by
(\ref{e:3.11}) we only need to prove
\begin{eqnarray}
&&\hat
m^-_{\tau}(\tilde\gamma):=\lim_{k\to\infty}\frac{m^-_{k\tau}(\tilde\gamma^k)}{k}
\qquad \hbox{exists},\label{e:3.17} \\
&&{\rm max}\left\{0, k\hat m^-_{\tau}(\tilde\gamma)-n\right\}\le
m^-_{k\tau}(\tilde\gamma^k)\nonumber\\
&&\hspace{37mm} \le k\hat m_{\tau}(\tilde\gamma)+n-
m^0_{k\tau}(\tilde\gamma^k) \quad\forall k\in\N.\label{e:3.18}
\end{eqnarray}

\noindent{\it Step 2.\quad Reduce to the case of
Lemma~\ref{lem:3.2}}. Note that
\begin{eqnarray*}
d{\tilde L}_\tau(\tilde \gamma) (\tilde \xi)
 \!\!\!&=&\!\!\!\int_0^\tau  \left( D_{\tilde q} \tilde
L\left(t,\tilde\gamma(t),\dot{\tilde\gamma}(t)\right)(\tilde\xi(t))
+
D_{\tilde v} \tilde L\left(t,\tilde\gamma(t),\dot{\tilde\gamma}(t)\right)(\tilde\xi(t)) \right)  \, dt\nonumber \\
\!\!\!&=&\!\!\! \int_0^\tau \left(D_{\tilde q} \tilde
L\left(t,\tilde\gamma(t),\dot{\tilde\gamma}(t)\right) -
  \frac{d}{dt} D_{\tilde v} \tilde
  L\left(t,\tilde\gamma(t),\dot{\tilde\gamma}(t)\right)\right)
  \cdot\tilde\xi(t) \, dt
\end{eqnarray*}
for any $\tilde\xi\in W^{1,2}(S_\tau, \R^n)$. Since $d{\tilde
L}_\tau(\tilde \gamma)=0$, we have also
\begin{eqnarray*}
 d^2 \tilde{L}_\tau
  (\tilde\gamma)(\tilde\xi,\tilde\eta)
   = \int_0^\tau \Bigl(\!\! \!\!\!&&\!\!\!\!\!D_{\tilde v\tilde v}
    \tilde L\left(t,\tilde\gamma(t),\dot{\tilde\gamma}(t)\right)
\bigl(\dot{\tilde\xi}(t), \dot{\tilde \eta}(t)\bigr) \\
&&+ D_{\tilde q\tilde v} \tilde
  L\left(t,\tilde\gamma(t), \dot{\tilde\gamma}(t)\right)
\bigl(\tilde\xi(t), \dot{\tilde\eta}(t)\bigr)\nonumber \\
&& + D_{\tilde v\tilde q} \tilde
  L\left(t,\tilde\gamma(t),\dot{\tilde\gamma}(t)\right)
\bigl(\dot{\tilde\xi}(t), \tilde\eta(t)\bigr) \\
&&+  D_{\tilde q\tilde q} \tilde L\left(t,\tilde\gamma(t),
\dot{\tilde\gamma}(t)\right) \bigl(\tilde\xi(t),\tilde
\eta(t)\bigr)\Bigr) \, dt
\end{eqnarray*}
for any $\tilde\xi,\tilde\eta\in W^{1,2}(S_\tau, \R^n)$. Set
\begin{equation}\label{e:3.19}
\left.\begin{array}{ll} &\hat P(t)=D_{\tilde v\tilde v} \tilde
L\left(t,\tilde\gamma(t),\dot{\tilde\gamma}(t)\right),\\
&\hat Q(t)=D_{\tilde q\tilde v} \tilde
  L\left(t,\tilde\gamma(t), \dot{\tilde\gamma}(t)\right),\\
&\hat R(t)=D_{\tilde q\tilde q} \tilde L\left(t,\tilde\gamma(t),
\dot{\tilde\gamma}(t)\right)
\end{array}\right\}
\end{equation}
 and
\begin{equation}\label{e:3.20}
 \hat L(t, \tilde y, \tilde v)=\frac{1}{2}\hat P(t)\tilde v\cdot\tilde v + \hat Q(t)\tilde y\cdot\tilde v +
\frac{1}{2}\hat R(t)\tilde y\cdot\tilde y.
\end{equation}
Clearly,
they satisfy the conditions of Lemma~\ref{lem:2.2}, and $\tilde
y=0\in W^{1,2}(S_{\tau},\R^n)$ is a critical point of the functional
$$
\hat f_{\tau}(\tilde y)=\int^{\tau}_0\hat L\left(t, \tilde y(t),\dot
{\tilde y}(t)\right)dt
$$
on $W^{1,2}(S_{\tau},\R^n)$. It is also easily
checked that
$$
d^2 \hat{f}_\tau
  (0)(\tilde\xi,\tilde\eta)
=d^2 \tilde{L}_\tau
  (\tilde\gamma)(\tilde\xi,\tilde\eta)\quad\forall
\tilde\xi,\tilde\eta\in W^{1,2}(S_\tau, \R^n).
$$
It follows that
$$
m^-_{k\tau}(\hat f_{k\tau}, 0)=m^-_{k\tau}(\tilde\gamma^k)
\quad\hbox{and}\quad m^0_{k\tau}(\hat f_{k\tau},
0)=m^0_{k\tau}(\tilde\gamma^k)\quad\forall k\in\N.
$$
These and Lemma~\ref{lem:3.2} together give the desired
(\ref{e:3.17}) and (\ref{e:3.18}). $\Box$\vspace{3mm}

\noindent{\bf 3.2. The case of even periodic solutions}.\quad Let
$M$ and $L$ be as in \S3.1. But we also assume that $L$ satisfies
(L4). Note that the \textsf{even periodic solutions are always
contractible}. Let ${\cal L}^E_{k\tau}$ denote the restriction of
${\cal L}_{k\tau}$ on $EH_{k\tau}$. As noted in the introduction, if
$\gamma\in EH_\tau$ is a critical point of ${\cal L}_{\tau}^E$ on
$EH_\tau$ then $\gamma^k$ is a critical point of ${\cal L}_{k\tau}$
on $H_{k\tau}$ for each $k\in\N$. Let
$$m^-_{1,
k\tau}(\gamma^k)\quad{\rm and}\quad m^0_{1, k\tau}(\gamma^k)
$$
denote the \verb"Morse index" and \verb"nullity" of ${\cal
L}^E_{k\tau}$ on $EH_{k\tau}$ respectively. Then $0\le m^0_{1,
k\tau}(\gamma^k)\le m^0_{k\tau}(\gamma^k)\le 2n$ for any $k$. We
shall prove

\begin{theorem}\label{th:3.3}
 Let $L$ satisfy the conditions (L1)-(L4). Then for any critical point $\gamma$ of ${\cal
L}^E_\tau$ on $EH_\tau$, the mean Morse index
\begin{equation}\label{e:3.21}
\hat m^-_{1,\tau}(\gamma):=\lim_{k\to\infty}\frac{m^-_{1,
k\tau}(\gamma^k)}{k}
\end{equation}
 exists, and it holds that
\begin{equation}\label{e:3.22}
m^-_{1, k\tau}(\gamma^k)+ m^0_{1, k\tau}(\gamma^k)\le n \quad\forall
k\in\N \quad{\rm if}\quad \hat m^-_{1,\tau}(\gamma)=0.
\end{equation}
\end{theorem}

Firstly, by (\ref{e:2.10}) and (\ref{e:2.16}) the mean Morse index
\begin{equation}\label{e:3.23}
\hat m^-_{\tau}(f_\tau, \tilde
y):=\lim_{k\to\infty}\frac{m^-_{k\tau}(f_{k\tau}, \tilde y^k)}{k}
\end{equation}
 exists and equals to $\hat i_\tau(\Psi)$.
Under the assumptions of Lemma~\ref{lem:2.4},
 for each $k\in\N$,
 $\tilde y^k$ is a critical point of the \textsf{restriction} $f^E_{k\tau}$
of the functional $f_{k\tau}$  to $EW^{1,2}(S_{k\tau},\R^n)$, and it
follows from (\ref{e:2.12}), (\ref{e:2.18}), (\ref{e:2.20}) and
(\ref{e:2.22}) that
\begin{eqnarray}
&&\hat m^-_\tau(f^E_\tau,
 \tilde y):=\lim_{k\to +\infty}\frac{m^-_{k\tau}(f^E_{k\tau}, \tilde y^k)}{k}=\hat
 \mu_{1,\tau}(\Psi)=\frac{1}{2}\hat m^-_{\tau}(f_\tau, \tilde
y), \label{e:3.24}\\
&&\hat m^-_{2,\tau}(f_\tau,
 \tilde y):=\lim_{k\to +\infty}\frac{m^-_{2, k\tau}(f_{k\tau}, \tilde y^k)}{k}=
 \hat \mu_{2,\tau}(\Psi)=\frac{1}{2}\hat m^-_{\tau}(f_\tau, \tilde
y).\label{e:3.25}
\end{eqnarray}
Moreover, by (\ref{e:2.25}) and (\ref{e:2.26}), for any $k\in\N$ it
holds that
\begin{eqnarray*}
&&m^-_{2, k\tau}(f_{k\tau},
 \tilde y^k)+ m^-_{k\tau}(f^E_{k\tau},
 \tilde y^k)=m^-_{k\tau}(f_{k\tau},
 \tilde y^k), \\
&& m^0_{k\tau}(f_{k\tau}, \tilde y^k)=m^0_{k\tau}(f^E_{k\tau},
\tilde y^k)+ m^0_{2, k\tau}(f_{k\tau}, \tilde y^k).
\end{eqnarray*}
From these we derive that (\ref{e:3.5}) becomes
\begin{eqnarray}
  &&{\rm max}\{0, 2k\hat m^-_\tau(f^E_\tau, \tilde y)-n\}\le
m^-_{2, k\tau}(f_{k\tau},
 \tilde y^k)+ m^-_{k\tau}(f^E_{k\tau},
 \tilde y^k)\nonumber\\
&&\hspace{14mm}\le 2k\hat m^-_\tau(f^E_\tau, \tilde y)+ n-
m^0_{k\tau}(f^E_{k\tau}, \tilde y^k)- m^0_{2, k\tau}(f_{k\tau},
\tilde y^k)\label{e:3.26}
\end{eqnarray}
for any $k\in\N$. In particular, if $\hat m^-_\tau(f^E_\tau,
 \tilde y)=0$, then
\begin{equation}\label{e:3.27}
 m^-_{k\tau}(f^E_{k\tau},
 \tilde y^k)+ m^0_{k\tau}(f^E_{k\tau}, \tilde y^k)\le n\quad\forall k\in\N.
\end{equation}
(\cite[Th.3.7]{LuW2}).\vspace{2mm}

\noindent{\bf Proof of Theorem~\ref{th:3.3}}.\quad Since $\gamma$ is
even we can still choose $\gamma_0$ and $\Phi$ in (\ref{e:3.6}) to
be even, i.e. $\gamma_0(-t)=\gamma_0(t)$ and $\Phi(-t)=\Phi(t)$ for
any $t\in\R$. These imply
\begin{equation}\label{e:3.28}
\Xi(-t, \tilde q)=\Xi(t, \tilde q),\quad \frac{d}{dt}\Xi(-t,\tilde
q)=-\frac{d}{ds}\Xi(s,\tilde q)|_{s=-t}=\frac{d}{dt}\Xi(t,\tilde q).
\end{equation}
It follows that  the coordinate chart $\phi_{k\tau}$ in
(\ref{e:3.8}) naturally restricts to a coordinate chart on
$EH_{k\tau}$,
\begin{equation}\label{e:3.29}
\phi^E_{k\tau}: EW^{1,2}\bigl(S_{k\tau}, B^n_\rho(0)\bigr)\to
EH_{k\tau}
\end{equation}
which also satisfies
\begin{equation}\label{e:3.30}
\phi^E_{k\tau}\circ{\tilde\psi}^k=\psi^k\circ\phi^E_\tau\quad\forall
k\in\N.
\end{equation}
By (L4), (\ref{e:3.15}) and (\ref{e:3.28}) we have
\begin{eqnarray}\label{e:3.31}
\tilde L(-t, \tilde q, -\tilde v)&=&L\left(-t, \Xi(-t,\tilde q),
\frac{d}{d(-t)}\Xi(-t,\tilde q)+ d_{\tilde q}\Xi(-t,\tilde
q)(-\tilde
v)\right)\nonumber\\
&=&L\left(-t, \Xi(t,\tilde q), -\frac{d}{dt}\Xi(-t,\tilde q)-
d_{\tilde q}\Xi(t,\tilde q)(\tilde v)\right)\nonumber\\
&=&L\left(t, \Xi(t,\tilde q), \frac{d}{dt}\Xi(-t,\tilde q)+
d_{\tilde q}\Xi(t,\tilde q)(\tilde v)\right)\nonumber\\
&=&L\left(t, \Xi(t,\tilde q), \frac{d}{dt}\Xi(t,\tilde q)+ d_{\tilde
q}\Xi(t,\tilde q)(\tilde v)\right)\nonumber\\
&=&\tilde L(t, \tilde q, \tilde v).
\end{eqnarray}
That is, $\tilde L$ also satisfies (L4). It follows that for any
$k\in \N$, the functional
\begin{equation}\label{e:3.32}
\tilde{\cal L}^E_{k\tau}:EW^{1,2}\bigl(S_{k\tau},
B^n_\rho(0)\bigr)\to\R,\;\tilde{\cal L}^E_{k\tau}={\cal
L}^E_{k\tau}\circ\phi^E_{k\tau}
\end{equation}
is exactly the restriction of the functional $\tilde{\cal
L}_{k\tau}$ in (\ref{e:3.10}) to $EW^{1,2}(S_{k\tau}, B^n_\rho(0))$.
Hence the question is reduced to the case $M=\R^n$ again. That is,
we only need to prove
\begin{eqnarray}
&&\hat m^-_{1,\tau}(\tilde\gamma):=\lim_{k\to\infty}\frac{m^-_{1,
k\tau}(\tilde\gamma^k)}{k}\quad\hbox{exists},\label{e:3.33}\\
&&m^-_{1, k\tau}(\tilde\gamma^k)+ m^0_{1, k\tau}(\tilde\gamma^k)\le
n \quad\forall k\in\N \quad{\rm if}\quad \hat
m^-_{1,\tau}(\tilde\gamma)=0.\label{e:3.34}
\end{eqnarray}

By (\ref{e:3.31}) we have
\begin{eqnarray*}
&&D_{\tilde v\tilde v}\tilde L(-t, \tilde q, -\tilde v)=D_{\tilde
v\tilde v}\tilde
L(t, \tilde q, \tilde v),\\
&&D_{\tilde q\tilde v}\tilde L(-t, \tilde q, -\tilde v)=-D_{\tilde
q\tilde v}\tilde
L(t, \tilde q, \tilde v),\\
&&D_{\tilde q\tilde q}\tilde L(-t, \tilde q, -\tilde v)=D_{\tilde
q\tilde q}\tilde L(t, \tilde q, \tilde v)
\end{eqnarray*}
for any $(t, \tilde q, \tilde v)\in\R\times B^n_\rho(0)\times\R^n$.
Since $\tilde\gamma(-t)=\tilde\gamma(t)$ and
$\dot{\tilde\gamma}(-t)=-\dot{\tilde\gamma}(t)$, it follows from
this that $\hat P$, $\hat Q$ and $\hat R$ in (\ref{e:3.19}) satisfy
(\ref{e:2.17}). For $\hat L$ in (\ref{e:3.20}) and the functionals
$$
\hat f^E_{k\tau}(\tilde y):=\int^{k\tau}_0\hat L\left(t, \tilde
y(t),\dot {\tilde y}(t)\right)dt
 $$
 on $EW^{1,2}(S_{k\tau},\R^n)$, $k=1,2,\cdots$, we have
\begin{equation}\label{e:3.35}
 m^-_{k\tau}(\hat f^E_{k\tau},
0)=m^-_{1, k\tau}(\tilde\gamma^k) \quad\hbox{and}\quad
m^0_{k\tau}(\hat f^E_{k\tau}, 0)=m^0_{1,
k\tau}(\tilde\gamma^k)\quad\forall k\in\N.
\end{equation}
By (\ref{e:3.24}) and (\ref{e:3.27}) we get
\begin{equation}\label{e:3.36}
\hat m^-_\tau(\hat f^E_\tau,
 0):=\lim_{k\to +\infty}\frac{m^-_{k\tau}(\hat f^E_{k\tau}, 0)}{k}
\end{equation}
exists, and if $\hat m^-_\tau(\hat f^E_\tau,
 0)=0$,
\begin{equation}\label{e:3.37}
 m^-_{k\tau}(\hat f^E_{k\tau},
 0)+ m^0_{k\tau}(\hat f^E_{k\tau},
0)\le n\quad\forall k\in\N.
\end{equation}
 Now (\ref{e:3.35})-(\ref{e:3.37}) give
 (\ref{e:3.33}) and (\ref{e:3.34}), and therefore the desired
 (\ref{e:3.21}) and (\ref{e:3.22}). $\Box$

\section{Critical modules under iteration maps}
\setcounter{equation}{0}

In this section we shall study relations of critical modules under
iteration maps in three different cases. We first recall a few of
notions.  Let ${\cal M}$ be a $C^2$ Hilbert-Riemannian manifold and
$f\in C^1({\cal M}, \R)$ satisfies the Palais-Smale condition.
Denote by ${\cal K}(f)$ the set of critical points of $f$.
 Recall that a connected submanifold $N$ of
${\cal M}$ is a {\it critical submanifold} of $f$ if it is closed,
 consists entirely of critical points of $f$
 and $f|_N={\rm constant}$. Let $N\subset {\cal M}$ be an isolated critical
 submanifold of $f$ with $f|_N=c$, and $U$ be a neighborhood of $N$ such that $U\cap{\cal
K}(f)=N$. For $q\in\N\cup\{0\}$, recall that the $q^{th}$ critical
group with coefficient group $\K$ of $f$ at $N$ is defined by
\begin{equation}\label{e:4.1}
C_q(f,N;\K):=H_q\bigl(\{f\le c\}\cap U, (\{f\le c\}\setminus N)\cap
U; \K\bigr).
\end{equation}
 Hereafter $H_\ast(X,Y;\K)$ stands for the relative singular homology with
the abelian coefficient group
 $\K$  without special statements.
 The group $C_q(f,N;\K)$ does not depend on a special choice of such
neighborhoods $U$ up to isomorphisms. There also exists another
equivalent definition of critical groups, which is convenient in
many situations.

Let $V:({\cal M}\setminus{\cal K}(f))\to T{\cal M}$ be a
pseudo-gradient vector field for $f$ on ${\cal M}$. According to
\cite[pp.48, 74]{Ch} and \cite[Def.2.3]{Wa} or \cite{GM1},
 a pair of topological subspaces $(W, W^-)$ of ${\cal M}$ is called a
{\bf Gromoll-Meyer pair} with respect to $V$ for $N$, if\\
  (1) $W$ is a closed neighborhood of $N$ possessing the
\verb"mean value property", i.e., $\forall t_1<t_2$, $\eta(t_i)\in
W$, $i=1,2$, implies $\eta(t)\in W$ for all $t\in [t_1, t_2]$, where
$\eta(t)$ is the decreasing flow with respect to $V$. And there
exists $\epsilon>0$ such that $W\cap
f_{c-\epsilon}=f^{-1}[c-\epsilon, c)\cap{\cal K}(f)=\emptyset$,
$W\cap{\cal K}(f)=N$;\\
  (2) the set $W^-=\{p\in W\, |\, \eta(t,p)\notin W, \forall t>0\}$;\\
  (3) $W^-$ is a piecewise submanifold, and the flow $\eta$ is
transversal to $W^-$.\vspace{2mm}\\
By \cite[pp.74]{Ch} or \cite[\S 2]{Wa}, there exists an (arbitrarily
small) Gromoll-Meyer pair for $N$, $(W, W^-)$, and for such a pair
it holds that
\begin{equation}\label{e:4.2}
H_\ast(W, W^-;\K)\cong C_\ast(f,N;\K).
\end{equation}
Hence $H_\ast(W, W^-;\K)$ may be used to  give an equivalent
definition of $C_\ast(f,N;\K)$.  We need the following fact which
seems to be obvious, but is often neglected.

\begin{lemma}\label{lem:4.1}
Let ${\cal M}_1$ and ${\cal M}_2$ be $C^2$ Hilbert-Riemannian
manifolds, and $\Theta:{\cal M}_1\to {\cal M}_2$ be a homeomorphism.
Suppose that $f_i\in C^1({\cal M}_i, \R)$, $i=1,2$, satisfy the
Palais-Smale condition and $f_2=f_1\circ\Theta$. Let $N_1\subset
{\cal M}_1$ and $N_2=\Theta(N_1)\subset{\cal M}_2$ be isolated
critical submanifolds of $f_1$ and $f_2$ respectively. Assume that
$(W_1, W_1^-)$ is a Gromoll-Meyer pair of $N_1$ of $f_1$. Then
$$
C_\ast(f_2, N_2; \K)\cong H_\ast(\Theta(W_1), \Theta(W_1^-);\K)
$$
though $(\Theta(W_1), \Theta(W_1^-))$ is not necessarily a
Gromoll-Meyer pair of $N_2$ of $f_2$ (because $\Theta$ is only a
homeomorphism). Moreover, for $c=f_1|_{N_1}$ and $\epsilon>0$ it is
clear that
$$
(W_1, W_1^-)\subset \bigl(f_1^{-1}[c-\epsilon,  c+ \epsilon],
f_1^{-1}(c-\epsilon)\bigr)
$$
implies $ \bigl(\Theta(W_1), \Theta(W_1^-)\bigr)\subset \bigl(
f_2^{-1}[c-\epsilon, c+ \epsilon], f_2^{-1}(c-\epsilon)\bigr)$.
\end{lemma}

\noindent{\bf Proof.}\quad Take a small open neighborhood $U$ of
$N_1$ so that $U\subset W_1$.  Since $\Theta\bigl(\{f_1\le c\}\cap
U\bigr)= \{f_2\le c\}\cap U$ and $\Theta\bigl((\{f_1\le c\}\setminus
N_1)\cap U\bigr)=(\{f_2\le c\}\setminus N_2)\cap \Theta(U)$,  we
have isomorphisms
\begin{eqnarray*}
&&\Theta_\ast: H_\ast(W_1, W_1^-; \K)\to H_\ast(\Theta(W_1),
\Theta(W_1^-);\K),\\
&&\Theta_\ast: H_\ast\bigl(\{f_1\le c\}\cap U, (\{f_1\le
c\}\setminus N_1)\cap U; \K\bigr)\to\\
&&\hspace{20mm} H_\ast\bigl(\{f_2\le c\}\cap \Theta(U), (\{f_2\le
c\}\setminus N_2)\cap \Theta(U); \K\bigr)\\
&&\hspace{20mm} =C_\ast(f_2, N_2;\K).
\end{eqnarray*}
By (\ref{e:4.1}) and (\ref{e:4.2}), $H_\ast(W_1, W_1^-; \K)\cong
H_\ast\bigl(\{f_1\le c\}\cap U, (\{f_1\le c\}\setminus N_1)\cap U;
\K\bigr)$. The desired conclusion is obtained. $\Box$.

It is this result that we may often treat $(\Theta(W_1),
\Theta(W_1^-))$ as a Gromoll-Meyer pair without special statements.
For conveniences we call it a {\bf topological Gromoll-Meyer} of
$f_2$ at $N_2$.   The usual Gromoll-Meyer pair can be viewed the
special case of it. Moreover, if $\Gamma:{\cal M}_2\to {\cal M}_3$
is a $C^1$-diffeomorphism onto another $C^2$ Hilbert-Riemannian
manifold ${\cal M}_3$, then $(\Gamma\circ\Theta(W_1),
\Gamma\circ\Theta(W_1^-))$ is also a topological Gromoll-Meyer pair
of $f_3=f_2\circ\Gamma^{-1}$ at $N_3=\Gamma(N_2)$. (\ref{e:4.2}) and
Lemma~\ref{lem:4.1} show that the topological Gromoll-Meyer may be
used to give an equivalent definition of the critical group.

To understand the Note at the end of proof of Theorem~5.1 of
\cite[pp. 44]{Ch} we add a lemma, which is need in this paper.

\begin{lemma}\label{lem:4.2}
Let $H_i$ be Hilbert spaces with origins $\theta_i$, , $i=1,2,3$.
For $\varepsilon>0$ let $f\in C^2(B_\varepsilon(\theta_1)\times
B_\varepsilon(\theta_2)\times B_\varepsilon(\theta_3), \R)$. Assume
that $d_3f(x_1, \theta_2, \theta_3)=0$ for $x_1\in
B_\varepsilon(\theta_1)$ and that
$d^2_3f(\theta_1,\theta_2,\theta_3):H_3\to H_3$ is a Banach space
isomorphism. Then there exist a small $0<\delta\ll\varepsilon$ and
$C^1$-map $h:B_\delta(\theta_1)\times B_\delta(\theta_2)\to H_3$
such that\\
(i) $d_3f(x_1, x_2, h(x_1, x_2))=\theta_3$ for all $(x_1, x_2)\in
B_\delta(\theta_1)\times B_\delta(\theta_2)$,\\
(ii) $g: B_\delta(\theta_1)\times B_\delta(\theta_2)\to\R,\; (x_1,
x_2)\mapsto g(x_1, x_2)= f(x_1, x_2, h(x_1, x_2))$ is $C^2$.
\end{lemma}

\noindent{\bf Proof.}\quad Applying the implicit function theorem to
the map
$$
d_3f:B_\varepsilon(\theta_1)\times B_\varepsilon(\theta_2)\times
B_\varepsilon(\theta_3)\to H_3
$$
we get a $0<\delta\ll \varepsilon$ and a $C^1$-map
$h:B_\delta(\theta_1)\times B_\delta(\theta_2)\to H_3$ such that
$h(\theta_1, \theta_2)=\theta_3$ and
$$
d_3f(x_1, x_2, h(x_1, x_2))=0\;\forall (x_1, x_2)\in
B_\delta(\theta_1)\times B_\delta(\theta_2).
$$
Set $g(x_1, x_2)= f(x_1, x_2, h(x_1, x_2))$. Then
\begin{eqnarray*}
dg(x_1, x_2)&=&d_{(1, 2)}f(x_1, x_2, h(x_1, x_2))+ d_3f(x_1, x_2,
h(x_1, x_2))\circ d_{(x_1, x_2)}h(x_1, x_2)\\
&=&d_{(1, 2)}f(x_1, x_2, h(x_1, x_2))
\end{eqnarray*}
because  $d_3f(x_1, x_2, h(x_1, x_2))=0$, where $d_{(1,2)}$ denotes
the differential for the first two variables of $f$. Hence
\begin{eqnarray*}
d^2g(x_1, x_2)&=& d^2_{(1, 2)}f(x_1, x_2, h(x_1, x_2))\\
&+& d_3d_{(1,2)}f(x_1, x_2, h(x_1, x_2))\circ d_{(x_1, x_2)}h(x_1,
x_2).
\end{eqnarray*}
The desired claims are proved. $\Box$\vspace{2mm}

\noindent{\bf 4.1.} The arguments in this section are following
Section 3 in \cite{Lo2}. However, since our arguments are on a
Hilbert manifold, rather than Hilbert space, some new techniques are
needed. The precise proofs are also given for reader's convenience.
In this subsection we always assume: $M$ is $C^3$-smooth, $L$ is
$C^2$-smooth and satisfies (L1)-(L3).

\begin{lemma}\label{lem:4.3}
Let $\gamma\in H_\tau(\alpha)$ be an isolated critical point of the
functional ${\cal L}_\tau$ on $H_\tau(\alpha)$ such that $\gamma^k$
is an isolated critical point of the functional ${\cal L}_{k\tau}$
in $H_{k\tau}(\alpha^{k})$ for some $k\in\N$. Suppose that
$\gamma^\ast T^\ast M\to S_\tau$ is trivial. Then there exist
Gromoll-Meyer pairs
 $\bigl(W(\gamma),  W(\gamma)^-\bigr)$ of ${\cal L}_\tau$ at
 $\gamma$ and $\bigl(W(\gamma^k),  W(\gamma^k)^-\bigr)$ of ${\cal L}_{k\tau}$ at
 $\gamma^k$ such that
 \begin{equation}\label{e:4.3}
\bigl(\psi^k(W(\gamma)),  \psi^k(W(\gamma)^-)\bigr)\subset
\bigl(W(\gamma^k),  W(\gamma^k)^-\bigr).
 \end{equation}
\end{lemma}

\noindent{\bf Proof.}\quad  For each $j\in \N$, let
\begin{equation}\label{e:4.4}
\phi_{j\tau}: W^{1,2}(S_{j\tau}, B^n_\rho(0))\to
H_{j\tau}(\alpha^{j})\quad {\rm and}\quad \tilde{\cal
L}_{j\tau}={\cal L}_{j\tau}\circ\phi_{j\tau}
\end{equation}
as in (\ref{e:3.8}) and (\ref{e:3.10}). They satisfy (\ref{e:3.9}),
i.e.
$\phi_{j\tau}\circ{\tilde\psi}^j=\psi^j\circ\phi_\tau\quad\forall
j\in\N$, where $\psi^j: H_\tau(\alpha)\to H_{j\tau}(\alpha^{j})$ and
$\tilde\psi^j: W^{1,2}(S_{\tau}, \R^n)\to W^{1,2}(S_{j\tau}, \R^n)$
are the iteration maps. Let $\tilde\gamma=(\phi_\tau)^{-1}(\gamma)$.
Then $\phi_{j\tau}(\tilde\gamma^j)=\gamma^j$ for any $j\in\N$.

Let $\|\cdot\|_\tau$ and $\|\cdot\|_{k\tau}$ denote the norms in
$W^{1,2}(S_{\tau}, \R^n)$ and $W^{1,2}(S_{k\tau}, \R^n)$
respectively. By the construction on page 49 of \cite{Ch}, we set
\begin{eqnarray*}
&&\tilde W(\tilde\gamma):={\cal L}_\tau^{-1}[c-\varepsilon, c+
\varepsilon]\cap\bigl\{x\in W^{1,2}(S_{\tau}, \R^n)\,|\,
\lambda{\cal L}_\tau(x)+ \|x\|^2_\tau\le \mu\bigr\},\\
&&\tilde W(\tilde\gamma)^-:={\cal
L}_\tau^{-1}(c-\varepsilon)\cap\bigl\{x\in W^{1,2}(S_{\tau},
\R^n)\,|\,
\lambda{\cal L}_\tau(x)+ \|x\|^2_\tau\le \mu\bigr\},\\
&&\tilde W(\tilde\gamma^k):={\cal L}_{k\tau}^{-1}[kc-k\varepsilon,
kc+ k\varepsilon]\cap\bigl\{y\in W^{1,2}(S_{k\tau}, \R^n)\,|\,
\lambda{\cal L}_{k\tau}(y)+ \|y\|^2_{k\tau}\le k\mu\bigr\},\\
&&\tilde W(\tilde\gamma^k)^-:={\cal L}_{k\tau}^{-1}(kc-
k\varepsilon)\cap\bigl\{y\in W^{1,2}(S_{k\tau}, \R^n)\,|\,
\lambda{\cal L}_{k\tau}(y)+ \|y\|^2_{k\tau}\le k\mu\bigr\},
\end{eqnarray*}
where positive numbers $\lambda, \mu, \varepsilon$ and $k\lambda,
k\mu, k\varepsilon$ are such that the conditions as in (5.13)-(5.15)
 on page 49 of \cite{Ch} hold. Then $\bigl(\tilde W(\tilde\gamma),
\tilde W(\tilde\gamma)^-\bigr)$ and $\bigl(\tilde W(\tilde\gamma^k),
\tilde W(\tilde\gamma^k)^-\bigr)$ are
 Gromoll-Meyer pairs of $\tilde{\cal L}_\tau$ at $\tilde\gamma$ and
 of $\tilde{\cal L}_{k\tau}$ at $\tilde\gamma^k$, and
\begin{equation}\label{e:4.5}
\bigl(\tilde\psi^k(\tilde W(\tilde\gamma)), \tilde\psi^k(\tilde
W(\tilde\gamma)^-)\bigr)\subset \bigl(\tilde W(\tilde\gamma^k),
\tilde W(\tilde\gamma^k)^-\bigr).
 \end{equation}
Define
\begin{equation}\label{e:4.6}
\left.\begin{array}{ll} &\bigl(W(\gamma),
W(\gamma)^-\bigr):=\bigl(\phi_\tau(\tilde
W(\tilde\gamma)), \phi_\tau(\tilde W(\tilde\gamma)^-)\bigr),\\
&\bigl(W(\gamma^k), W(\gamma^k)^-\bigr):=\bigl(\phi_{k\tau}(\tilde
W(\tilde\gamma^k)), \phi_{k\tau}(\tilde W(\tilde\gamma^k)^-)\bigr).
\end{array}\right\}
\end{equation}
Since $\phi_{k\tau}\circ{\tilde\psi}^k=\psi^k\circ\phi_\tau$,
(\ref{e:4.3}) follows from (\ref{e:4.5}). $\Box$

When $\gamma$ and $\gamma^k$ are isolated, according to the
definition of critical groups in (\ref{e:4.1}) it is easy to see
that the iteration map $\psi^k: H_\tau(\alpha)\to
H_{k\tau}(\alpha^{k})$ induces homomorphisms
$$
(\psi^k)_\ast: C_\ast({\cal L}_\tau, \gamma;\K)\to C_\ast({\cal
L}_{k\tau}, \gamma^k;\K).
$$
Lemma~\ref{lem:4.3} shows that the homomorphisms are still
well-defined when  the critical groups $C_\ast({\cal L}_\tau,
\gamma;\K)$ and $C_\ast({\cal L}_{k\tau}, \gamma^k;\K)$ are defined
by (\ref{e:4.2}). Later similar cases are always understand in this
way. Our purpose is to prove:

\begin{theorem}\label{th:4.4}
Let $\gamma\in H_\tau(\alpha)$ be an isolated critical point of the
functional ${\cal L}_\tau$ on $H_\tau(\alpha)$ such that
$\gamma^\ast TM\to S_\tau$ is trivial. Suppose that for some
$k\in\N$ the iteration $\gamma^k$ is also an isolated critical point
of the functional ${\cal L}_{k\tau}$ in $H_{k\tau}(\alpha^k)$, and
\begin{equation}\label{e:4.7}
m^-_{k\tau}(\gamma^k)=m^-_{\tau}(\gamma)\quad{\rm and}\quad
  m^0_{k\tau}(\gamma^k)=m^0_{\tau}(\gamma).
  \end{equation}
Then for $c={\cal L}_\tau(\gamma)$ and any $\epsilon>0$ there exist
topological Gromoll-Meyer pairs of ${\cal L}_\tau$ at $\gamma\in
H_\tau(\alpha)$ and
 of ${\cal L}_{k\tau}$ at $\gamma^k\in H_{k\tau}(\alpha^k)$,
\begin{eqnarray*}
&&(W_\tau,  W^-_\tau)\subset\bigl(({\cal L}_\tau)^{-1}[c-\epsilon,
c+ \epsilon], ({\cal L}_\tau)^{-1}(c-\epsilon)\bigr)
 \quad\hbox{and}\\
&&(W_{k\tau}, W^-_{k\tau})\subset\bigl(({\cal
L}_{k\tau})^{-1}[kc-k\epsilon, kc+ k\epsilon], ({\cal
L}_{k\tau})^{-1}(kc-k\epsilon)\bigr),
\end{eqnarray*}
such that
\begin{equation}\label{e:4.8}
 (\psi^k(W_\tau), \psi^k(W^-_\tau))\subset
(W_{k\tau}, W^-_{k\tau})
\end{equation}
and that the homomorphism
\begin{eqnarray}\label{e:4.9}
&&(\psi^k)_\ast: C_\ast({\cal L}_\tau, \gamma;\K):=H_\ast(W_\tau,
W^-_\tau; \K)\nonumber\\
&&\hspace{40mm}\to C_\ast({\cal L}_{k\tau},
\gamma^k;\K):=H_\ast(W_{k\tau}, W^-_{k\tau}; \K)
\end{eqnarray}
is an isomorphism.  Specially, $(\psi^1)_\ast=id$, and
$(\psi^k)_\ast\circ(\psi^l)_\ast=(\psi^{kl})_\ast$ if the iterations
$\gamma^l$ and $\gamma^{kl}$ are also isolated, and
\begin{equation}\label{e:4.10}
\left.\begin{array}{ll}
&m^-_{kl\tau}(\gamma^{kl})=m^-_{l\tau}(\gamma^l)=m^-_{\tau}(\gamma),\\
&m^0_{kl\tau}(\gamma^{kl})= m^0_{l\tau}(\gamma^l)=m^0_{\tau}(\gamma)
\end{array}\right\}.
\end{equation}
\end{theorem}

When $M=\R^n$, this theorem was proved by \cite[Th.3.7]{Lo2}. We
shall reduce the proof of Theorem~\ref{th:4.4} to that case.

 Using the chart in (\ref{e:4.4}) let $\tilde\gamma=(\phi_{\tau})^{-1}(\gamma)$.
Then $\tilde\gamma^j=(\phi_{j\tau})^{-1}(\gamma^j)$ for each
$j\in\N$. Then $\tilde\gamma^j$ are isolated critical points of
$\tilde{\cal L}_{j\tau}={\cal L}_{j\tau}\circ\phi_{j\tau}$ in
$W^{1,2}(S_{j\tau}, \R^n)$, $j=1,k, l, kl$. Moreover,
$m^-_{j\tau}(\tilde\gamma^j)=m^-_{\tau}(\tilde\gamma)$ and
  $m^0_{k\tau}(\tilde\gamma^j)=m^0_{\tau}(\tilde\gamma)$ for $j=k,l,
  kl$. Let
$\bigl(\tilde W(\tilde\gamma), \tilde W(\tilde\gamma)^-\bigr)$ and
$\bigl(\tilde W(\tilde\gamma^k), \tilde W(\tilde\gamma^k)^-\bigr)$
be
 Gromoll-Meyer pairs of $\tilde{\cal L}_\tau$ at $\tilde\gamma$ and
 of $\tilde{\cal L}_{k\tau}$ at $\tilde\gamma^k$, satisfying (\ref{e:4.5}).
Define
\begin{eqnarray*}
&&C_\ast(\tilde {\cal L}_{\tau}, \tilde\gamma;
\K)=H_\ast\bigl(\tilde
W(\tilde \gamma),\tilde W(\tilde\gamma)^-; \K\bigr),\\
&&C_\ast({\cal L}_{\tau}, \gamma; \K)=H_\ast\bigl(W(\gamma),
W(\gamma)^-; \K\bigr),\\
&&C_\ast(\tilde {\cal L}_{k\tau}, \tilde\gamma^k;
\K)=H_\ast\bigl(\tilde
W(\tilde \gamma^k),\tilde W(\tilde\gamma^k)^-; \K\bigr),\\
&&C_\ast({\cal L}_{k\tau}, \gamma^k; \K)=H_\ast\bigl(
W(\gamma^k),\tilde W(\gamma^k)^-; \K\bigr).
\end{eqnarray*}
Since $\phi_{k\tau}\circ{\tilde\psi}^k=\psi^k\circ\phi_\tau$, we
have
$(\phi_{k\tau})_\ast\circ({\tilde\psi}^k)_\ast=(\psi^k)_\ast\circ(\phi_\tau)_\ast$.
Clearly,
\begin{eqnarray*}
&&(\phi_\tau)_\ast:C_\ast(\tilde {\cal L}_{\tau}, \tilde\gamma;
\K)\to
C_\ast({\cal L}_{\tau}, \gamma; \K)\qquad{\rm and}\\
&&(\phi_{k\tau})_\ast:C_\ast(\tilde {\cal L}_{k\tau},
\tilde\gamma^k; \K)\to C_\ast({\cal L}_{k\tau}, \gamma^k; \K)
\end{eqnarray*}
are isomorphisms. Hence we only need to prove that
\begin{equation}\label{e:4.11}
(\tilde\psi^k)_*:  C_\ast(\tilde {\cal L}_{\tau}, \tilde\gamma;
\K)\longrightarrow
         C_\ast(\tilde {\cal L}_{k\tau}, \tilde\gamma^k; \K)
\end{equation}
is an isomorphism which maps  generators to the generators. This is
exactly one proved by \cite[Th.3.7]{Lo2}. Theorem 3.7 in \cite{Lo2}
also gives that  $(\tilde\psi^1)_*=id$ and $(\tilde\psi^k)_*\circ
(\tilde\psi^l)_*=(\tilde\psi^{kl})_*$. So other conclusions follow
immediately.

For later conveniences we outline the arguments therein.  Let
\begin{eqnarray*}
W^{1,2}(S_{k\tau}, \R^n)\!\!\!&=&\!\!\!M^0(\tilde\gamma_k)\oplus
M(\tilde\gamma_k)^-\oplus M(\tilde\gamma)^+\\
&=&\!\!\!M^0(\tilde\gamma_k)\oplus M(\tilde\gamma_k)^\bot
\end{eqnarray*}
be the orthogonal decomposition of  the space $W^{1,2}(S_{k\tau},
\R^n)$  according to the null, negative, and positive definiteness
of the quadratic form $\tilde{\cal L}_{k\tau}''(\tilde\gamma^k)$.
The generalized Morse lemma (\cite[Th.5.1, pp. 44]{Ch} yields a
homeomorphism $\tilde\Theta_{k\tau}$ from some open neighborhood
$\tilde U_{k\tau}$ of $0$ in $W^{1,2}(S_{k\tau}, \R^n)$ to
$\tilde\Theta_{k\tau}(\tilde U_{k\tau})\subset W^{1,2}(S_{k\tau},
\R^n)$ with $\tilde\Theta_{k\tau}(0)=\tilde\gamma^k$, and a map
$\tilde h_{k\tau}\in C^1\left(\tilde U_{k\tau}\cap
M(\tilde\gamma^k)^0, M(\tilde\gamma^k)^\bot\right)$ such that
\begin{eqnarray}\label{e:4.12}
\tilde {\cal
L}_{k\tau}(\tilde\Theta_{k\tau}(\eta+\xi))\!\!\!&=&\!\!\!\tilde
{\cal L}_{k\tau}\bigl(\tilde\gamma^k+ \eta+ \tilde
h_{k\tau}(\eta)\bigr)
+\displaystyle\frac{1}{2}\bigl(\tilde {\cal L}_{k\tau}''(\tilde\gamma^k)\xi,\xi\bigr)\nonumber\\
&\equiv &\!\!\!\tilde\alpha_{k\tau}(\eta)+ \tilde\beta_{k\tau}(\xi)
\end{eqnarray}
for any $\eta+\xi\in \tilde U_{k\tau}\cap
\left(M(\tilde\gamma^k)^0\oplus M(\tilde\gamma_k)^\bot\right)$.
({\bf Note}: $\tilde\beta_{k\tau}$ is $C^\infty$,
$\tilde\alpha_{k\tau}$ is $C^2$ as noted at the end of proof of
Theorem~5.1 on the page 44 of \cite{Ch}. Carefully checking the
beginning  proof therein one can easily derive this from
Lemma~\ref{lem:4.2}). It is easy to prove that
\begin{equation}\label{e:4.13}
 \tilde\psi^k\bigl(\tilde{\cal L}^{\prime}_{\tau}(x)\bigr)= \tilde {\cal
L}^{\prime}_{k\tau}\bigl(\tilde\psi^k(x)\bigr)\quad{\rm
and}\quad\tilde\psi^k\bigl(\tilde {\cal L}^{\prime\prime}_{\tau}(
x)\xi\bigr)= \tilde{\cal
L}^{\prime\prime}_{k\tau}\bigl(\tilde\psi^k(x)\bigr)\tilde\psi^k(\xi)
\end{equation}
for any $\tau, k\in \N$, $x\in W^{1,2}\bigl(S_\tau,
B^n_\rho(0)\bigr)$ and $\xi\in W^{1,2}(S_\tau, \R^n)$, and that
\begin{equation}\label{e:4.14}
 \tilde\alpha_{k\tau}(\tilde\psi^k(\eta))= k\tilde\alpha(\eta)\quad{\rm
and}\quad\tilde\beta_{k\tau}(\tilde\psi^k(\xi))=
k\tilde\beta_\tau(\xi)
\end{equation}
for any $\eta\in\tilde U_\tau\cap M^0(\tilde\gamma)$ and
$\xi\in\tilde U_\tau\cap M^\bot(\tilde\gamma)$.

\begin{lemma}\label{lem:4.5}\!\!\!{\rm (\cite[Lem. 3.2,
3.3]{Lo2})}\quad The iteration map $\tilde\psi^k:
M^\ast(\tilde\gamma)\to M^\ast(\tilde\gamma^k)$ for $\ast=0, -, +$
is linear, continuous and injective. If
$m^-_{k\tau}(\tilde\gamma^k)=m^-_{\tau}(\tilde\gamma)$, the map
$\tilde\psi^k:M^-(\tilde\gamma)\to M^-(\tilde\gamma^k)$ is a linear
diffeomorphism. If
$m^0_{k\tau}(\tilde\gamma^k)=m^0_{\tau}(\tilde\gamma)$, then the map
$\tilde\psi^k:M^0(\tilde\gamma)\to M^0(\tilde\gamma^k)$ is a linear
diffeomorphism, and $\tilde U_{k\tau}$, the homeomorphism
$\tilde\Theta_{k\tau}$ and map
 $\tilde h_{k\tau}\in C^1\left(\tilde U_{k\tau}\cap M(\tilde\gamma^k)^0,
M(\tilde\gamma^k)^\bot\right)$ are chosen to satisfy:
\begin{eqnarray}
&& \tilde U_{k\tau}\cap\tilde\psi^k(W^{1,2}(S_\tau,
\R^n))=\tilde\psi^k(\tilde U_{\tau}),\label{e:4.15}\\
&&\tilde\Theta_{k\tau}\circ\tilde\psi^k=\tilde\psi^k\circ\tilde\Theta_\tau:\tilde
U_\tau\to \tilde\Theta_\tau(\tilde U_\tau\cap
M^0(\tilde\gamma^k)),\label{e:4.16}\\
&& \tilde h_{k\tau}(\tilde\psi^k(\eta))=\tilde\psi^k(\tilde
h_\tau(\eta))\quad\forall \eta\in\tilde U_\tau\cap
M(\tilde\gamma).\label{e:4.17}
\end{eqnarray}
\end{lemma}

Let $(W_0, W_0^-)$ and $(W_1, W_1^-)$ be Gromoll-Meyer pairs of
$\tilde\alpha_\tau$ and $\tilde\beta_\tau$ at their origins
respectively. By \cite[Prop.3.5. $2^\circ$]{Lo2},
$(\tilde\psi^k(W_0), \tilde\psi^k(W_0^-))$ is a Gromoll-Meyer pair
of $\tilde\alpha_{k\tau}$ at the origin. The Gromoll-Meyer pair
$(W_1, W_1^-)$ can also be chosen to satisfy
\begin{equation}\label{e:4.18}
(\tilde\psi^k(W_1), \tilde\psi^k(W_1^-))\subset (V, V^-)
\end{equation}
for some Gromoll-Meyer pair $(V, V^-)$ of $\tilde\beta_{k\tau}$ at
the origin. By \cite[Lem.5.1. pp.51]{Ch}
\begin{eqnarray}
&&\left( W_0\times W_1, (W_0\times W_1^-)\cup (W_0^-\times
W_1)\right),\label{e:4.19}\\
&& \left( \tilde\psi^k(W_0)\times V, (\tilde\psi^k(W_0)\times
V^-)\cup (\tilde\psi^k(W_0^-)\times V)\right)\label{e:4.20}
\end{eqnarray}
 are Gromoll-Meyer pairs of
$\tilde\alpha_\tau+ \tilde\beta_\tau$ and $\tilde\alpha_{k\tau}+
\tilde\beta_{k\tau}$ at their origins respectively, and also satisfy
\begin{equation}\label{e:4.21}
\left.\begin{array}{ll}
 &\left( \tilde\psi^k(W_0\times W_1), \tilde\psi^k((W_0\times
W_1^-)\cup (W_0^-\times
W_1))\right)\subset\vspace{2mm}\\
& \left( \tilde\psi^k(W_0)\times V, (\tilde\psi^k(W_0)\times
V^-)\cup (\tilde\psi^k(W_0^-)\times V)\right).
\end{array}\right\}
\end{equation}
Note that
\begin{eqnarray}
&(\widehat W_\tau, \widehat W^-_\tau):=\tilde\Theta_\tau\left(
W_0\times W_1, (W_0\times W_1^-)\cup (W_0^-\times
W_1)\right),\label{e:4.22}\\
& (\widehat W_{k\tau}, \widehat
W^-_{k\tau}):=\tilde\Theta_{k\tau}\left( \tilde\psi^k(W_0)\times V,
(\tilde\psi^k(W_0)\times V^-)\cup (\tilde\psi^k(W_0^-)\times
V)\right)\label{e:4.23}
\end{eqnarray}
 are topological Gromoll-Meyer pairs of
$\tilde{\cal L}_\tau$ at $\tilde\gamma$ and $\tilde{\cal L}_{k\tau}$
at $\tilde\gamma^k$ respectively. Let
\begin{eqnarray*}
&&C_\ast(\tilde\alpha_{\tau}+ \tilde\beta_{\tau}, 0; \K):=H_\ast
(W_0\times W_1, (W_0\times W_1^-)\cup (W_0^-\times W_1);\K),\\
&&C_\ast(\tilde{\cal L}_{\tau}, 0; \K):=H_\ast(\widehat W_\tau,
\widehat W^-_\tau;\K),\\
&&C_\ast(\tilde\alpha_{k\tau}+ \tilde\beta_{k\tau}, 0;
\K):=H_\ast\left( \tilde\psi^k(W_0)\times V,
(\tilde\psi^k(W_0)\times V^-)\cup
(\tilde\psi^k(W_0^-)\times V);\K\right),\\
&&C_\ast(\tilde{\cal L}_{k\tau}, 0; \K):=H_\ast(\widehat W_{k\tau},
\widehat W^-_{k\tau};\K).
\end{eqnarray*}
We have the isomorphisms on critical modules,
\begin{eqnarray*}
&&(\tilde\Theta_{\tau})_\ast:C_\ast(\tilde\alpha_{\tau}+
\tilde\beta_{\tau}, 0; \K)\cong C_\ast(\tilde{\cal L}_{\tau},
\tilde\gamma; \K),\\
&&(\tilde\Theta_{k\tau})_\ast:C_\ast(\tilde\alpha_{k\tau}+
\tilde\beta_{k\tau}, 0; \K)\cong C_\ast(\tilde{\cal L}_{k\tau},
\tilde\gamma^k; \K).
\end{eqnarray*}
By (\ref{e:4.21})  we have a homomorphism
\begin{equation}\label{e:4.24}
(\tilde\psi^k)_\ast:C_\ast(\tilde\alpha_{\tau}+ \tilde\beta_{\tau},
0; \K)\to C_\ast(\tilde\alpha_{k\tau}+ \tilde\beta_{k\tau}, 0; \K).
\end{equation}
Moreover, (\ref{e:4.16}) and (\ref{e:4.21}) show that
\begin{equation}\label{e:4.25}
(\tilde\psi^k(\widehat W_\tau),
\tilde\psi^k(\widehat W^-_\tau))\subset (\widehat W_{k\tau},
\widehat W^-_{k\tau})
\end{equation}
and therefore the homomorphism
$$
(\tilde\psi^k)_\ast:C_\ast(\tilde{\cal L}_{\tau}, 0; \K)\to
C_\ast(\tilde{\cal L}_{k\tau}, 0; \K)
$$
satisfy
\begin{equation}\label{e:4.26}
(\tilde\psi^k)_\ast\circ(\tilde\Theta_{\tau})_\ast=(\tilde\Theta_{k\tau})_\ast\circ(\tilde\psi^k)_\ast.
\end{equation}
Hence the problem is reduced to  prove:

\begin{lemma}\label{lem:4.6}
The Gromoll-Meyer pairs $(W_1, W_1^-)$ and $(V, V^-)$  in
(\ref{e:4.18}) can be chosen such that
\begin{equation}\label{e:4.27}
(\tilde\psi^k)_\ast:C_\ast(\tilde\alpha_{\tau}+ \tilde\beta_{\tau},
0; \K)\to C_\ast(\tilde\alpha_{k\tau}+ \tilde\beta_{k\tau}, 0; \K)
\end{equation}
is an isomorphism.
\end{lemma}

\noindent{\bf Proof}.\quad For $j=1, k$, decompose $\xi\in
M(\tilde\gamma_j)^\bot=M(\tilde\gamma_j)^-\oplus
M(\tilde\gamma_j)^+$ into $\xi=\xi^- + \xi^+$ and write
$$
\tilde\beta_{j\tau}(\xi)=\tilde\beta_{j\tau}(\xi^-)+
\tilde\beta_{j\tau}(\xi^+)= \tilde\beta_{j\tau}^-(\xi^-)+
\tilde\beta_{j\tau}^+(\xi^+).
$$
Then $\tilde\beta_{j\tau}^-$ and $\tilde\beta_{j\tau}^+$ are
negative and positive definite quadratic forms on
$M(\tilde\gamma_j)^-$ and $M(\tilde\gamma_j)^+$ with Morse indexes
$m^-(\tilde\gamma^j)$ and $0$ respectively, $j=1,k$. The
(\ref{e:4.12})-(\ref{e:4.14}) imply
$$
\tilde\beta^-_{k\tau}(\tilde\psi^k(\xi^-))=
k\tilde\beta^-_\tau(\xi^-)\quad{\rm and}\quad
\tilde\beta^-_{k\tau}(\tilde\psi^k(\xi^+))=
k\tilde\beta^-_\tau(\xi^+)
$$
for any $\xi^-\in M^-(\tilde\gamma)$ and $\xi^+\in
M^+(\tilde\gamma)$. Since
$m^-_{k\tau}(\tilde\gamma^k)=m^-_{\tau}(\tilde\gamma)$, by
Lemma~\ref{lem:4.5} the map $\tilde\psi^k:M^-(\tilde\gamma)\to
M^-(\tilde\gamma^k)$ is a linear diffeomorphism. Let $(W_{11},
W_{11}^-)$ be a Gromoll-Meyer pair of $\tilde\beta^-_\tau$ at the
origin. Then
\begin{equation}\label{e:4.28}
(\tilde\psi^k(W_{11}), \tilde\psi^k(W^-_{11}))
\end{equation}
is a Gromoll-Meyer pair of $\tilde\beta^-_{k\tau}$  at the origin.
For $\delta>0$ sufficiently small, set
\begin{eqnarray*}
&&W_{12}:=\{\xi^+\in M(\tilde\gamma)^+\,|\,
\|\xi^+\|_\tau\le\delta\,\},\\
&& W^-_{12}:=\{\xi^+\in
M(\tilde\gamma)^+\,|\, \|\xi^+\|_\tau=\delta\,\},\\
&&V_{12}:=\{\xi^+\in M(\tilde\gamma^k)^+\,|\, \|\xi^+\|_{k\tau}\le
\sqrt{k}\delta\,\},\\
&& V^-_{12}:=\{\xi^+\in M(\tilde\gamma^k)^+\,|\,
\|\xi^+\|_{k\tau}=\sqrt{k}\delta\,\}.
\end{eqnarray*}
 It is easily checked that $(W_{12}, W^-_{12})$ and
$(V_{12}, V^-_{12})$ are Gromoll-Meyer pairs of
$\tilde\beta^+_{\tau}$ and $\tilde\beta^+_{k\tau}$ at their origins
respectively, and that
\begin{equation}\label{e:4.29}
\left(\tilde\psi^k(W_{12}), \tilde\psi^k(W^-_{12})\right)\subset
(V_{12}, V^-_{12}).
\end{equation}
By \cite[Lem.5.1. pp.51]{Ch}, we may take
\begin{eqnarray}
&&(W_1, W_1^-):=\left( W_{11}\times W_{12}, (W_{11}\times
W_{12}^-)\cup (W_{11}^-\times
W_{12})\right),\label{e:4.30}\\
&& (V, V^-):=\left( \tilde\psi^k(W_{11})\times V_{12},
(\tilde\psi^k(W_{11})\times V_{12}^-)\cup
(\tilde\psi^k(W_{11}^-)\times V_{12})\right).\label{e:4.31}
\end{eqnarray}
Then $(W_0\times W_1, (W_0\times W_1^-)\cup (W_0^-\times W_1))$
becomes $(W,  W^-)$, and
\begin{equation}\label{e:4.32}
C_\ast(\tilde\alpha_{\tau}+ \tilde\beta_{\tau}, 0; \K)=H_\ast(W,
W^-;\K),
\end{equation}
where $W:=W_0\times W_{11}\times W_{12}$ and
\begin{equation}\label{e:4.33}
 W^-:=\bigl(W_0\times (W_{11}\times W_{12}^-)\cup (W_{11}^-\times
W_{12})\bigr)\cup (W_0^-\times W_{11}\times W_{12}).
\end{equation}
Moreover, $\bigl( \tilde\psi^k(W_0)\times V,
(\tilde\psi^k(W_0)\times V^-)\cup (\tilde\psi^k(W_0^-)\times
V)\bigr)$ becomes $(U, U^-)$, and
\begin{equation}\label{e:4.34}
C_\ast(\tilde\alpha_{k\tau}+ \tilde\beta_{k\tau}, 0; \K)=H_\ast(U,
U^-;\K),
\end{equation}
where $U=\tilde\psi^k(W_0)\times  \tilde\psi^k(W_{11})\times V_{12}$
and
\begin{eqnarray}\label{e:4.35}
&&U^-=\bigl(\tilde\psi^k(W_0)\times (\tilde\psi^k(W_{11})\times
V_{12}^-)\cup (\tilde\psi^k(W_{11}^-)\times
V_{12})\bigr)\nonumber\\
 &&\hspace{40mm}\cup \bigl(\tilde\psi^k(W_0^-)\times
\tilde\psi^k(W_{11})\times V_{12}\bigr).
\end{eqnarray}
Note that $\tilde\psi^k(W)=\tilde\psi^k(W_0)\times
\tilde\psi^k(W_{11})\times \tilde\psi^k(W_{12})$ and
\begin{eqnarray}\label{e:4.36}
&&\tilde\psi^k(W^-)=\bigl(\tilde\psi^k(W_0)\times
(\tilde\psi^k(W_{11})\times \tilde\psi^k(W_{12}^-))\cup
(\tilde\psi^k(W_{11}^-)\times
\tilde\psi^k(W_{12}))\bigr)\nonumber\\
 &&\hspace{40mm}\cup \bigl(\tilde\psi^k(W_0^-)\times
\tilde\psi^k(W_{11})\times \tilde\psi^k(W_{12})\bigr).
\end{eqnarray}
Since $\tilde\psi^k: M^+(\tilde\gamma)\to M^+(\tilde\gamma^k)$ is a
linear, continuous and injection, by (\ref{e:4.29}) and the
constructions of $(V_{12}, V_{12}^-)$ and $(W_{12}, W_{12}^-)$ it is
readily checked that $\bigl(\tilde\psi^k(W_{12}),
\tilde\psi^k(W^-_{12})\bigr)$ is a deformation retract of $(V_{12},
V^-_{12})$. It follows that
$$
\bigl(\tilde\psi^k(W), \tilde\psi^k(W^-)\bigr)\subset (U, U^-)
$$
is a deformation retract of $(U, U^-)$. Hence
$$
(\tilde\psi^k)_\ast: H_\ast(W, W^-;\K\bigr)\to H_\ast(U, U^-;\K)
$$
and therefore, by (\ref{e:4.32}) and (\ref{e:4.34}), the
homomorphism $(\tilde\psi^k)_\ast$ in (\ref{e:4.27}) is an
isomorphism.

\textsf{We may also prove the conclusion as follows}.  By the
arguments at the middle of \cite[pp. 51]{Ch} we can use K\"{u}nneth
formula to arrive
\begin{eqnarray}
&&C_\ast(\tilde\alpha_{\tau}+ \tilde\beta_{\tau}, 0; \K)=H_\ast(W_0,
W_0^-;\K)\otimes\nonumber\\
&&\qquad\qquad H_\ast(W_{11}, W_{11}^-;\K)\otimes H_\ast(W_{12}, W_{12}^-;\K),\label{e:4.37}\\
&&C_\ast(\tilde\alpha_{k\tau}+ \tilde\beta_{k\tau}, 0; \K)=H_\ast(
\tilde\psi^k(W_0),
 \tilde\psi^k(W_0^-);\K)\otimes\nonumber\\
 &&\qquad\qquad H_\ast(
\tilde\psi^k(W_{11}),
 \tilde\psi^k(W_{11}^-);\K)\otimes H_\ast(V_{12}, V_{12}^-;\K).\label{e:4.38}
\end{eqnarray}
Now $m^-_{k\tau}(\tilde\gamma^k)=m^-_{\tau}(\tilde\gamma)$ and
$m^0_{k\tau}(\tilde\gamma^k)=m^0_{\tau}(\tilde\gamma)$ imply that
\begin{eqnarray*}
&&(\tilde\psi^k)_\ast:H_\ast(W_0, W_0^-;\K)\to H_\ast(
\tilde\psi^k(W_0),
 \tilde\psi^k(W_0^-);\K),\\
&&(\tilde\psi^k)_\ast:H_\ast(W_{11}, W_{11}^-;\K)\to H_\ast(
\tilde\psi^k(W_{11}),
 \tilde\psi^k(W_{11}^-);\K)
\end{eqnarray*}
are isomorphisms. Since  $\bigl(\tilde\psi^k(W_{12}),
\tilde\psi^k(W^-_{12})\bigr)$ is a deformation retract of $(V_{12},
V^-_{12})$ as above, it follows that
\begin{equation}\label{e:4.39}
(\tilde\psi^k)_\ast:H_\ast\bigl(\tilde\psi^k(W_{12}),
\tilde\psi^k(W^-_{12});\K\bigr)\to H_\ast(V_{12}, V^-_{12};\K)
\end{equation}
is an isomorphism. By (\ref{e:4.37}) and (\ref{e:4.38}) we get the
proof of Lemma~\ref{lem:4.6}. $\Box$\vspace{2mm}

For $(\widehat W_\tau, \widehat W^-_\tau)$ in (\ref{e:4.22}) and
$(\widehat W_{k\tau}, \widehat W^-_{k\tau})$ in (\ref{e:4.23}),
where the Gromoll-Meyer pairs $(W_1, W_1^-)$ and $(V, V^-)$  in
(\ref{e:4.18}) are also required to satisfy Lemma~\ref{lem:4.6}. Set
$$
(W_\tau,  W^-_\tau):=\bigl(\phi_\tau(\widehat W_\tau),
\phi_\tau(\widehat W^-_\tau)\bigr)\quad{\rm and}\quad
 (W_{k\tau}, W^-_{k\tau}):=\bigl(\phi_{k\tau}(\widehat W_{k\tau}),
\phi_{k\tau}(\widehat W^-_{k\tau})\bigr).
$$
Since $\phi_{k\tau}\circ{\tilde\psi}^k=\psi^k\circ\phi_\tau$, by
(\ref{e:4.25}) we have $(\psi^k(W_\tau), \psi^k(W^-_\tau))\subset
(W_{k\tau}, W^-_{k\tau})$ and that the homomorphism
$$
(\psi^k)_\ast: H_\ast(W_\tau, W^-_\tau; \K)\to H_\ast(W_{k\tau},
W^-_{k\tau}; \K)
$$
is an isomorphism. Consequently, $(W_\tau,  W^-_\tau)$ and
$(W_{k\tau}, W^-_{k\tau})$  are desired topological Gromoll-Meyer
pairs.

The other conclusions are also easily proved. So
Theorem~\ref{th:4.4} holds. $\Box$\vspace{2mm}

\noindent{\bf 4.2.} In this subsection we always assume: $M$ is
$C^3$-smooth, $L$ is $C^2$-smooth and satisfies (L1)-(L4).  Let
$\gamma\in EH_\tau$ be an isolated critical point of the functional
${\cal L}^E_\tau$ on $EH_\tau$, and
\begin{equation}\label{e:4.40}
\phi^E_{k\tau}: EW^{1,2}(S_{k\tau}, B^n_\rho(0))\to EH_{k\tau}\quad
{\rm and}\quad \tilde{\cal L}^E_{k\tau}={\cal
L}_{k\tau}\circ\phi_{k\tau}
\end{equation}
be as in (\ref{e:3.29}) and (\ref{e:3.32})  for each $k\in \N$. They
satisfy (\ref{e:3.30}), i.e.
$\phi^E_{k\tau}\circ{\tilde\psi}^k=\psi^k\circ\phi^E_\tau$ for any
$k\in\N$, where $\psi^k: EH_\tau\to EH_{k\tau}$ and
$$
\tilde\psi^k: EW^{1,2}(S_{\tau}, \R^n)\to EW^{1,2}(S_{k\tau},
\R^n)
$$
are the iteration maps. Let
$\tilde\gamma=(\phi^E_\tau)^{-1}(\gamma)$ and thus
$\phi^E_{k\tau}(\tilde\gamma^k)=\gamma^k$ for any $k\in\N$.
 Suppose
that $\gamma^k$ and therefore $\tilde\gamma^k$ are also isolated.
Denote by
$$
C_q(\tilde {\cal L}^E_{k\tau}, \tilde\gamma^k;\K)=H_q\bigl(\tilde
W(\tilde \gamma^k)_E,\tilde W(\tilde\gamma^k)_E^-; \K\bigr)
$$
  the  critical module of $\tilde {\cal L}^E_{k\tau}$ at  $\tilde\gamma^k$
 via the relative singular homology with coefficients in $\K$,
  where $\bigl(\tilde W(\tilde\gamma^k)_E, \tilde W(\tilde\gamma^k)^-_E\bigr)$ is
a Gromoll-Meyer pair via some pseduo-gradient vector field of
$\tilde {\cal L}^E_{k\tau}$ near $\tilde\gamma^k$ in
$EW^{1,2}(S_{k\tau}, \R^n)$. Let
\begin{eqnarray*}
EW^{1,2}(S_{k\tau}, \R^n)\!\!\!&=&\!\!\!M^0(\tilde\gamma_k)_E\oplus
M(\tilde\gamma_k)^-_E\oplus M(\tilde\gamma)^+_E\\
&=&\!\!\!M^0(\tilde\gamma_k)_E\oplus M(\tilde\gamma_k)^\bot_E
\end{eqnarray*}
be the orthogonal decomposition of  the space $EW^{1,2}(S_{k\tau},
\R^n)$  according to the null, negative, and positive definiteness
of the quadratic form $(\tilde{\cal L}^E_{k\tau})''(\tilde\gamma)$.
As above we can use the generalized Morse lemma to get a
homeomorphism $\tilde\Theta^E_{k\tau}$ from some open neighborhood
$\tilde U^E_{k\tau}$ of $0$ in $EW^{1,2}(S_{k\tau}, \R^n)$ to
$\tilde\Theta^E_{k\tau}(\tilde U^E_{k\tau})\subset
EW^{1,2}(S_{k\tau}, \R^n)$ with
$\tilde\Theta^E_{k\tau}(0)=\tilde\gamma^k$, and a map $\tilde
h^E_{k\tau}\in C^1\bigl(\tilde U^E_{k\tau}\cap
M(\tilde\gamma^k)^0_E, M(\tilde\gamma^k)_E^\bot\bigr)$ such that
$$
\begin{array}{rcl}
\tilde {\cal
L}^E_{k\tau}\bigl(\tilde\Theta^E_{k\tau}(\eta+\xi)\bigr)\!\!\!&=&\!\!\!\tilde
{\cal L}^E_{k\tau}\bigl(\tilde\gamma^k+\eta+ \tilde
h^E_{k\tau}(\eta)\bigr)
+\displaystyle\frac{1}{2}\bigl((\tilde {\cal L}^E_{k\tau})''(\tilde\gamma^k)\xi,\xi\bigr)\\
&\equiv &\!\!\!\tilde\alpha^E_{k\tau}(\eta)+
\tilde\beta^E_{k\tau}(\xi)
\end{array}
$$
for any $\eta+\xi\in \tilde U^E_{k\tau}\cap
(M(\tilde\gamma^k)_E^0\oplus M(\tilde\gamma_k)_E^\bot)$, where
$\tilde\beta^E_{k\tau}$ and $\tilde\alpha^E_{k\tau}$ are
respectively $C^\infty$ and $C^2$ as noted below (\ref{e:4.12}).
Then $\tilde\Theta^E_{k\tau}$ induces isomorphisms on critical
modules,
\begin{equation}\label{e:4.41}
(\tilde\Theta^E_{k\tau})_\ast:C_\ast(\tilde\alpha^E_{k\tau}+
\tilde\beta^E_{k\tau}, 0;\K)\cong C_\ast(\tilde{\cal L}^E_{k\tau},
\tilde\gamma^k;\K).
\end{equation}
 Note that
\begin{equation}\label{e:4.42}
(W(\gamma^k)_E, W^-(\gamma^k)_E):=\left(\phi^E_{k\tau}\bigl(\tilde
W(\tilde\gamma^k)_E\bigr), \phi^E_{k\tau}\bigl(\tilde
W^-(\gamma^k)_E\bigr)\right)
\end{equation}
 is  a Gromoll-Meyer pair of ${\cal L}^E_{k\tau}$ at $\gamma^k$.
 Define the  critical modules
\begin{equation}\label{e:4.43}
 C_\ast({\cal L}^E_{k\tau}, \gamma^k; \K):=H_\ast\bigl(W(\gamma^k)_E, W^-(\gamma^k)_E;
 \K\bigr).
\end{equation}
Then corresponding to Theorem~\ref{th:4.4} we have the following
generalization of \cite[Lemma 4.1]{LuW2}.

\begin{theorem}\label{th:4.7}
Let $\gamma\in EH_\tau$ be an isolated critical point of the
functional ${\cal L}^E_\tau$ on $EH_\tau$. If  the iteration
 $\gamma^k$ is also isolated  for some $k\in\N$, and
 $$m^-_{1, k\tau}(\gamma^k)=m^-_{1, \tau}(\gamma)\quad{\rm and}\quad
  m^0_{1, k\tau}(\gamma^k)=m^0_{1, \tau}(\gamma),
  $$
then for $c={\cal L}^E_\tau(\gamma)$ and any $\epsilon>0$ there
exist topological Gromoll-Meyer pairs of ${\cal L}^E_\tau$ at
$\gamma\in EH_\tau$ and
 of ${\cal L}^E_{k\tau}$ at $\gamma^k\in EH_{k\tau}$,
\begin{eqnarray*}
&&(W_\tau,  W^-_\tau)\subset\bigl(({\cal L}^E_\tau)^{-1}[c-\epsilon,
c+ \epsilon], ({\cal L}^E_\tau)^{-1}(c-\epsilon)\bigr)
 \quad\hbox{and}\\
&&(W_{k\tau}, W^-_{k\tau})\subset\bigl(({\cal
L}^E_{k\tau})^{-1}[kc-k\epsilon, kc+ k\epsilon], ({\cal
L}^E_{k\tau})^{-1}(kc-k\epsilon)\bigr),
\end{eqnarray*}
such that
\begin{equation}\label{e:4.44}
 (\psi^k(W_\tau), \psi^k(W^-_\tau))\subset
(W_{k\tau}, W^-_{k\tau})
\end{equation}
and that  the iteration map $\psi^k: EH_\tau\to EH_{k\tau}$ induces
isomorphisms
\begin{eqnarray}\label{e:4.45}
&&(\psi^k)_\ast: C_\ast({\cal L}^E_\tau, \gamma;\K):=H_\ast(W_\tau,
W^-_\tau; \K)\nonumber\\
&&\hspace{40mm}\to C_\ast({\cal L}^E_{k\tau},
\gamma^k;\K):=H_\ast(W_{k\tau}, W^-_{k\tau}; \K).
\end{eqnarray}
 Specially, $(\psi^1)_\ast=id$,
and $(\psi^k)_\ast\circ(\psi^l)_\ast=(\psi^{kl})_\ast$ if the
iterations $\gamma^l$ and $\gamma^{kl}$ are also isolated, and
$$
\left.\begin{array}{ll} &m^-_{1, kl\tau}(\gamma^{kl})=m^-_{1,
l\tau}(\gamma^l)=m^-_{1, \tau}(\gamma),\\
& m^0_{1, kl\tau}(\gamma^{kl})= m^0_{1, l\tau}(\gamma^l)=m^0_{1,
\tau}(\gamma)
\end{array}\right\}.
$$
 \end{theorem}

\noindent{\bf 4.3.}  Let us consider the case  $L$ is independent
$t$. In this subsection we always assume: $M$ is $C^5$-smooth, $L$
is $C^4$-smooth and satisfies (L1)-(L3). The goal is to generalize
\cite[Th.2.5]{LoLu} to the present general case. However, unlike the
last two cases we cannot choose a local coordinate chart around a
critical orbit. For $\tau>0$, let $
S_\tau:=\R/\tau\Z=\{[s]_{\tau}\,|\,[s]_{\tau}=s+\tau\Z,\, s\in\R\}$
 and the functional ${\cal
 L}_\tau:H_\tau(\alpha)\to\R$ be still defined by (\ref{e:1.16}).
 By \cite[Chp.2, \S2.2]{Kl},  there
exist  equivariant and also isometric operations of
$S_{\tau}$-action on $H_{\tau}(\alpha)$ and $TH_\tau(\alpha)$:
\begin{equation}\label{e:4.46}
\left.\begin{array}{ll} &[s]_{\tau}\cdot \gamma(t)=\gamma(s+t),
\quad \forall [s]_{\tau}\in S_{\tau},
         \; \gamma\in H_{\tau}(\alpha),\\
 &[s]_{\tau}\cdot \xi(t)=\xi(s+t),
\quad \forall [s]_{\tau}\in S_{\tau},
         \; \xi\in T_{\gamma}H_{\tau}(\alpha)
 \end{array}\right\}
         \end{equation}
which are  continuous, but not differentiable. Clearly, ${\cal
L}_{\tau}$ is invariant under this action. Since under our
assumptions each critical point $\gamma$ of ${\cal L}_\tau$ is
$C^4$-smooth, by \cite[p. 499]{GM2}, the orbit $S_{\tau}\cdot
\gamma$ is a $C^3$-submanifold in $H_{\tau}(\alpha)$.
 It is easily checked that $S_{\tau}\cdot
\gamma$ is a $C^3$-smooth critical submanifold of ${\cal L}_\tau$.
Seemingly, the theory of \cite{Wa} cannot be applied to this case
because the action of $S_\tau$ is only continuous. However, as
pointed out in the second paragraph of \cite[pp. 500]{GM2} this
theory still hold since critical orbits are smooth and $S_\tau$ acts
by isometries.

 For any
$k\in\N$, there is a natural $k$-fold cover $\varphi_k$ from
$S_{k\tau}$ to $S_{\tau}$ defined by
\begin{equation}\label{e:4.47}
\varphi_k: [s]_{k\tau}\mapsto [s]_{\tau}.
 \end{equation}
 It is easy to check that the $S_{\tau}$-action
on $H_{\tau}(\alpha)$, the $S_{k\tau}$-action on
$H_{k\tau}(\alpha^k)$, and the $k$-th iteration map $\psi^k$
defined above (\ref{e:3.9}) satisfy:
\begin{equation}\label{e:4.48}
\left.\begin{array}{rcl}
  ([s]_{\tau}\cdot \gamma)^k &=& [s]_{k\tau}\cdot \gamma^k,  \\
  {\cal L}_{k\tau}([s]_{k\tau}\cdot \gamma^k) &=& k{\cal L}_{\tau}([s]_{\tau}\cdot \gamma)
     = k{\cal L}_{\tau}(\gamma)
\end{array}\right\}
 \end{equation}
 for all $\gamma\in H_{\tau}(\alpha)$, $k\in\N$,  and $s\in \R$.

 Let $\gamma\in H_{\tau}(\alpha)$ be a \textsf{non-constant critical point} of ${\cal L}_{\tau}$ with
minimal period $\tau/m$ for some $m\in\N$. Denote by ${\cal
O}=S_{\tau}\cdot \gamma=S_{\tau/m}\cdot \gamma$. It is a
$1$-dimensional $C^3$-submanifold diffeomorphic to the circle.
 Let $c={\cal L}_\tau|_{\cal O}$. \textsf{Assume that ${\cal O}$ is isolated}. We
 may  take a neighborhood $U$ of ${\cal O}$ such that
 ${\cal K}({\cal L}_\tau)\cap U={\cal
 O}$. By (\ref{e:4.1}) we have critical group
 $C_\ast({\cal L}_\tau, {\cal O}; \K)$
of ${\cal L}_\tau$ at ${\cal O}$.   For every $s\in [0,\tau/m]$ the
tangent space $T_{s\cdot \gamma}(S_{\tau}\cdot \gamma)$ is
$\R(s\cdot \gamma)^\cdot$, and the fiber $N({\cal O})_{s\cdot
\gamma}$ at $s\cdot \gamma$ of the normal bundle $N({\cal O})$ of
${\cal O}$ is a subspace of codimension $1$ which is orthogonal to
$(s\cdot \gamma)^\cdot$ in $T_{s\cdot \gamma}H_{\tau}(\alpha)$, i.e.
$$
N({\cal O})_{s\cdot \gamma}=\left\{ \xi\in T_{s\cdot
\gamma}H_{\tau}(\alpha)\,|\, \langle\!\langle\xi, (s\cdot
\gamma)^\cdot\rangle\!\rangle_1=0\,\right\}.
$$
Since $H_\tau(\alpha)$ is $C^4$-smooth and ${\cal O}$ is a
$C^3$-smooth submanifold, $N({\cal O})$ is $C^2$-smooth manifold.
\footnote{This is the reason that we require higher smoothness of
$M$ and $L$.} Notice that $N({\cal O})$ is invariant under the
$S_{\tau}$-actions in (\ref{e:4.20}) and each $[s]_{\tau}$ gives an
isometric bundle map
\begin{equation}\label{e:4.49}
N({\cal O})\to N({\cal O}),\,(z, v)\mapsto ([s]_{\tau}\cdot z,
[s]_{\tau}\cdot v).
\end{equation}
Under the present case it is easily checked that  ${\cal L}_{\tau}$
satisfies the Assumption 7.1 on the page 71 of [Ch], that is, there
exists  $\epsilon>0$ such that
\begin{equation}\label{e:4.50}
\sigma\bigl({\cal L}''_{\tau}(x)\bigr)\cap([-\epsilon,\epsilon]
\setminus\{0\})=\emptyset,\quad\dim{\rm ker}({\cal L}''_{\tau}(x))
={\rm constant}
\end{equation}
for any $x\in{\cal O}$. Then Lemma 7.4  of \cite[pp. 71]{Ch} gives
the orthogonal $C^2$-smooth bundle decomposition
\begin{equation}\label{e:4.51}
 N({\cal O})=N({\cal O})^+\oplus N({\cal O})^-\oplus N({\cal O})^0,
\quad N({\cal O})^*=P_* N({\cal O})
\end{equation}
for $*=+,-,0$. Here $P_\ast:N({\cal O})\to N({\cal O})^\ast$,
$\ast=+, 0, -$, are orthogonal bundle projections. Each $N({\cal
O})^\ast$ is a $C^2$-smooth submanifold. It is not hard to check
that ${\cal L}'_{\tau}$ and ${\cal L}''_{\tau}$ satisfy
$$
{\cal L}_{\tau}^\prime([s]_{\tau}\cdot x)=[s]_{\tau}\cdot {\cal L}_{\tau}^\prime(x)\quad
{\rm and}\quad {\cal L}_{\tau}^{\prime\prime}([s]_{\tau}\cdot
x)([s]_{\tau}\cdot \xi)=[s]_{\tau}\cdot
      ({\cal L}_{\tau}^{\prime\prime}(x)(\xi))
      $$
for all $x\in H_{\tau}(\alpha)$, $\xi\in T_xH_\tau(\alpha)$ and
$[s]_{\tau}\in S_{\tau}$. It follows that the bundle map
(\ref{e:4.49}) preserves the decomposition (\ref{e:4.51}). In
particular, we obtain
$$ \bigl({\rm rank}N({\cal O})^-,{\rm rank}N({\cal O})^0\bigr)
    = \bigl(m^-_{\tau}(x), m^0_{\tau}(x)-1\bigr)  \quad \forall x\in{\cal O},
     $$
 where $m^-_{\tau}(x)$ and $m^0_{\tau}(x)$ are Morse index and
 nullity of ${\cal L}_\tau$ at $x$ respectively.
Define
\begin{equation}\label{e:4.52}
 \bigl(m^-_{\tau}({\cal O}), m^0_{\tau}({\cal O})\bigr):= \bigl({\rm rank}N({\cal O})^-,
    {\rm rank}N({\cal O})^0\bigr).
\end{equation}
Then
\begin{equation}\label{e:4.53}
 \bigl(m^-_{\tau}({\cal O}), m^0_{\tau}({\cal O})\bigr) =
\bigl(m^-_{\tau}(x), m^0_{\tau}(x)-1\bigr)
       \quad \forall x\in{\cal O}.
\end{equation}
For a single point critical orbit ${\cal O}=\{\gamma\}$, i.e.,
$\gamma$ is constant, we define
\begin{equation}\label{e:4.54}
\bigl(m^-_{\tau}({\cal O}), m^0_{\tau}({\cal O})\bigr): =
\bigl(m^-_{\tau}(\gamma), m^0_{\tau}(\gamma)\bigr).
\end{equation}

Note that for sufficiently small $\varepsilon>0$ the set
$$
N({\cal O})(\varepsilon):=\bigl\{(y,v)\in N({\cal O})|\,y\in
  {\cal O},\,\|v\|_1<\varepsilon\bigr\}
   $$
is contained in an open neighborhood of the zero section of the
tangent bundle $TH_\tau(\alpha)$. By \cite[Th.1.3.7, pp. 20]{Kl} we
have a $C^2$-embedding  from $N({\cal O})(\varepsilon)$ to an open
neighborhood of  the diagonal of $H_\tau(\alpha)\times
H_\tau(\alpha)$,
$$
N({\cal O})(\varepsilon)\to H_\tau(\alpha)\times H_\tau(\alpha),\;
(y, v)\mapsto (y, \exp_yv),
$$
where $\exp$ is the exponential map of the chosen Riemannian metric
on $M$ and $(\exp_yv)(t)=\exp_{y(t)}{v(t)}\;\forall t\in\R$. This
yields  a $C^2$ diffeomorphism from $N({\cal O})(\varepsilon)$ to an
open neighborhood $Q_{\varepsilon}({\cal O})$ of ${\cal O}$,
\begin{equation}\label{e:4.55}
\Psi_\tau: N({\cal O})(\varepsilon)\to Q_{\varepsilon}({\cal O}),\;
\Psi_\tau(y, v)(t)=\exp_{y(t)}{v(t)}\;\forall t\in\R,
\end{equation}
(Note that it is not the exponential map of the Levi-Civita
connection derived the Riemannian metric $\langle\!\langle\, ,\,
\rangle\!\rangle_\tau$ on $H_\tau(\alpha)$.)
 Clearly,
\begin{equation}\label{e:4.56}
\Psi_\tau(y,0)=y\;\forall y\in{\cal
O}\quad{\rm and}\quad \Psi_\tau([s]_\tau\cdot y,\, [s]_\tau\cdot
v)=[s]_\tau\cdot\Psi_\tau(y,  v)
\end{equation}
for any $(y, v)\in N({\cal O})(\varepsilon)$ and $[s]_\tau\in
S_\tau$. It follows that $Q_{\varepsilon}({\cal O})$ is a
$S_{\tau}$-invariant neighborhood of ${\cal O}$, and that
$\Psi_\tau$ is $S_{\tau}$-equivariant. We also require
$\varepsilon>0$ so small that $Q_{\varepsilon}({\cal O})$ contains
no other critical orbit besides ${\cal O}$, and that
$\Psi_\tau\bigr(\{y\}\times N({\cal O})_y(\varepsilon)\bigl)$ and
${\cal O}$ have a unique intersection point $y$ (after identifying
${\cal O}$ with the zero section  $N({\cal O})(\varepsilon)$). Then
${\cal L}_\tau\circ\Psi_\tau|_{N({\cal O})_y(\varepsilon)}$
possesses $y$ as an isolated critical point. Checking the proofs of
Theorem~7.3 and Corollary~7.1 in \cite[pp. 72]{Ch}, and replacing
$f\circ\exp|_{\xi_x}$ and $\exp_x\phi_x$ therein by ${\cal
L}_\tau\circ\Psi_\tau|_{N({\cal O})_x(\varepsilon)}$ and
$\Psi_\tau|_{N({\cal O})_x(\varepsilon)}\circ\phi_x$ for $x\in{\cal
O}$, one easily gets:

\begin{lemma}\label{lem:4.8}
For sufficiently small $0<\epsilon<\varepsilon$, there exist a
$S_\tau$-equivariant homeomorphism $\Phi_\tau$ from  $N({\cal
O})(\epsilon)$ to a $S_{\tau}$-invariant open neighborhood
$\Omega_\epsilon({\cal O})\subset Q_\varepsilon({\cal O})$ of ${\cal
O}$, and a $C^1$-map $h_\tau:N({\cal O})^0(\epsilon)\to N({\cal
O})^+(\epsilon)\oplus N({\cal O})^-(\epsilon)$ such that
\begin{eqnarray*}
&&{\cal L}_\tau\circ\Phi_\tau(y,
v)=\frac{1}{2}\left(\|P_+(y)v\|_1^2-\|P_-(y)v\|_1^2\right)\\
&&\hspace{23mm} + \, {\cal L}_\tau\circ\Psi_\tau\bigl((y, P_0(y)v) +
h_\tau(P_0(y)v)\bigr)
\end{eqnarray*}
for $(y, v)\in N({\cal O})(\epsilon)$, where $P_\ast$ is as in
(\ref{e:4.51}).
\end{lemma}

Let $N({\cal O})^\bot(\epsilon)= N({\cal O})^+(\epsilon)\oplus
N({\cal  O})^-(\epsilon)$ and write $v=v^0 + v^\bot$.  Set
\begin{equation}\label{e:4.57}
\left.\begin{array}{ll}
 &\Xi_\tau(y,
v^\bot)=\frac{1}{2}\left(\|P_+(y)v\|_1^2-\|P_-(y)v\|_1^2\right),\\
&\Upsilon_\tau(y, v^0)={\cal L}_\tau\circ\Psi_\tau\bigl((y, P_0(y)v)
+ h_\tau(P_0(y)v)\bigr)
\end{array}\right\}
\end{equation}
for $(y, v)\in N({\cal O})(\epsilon)$. Then define $F_\tau: N({\cal
O})(\epsilon)\to\R$ by
\begin{equation}\label{e:4.58}
F_\tau(y,v)={\cal L}_\tau\circ\Phi_\tau(y, v)=\Upsilon_\tau(y, v^0)+
\Xi_\tau(y, v^\bot)
\end{equation}
for all $(y, v)\in N({\cal O})(\epsilon)$.  ({\bf Note}: Though we
require the higher smoothness of $M$ and $L$ we do not know whether
or not  ${\cal L}_\tau$ has higher smoothness than order two unlike
the special $L$ considered in \cite{Lo2}. Hence from \cite[Th.7.3,
pp. 72]{Ch} we can only get that $\Phi_\tau$ is a homeomorphism.
However,  $N({\cal O})(\epsilon)$ is a $C^2$-bundle \footnote{The
requirements of the higher smoothness of $M$ and $L$ is used to
assure this.} and therefore
 \begin{equation}\label{e:4.59}
\hbox{both}\,\;\Xi_\tau\;\hbox{and}\; \Upsilon_\tau\;\hbox{are}\;
C^2.
\end{equation}
By the local trivialization of $N({\cal O})(\epsilon)$ the final
claim can be derived from Lemma~\ref{lem:4.2} and the proofs of
\cite[Th.5.1, pp. 44]{Ch} and \cite[Th.7.3, pp. 72]{Ch}.) Clearly,
both $\Upsilon_\tau$ and $\Xi_\tau$ are also $S_{\tau}$-invariant,
 and have the unique critical orbit ${\cal O}$ in
$N({\cal O})^\bot(\varepsilon)$ and $N({\cal O})^0(\varepsilon)$
respectively. Since $F_\tau$ is  $C^2$-smooth,
 we can follow \cite{Wa} to construct a Gromoll-Meyer
pair of ${\cal O}$ as a critical submanifold of $F_\tau$ on $N({\cal
O})(\varepsilon)$,
\begin{equation}\label{e:4.60}
(W({\cal O}), W({\cal O})^-).
\end{equation}
  (Note that different from \cite{Wa}
the present $S_{\tau}$-action on $N({\cal O})(\epsilon)$ is only
continuous; but the arguments there can still be carried out due to
the special property of our $S_{\tau}$-action in (\ref{e:4.20}) and
the definition of $F_{\tau}$.) In the present case, for any $y\in
{\cal O}$, $F_{\tau}|_{N({\cal O})_y(\epsilon)}$ has a unique
critical point $y$ in $N({\cal O})_y(\epsilon)$ ( the fibre of disk
bundle $N({\cal O})(\epsilon)$ at $y$), and
\begin{equation}\label{e:4.61}
\left(W({\cal O})_y, W({\cal O})_y^-\right):=\bigl(W({\cal O})\cap
N({\cal O})_y(\epsilon),\, W({\cal O})^-\cap
  N({\cal O})_y(\epsilon)\bigr)
   \end{equation}
is a  Gromoll-Meyer pair of  $F_{\tau}|_{N_y({\cal O})(\epsilon)}$
at its isolated critical point $y$ satisfying
\begin{equation}\label{e:4.62}
\bigl(W({\cal O})_{[s]_{\tau}\cdot y},\,  W({\cal
O})^-_{[s]_{\tau}\cdot y}\bigr)=
  \bigl([s]_{\tau}\cdot W({\cal O})_y,\,  [s]_{\tau}\cdot W({\cal O})^-_y\bigr)
\end{equation}
for any $[s]_{\tau}\in S_{\tau}$ and $y\in {\cal O}$. Clearly,
\begin{equation}\label{e:4.63}
\bigl(\widehat W({\cal O})), \widehat W({\cal
O})^-)\bigr):=\bigl(\Phi_\tau(W({\cal O})), \Phi_\tau(W({\cal
O})^-)\bigr)
\end{equation}
 is a topological Gromoll-Meyer pair of ${\cal L}_{\tau}$ at ${\cal O}$, which is also
$S_\tau$-invariant. Define
\begin{eqnarray}
&&C_\ast({\cal L}_\tau, {\cal O}; \K):=H_\ast\bigl(\widehat W({\cal
O})), \widehat W({\cal O})^-; \K\bigr),\label{e:4.64}\\
&&C_\ast(F_\tau, {\cal O}; \K):=H_\ast(W({\cal O}), W({\cal O})^-;
\K)\label{e:4.65}
\end{eqnarray}
via the relative singular homology.
 $\Phi_\tau$ induces an obvious isomorphism
\begin{equation}\label{e:4.66}
(\Phi_\tau)_\ast: C_\ast({\cal L}_\tau, {\cal O}; \K)
 \cong C_\ast(F_\tau, {\cal O};\K).
\end{equation}
Since the normal bundle $N({\cal O})$ is differentiably trivial, it
follows from \cite[(2.13), (2.14)]{Wa} (cf. also the shifting
theorem in [GM1] and [Ch]) that for any $q\in \{0\}\cup\N$,
\begin{eqnarray*}
  C_q(F_{\tau}, {\cal O}; \K)
&\cong& \oplus^q_{j=0}\left[C_{q-j}
  \left(F_{\tau}\Bigm|_{N({\cal O})_y(\epsilon)}, y; \K\right) \otimes
     H_j(S_{\tau}; \K)\right]  \\
    &\cong& \oplus^q_{j=0}\left[C_{q-j-m^-_{\tau}({\cal O})}
  \left(F_{\tau}\Bigm|_{N({\cal O})_y^0(\epsilon)}, y; \K\right) \otimes
     H_j(S_{\tau}; \K)\right]  \\
    &\cong& C_{q-1-m^-_{\tau}({\cal O})}
      \left(F_{\tau}\Bigm|_{N({\cal O})_y^0(\epsilon)},y; \K\right)\quad\forall y\in{\cal O}.
 \end{eqnarray*}
Here $C_{q-1-m^-_{\tau}({\cal O})} \bigl(F_{\tau}|_{N({\cal
O})_y^0(\epsilon)}, y; \K\bigr)$ is independent of the choice of
$y\in {\cal O}=S_{\tau}\cdot \gamma$. Taking $y=\gamma$ we obtain
\begin{equation}\label{e:4.67}
   C_*({\cal L}_{\tau},\, S_{\tau}\cdot \gamma;\K) \cong C_{*-1-m^-_{\tau}(S_{\tau}\cdot \gamma)}
\bigl(F_{\tau}|_{N(S_{\tau}\cdot \gamma)_\gamma^0(\epsilon)},
\gamma;\K\bigr).
\end{equation}

Suppose that $\psi^k({\cal O})=S_{k\tau}\cdot\gamma^k$ is also an
isolated critical orbit of the functional ${\cal L}_{k\tau}$ on
$H_{k\tau}(\alpha^k)$ for some $k\in\N$. Our purpose is to study the
relations between critical groups $C_\ast({\cal L}_\tau, {\cal O};
\K)$ and $C_\ast({\cal L}_{k\tau}, \psi^k({\cal O}); \K)$.

Let $N(S_{k\tau}\cdot\gamma^k)$ be the normal bundle of
$S_{k\tau}\cdot\gamma^k$ in $H_{k\tau}(\alpha^k)$ and
$$
N(S_{k\tau}\cdot\gamma^k)(\varepsilon)=\bigl\{(y,v)\in
N(S_{k\tau}\cdot\gamma^k)\,|\,y\in
  S_{k\tau}\cdot\gamma^k,\,\|v\|_1<\varepsilon\bigr\}.
   $$
Corresponding to (\ref{e:4.51}) there exist  natural orthogonal
bundle decompositions
\begin{eqnarray}
 &&N(\psi^k({\cal O}))=N(\psi^k({\cal O}))^+\oplus N(\psi^k({\cal O}))^-\oplus N(\psi^k({\cal
 O}))^0,\label{e:4.68}\\
 &&N(\psi^k({\cal O}))(\varepsilon)=N(\psi^k({\cal O}))^+(\varepsilon)
 \oplus N(\psi^k({\cal O}))^-(\varepsilon)\oplus N(\psi^k({\cal
 O}))^0(\varepsilon),\label{e:4.69}
\end{eqnarray}
where $N(\psi^k({\cal O}))^+(\varepsilon)= N(\psi^k({\cal
O}))(\varepsilon)\cap N(\psi^k({\cal O}))^\ast$ for $*=+,-,0$.

It is not hard to check that
\begin{equation}\label{e:4.70}
\psi^k(N({\cal O})(\varepsilon))\subset
N(S_{k\tau}\cdot\gamma^k)(\sqrt{k}\varepsilon)\quad{\rm and}\quad
\psi^k(N({\cal O})^\ast(\varepsilon))\subset
N(S_{k\tau}\cdot\gamma^k)^\ast
\end{equation}
for $\ast=+, 0, -$. By shrinking $\varepsilon>0$
 we have also a $C^2$-smooth $S_{\tau}$-equivariant diffeomorphism from
$N(S_{k\tau}\cdot\gamma^k)(\sqrt{k}\varepsilon)$ to a
$S_{k\tau}$-invariant open neighborhood
$Q_{\sqrt{k}\varepsilon}(S_{k\tau}\cdot\gamma^k)$ of
$S_{k\tau}\cdot\gamma^k$,
\begin{eqnarray}\label{e:4.71}
&&\Psi_{k\tau}: N(S_{k\tau}\cdot\gamma^k)(\sqrt{k}\varepsilon)\to
Q_{\sqrt{k}\varepsilon}(S_{k\tau}\cdot\gamma^k),\\
&&\qquad \Psi_{k\tau}(y, v)(t)=\exp_{y(t)}{v(t)}\;\forall
t\in\R.\nonumber
\end{eqnarray}

With the same arguments as above Lemma~\ref{lem:4.8}, by furthermore
shrinking $0<\epsilon<\varepsilon$, there exist a
$S_{k\tau}$-equivariant homeomorphism $\Phi_{k\tau}$ from from
$N(\psi^k({\cal O}))(\sqrt{k}\epsilon)$ to a $S_{k\tau}$-invariant
open neighborhood $\Omega_{\sqrt{k}\epsilon}(\psi^k({\cal
O}))\subset Q_{\sqrt{k}\varepsilon}(\psi^k({\cal O}))$ of
$\psi^k({\cal O})$, and a $C^1$-map
$$
h_{k\tau}:N(\psi^k({\cal O}))^0(\sqrt{k}\epsilon)\to N(\psi^k({\cal
O}))^+(\sqrt{k}\epsilon)\oplus N(\psi^k({\cal
O}))^-(\sqrt{k}\epsilon)
$$
such that
\begin{equation}\label{e:4.72}
{\cal L}_{k\tau}\circ\Phi_{k\tau}(y, v)=\Upsilon_{k\tau}(y, v^0)+
\Xi_{k\tau}(y, v^\bot)
\end{equation}
for $(y, v)\in N(\psi^k({\cal O}))(\sqrt{k}\epsilon)$, where
$v^\bot\in N(\psi^k({\cal O}))^\bot(\sqrt{k}\epsilon)=
N(\psi^k({\cal O}))^+(\sqrt{k}\epsilon)\oplus N(\psi^k({\cal
O}))^-(\sqrt{k}\epsilon)$ and
\begin{equation}\label{e:4.73}
\left.\begin{array}{ll}
 &\Xi_{k\tau}(y,
v^\bot)=\frac{1}{2}\left(\|v^+\|_1^2-\|v^-\|_1^2\right),\\
&\Upsilon_{k\tau}(y, v^0)={\cal L}_\tau\circ\Psi_\tau\bigl((y, v^0)
+ h(y, v^0)\bigr)
\end{array}\right\}
\end{equation}
have the similar properties to (\ref{e:4.59}). As in (\ref{e:4.58})
we define a $S_{k\tau}$-invariant, $C^2$-smooth function $F_{k\tau}:
N(\psi^k({\cal O}))(\sqrt{k}\epsilon)\to \R$ by
\begin{equation}\label{e:4.74}
F_{k\tau}(y,v) = {\cal
L}_{k\tau}(\Phi_{k\tau}(y,v))=\Upsilon_{k\tau}(y, v^0)+
\Xi_{k\tau}(y, v^\bot).
\end{equation}
It has the unique critical orbit $\psi^k({\cal O})$ in
$N(\psi^k({\cal O}))(\sqrt{k}\epsilon)$. Note that (\ref{e:4.55})
and (\ref{e:4.71}) imply
\begin{equation}\label{e:4.75}
\Psi_{k\tau}\circ\psi^k=\psi^k\circ\Psi_\tau.
\end{equation}
As in \cite[Prop. 2.3]{LoLu}, we can suitably modify the proof of
\cite[Lem. 3.3]{Lo2} to get:

\begin{lemma}\label{lem:4.9}
 Suppose that
$m^0_{k\tau}(\psi^k({\cal O}))= m^0_{\tau}({\cal O})$.   Then:
\begin{description}
  \item[(i)] The maps $h_{\tau}$ and $h_{k\tau}$  satisfy
\begin{equation}\label{e:4.76}
 h_{k\tau}(\psi^k(p))=\psi^k(h_{\tau}(p)),\quad\forall p=(y, v)\in
      N({\cal O})^0(\epsilon).
\end{equation}
\item[(ii)] The homeomorphisms $\Phi_{\tau}$ and
$\Phi_{k\tau}$ satisfy
\begin{equation}\label{e:4.77}
 \Phi_{k\tau}\circ \psi^k = \psi^k\circ \Phi_{\tau}
\end{equation}
as maps from $N({\cal O})(\epsilon)$ to $H_{k\tau}(\alpha^k)$.

 \item[(iii)] For $q\in N({\cal O})^0(\epsilon)$,
$p\in N({\cal O})^\bot(\epsilon)$,  there hold
\begin{equation}\label{e:4.78}
\Upsilon_{k\tau}(\psi^k(q))=k\Upsilon_{\tau}(q), \qquad
    \Xi_{k\tau}(\psi^k(p))=k\Xi_{\tau}(p).
    \end{equation}
    \end{description}
\end{lemma}
Indeed, the key in the proof of \cite[Lem. 3.3]{Lo2} is that the
maps $h_\tau$ and $h_{k\tau}$ are uniquely determined by the
implicit function theorem as showed in the proof of the Generalized
Morse lemma \cite[pp. 44]{Ch}. It follows from (\ref{e:4.78}) that
\begin{equation}\label{e:4.79}
F_{k\tau}\circ\psi^k=kF_\tau.
\end{equation}
By the construction of the Gromoll-Meyer pair in \cite{Wa} we can
construct such a pair of $F_{k\tau}$ at $\psi^k({\cal O})$ on
$N(\psi^k({\cal O}))(\sqrt{k}\epsilon)$, $(W(\psi^k({\cal O})),
W(\psi^k({\cal O}))$
 such that
\begin{equation}\label{e:4.80}
(\psi^k(W({\cal O})), \psi^k(W({\cal O})^-))\subset (W(\psi^k({\cal
O})), W(\psi^k({\cal O}))^-)
\end{equation}
for the pair $(W({\cal O}), W({\cal O})^-)$ in (\ref{e:4.60}). Set
\begin{equation}\label{e:4.81}
\bigl(\widehat W(\psi^k({\cal O})), \widehat W(\psi^k({\cal
O}))^-\bigr):=\bigl(\Phi_{k\tau}(W(\psi^k({\cal O}))),
\Phi_{k\tau}(W(\psi^k({\cal O}))^-)\bigr),
\end{equation}
which is a topological Gromoll-Meyer pair, and
\begin{eqnarray}
&&C_\ast({\cal L}_{k\tau}, \psi^k({\cal O});
\K):=H_\ast\bigl(\widehat W(\psi^k({\cal O})), \widehat
W(\psi^k({\cal O}))^-;\,\K\bigr),\label{e:4.82}\\
&&C_\ast(F_{k\tau}, \psi^k({\cal O}); \K):=H_\ast\bigl(
W(\psi^k({\cal O})), W(\psi^k({\cal O}))^-;\,\K\bigr).\label{e:4.83}
\end{eqnarray}
It follows from (\ref{e:4.77}) and (\ref{e:4.80}) that
\begin{equation}\label{e:4.84}
\bigl(\psi^k(\widehat W({\cal O})), \psi^k(\widehat W({\cal
O})^-)\bigr)\subset \bigl(\widehat W(\psi^k({\cal O})), \widehat
W(\psi^k({\cal O}))^-\bigr)
\end{equation}
and that $\psi^k$ induces  homomorphisms
\begin{eqnarray}
&&(\psi^k)_\ast: C_\ast({\cal L}_\tau, {\cal O}; \K)\to C_\ast({\cal
L}_{k\tau}, \psi^k({\cal O}); \K),\label{e:4.85}\\
&&(\psi^k)_\ast: C_\ast(F_\tau, {\cal O}; \K)\to C_\ast(F_{k\tau},
\psi^k({\cal O}); \K)\label{e:4.86}
\end{eqnarray}
satisfying
\begin{equation}\label{e:4.87}
(\psi^k)_\ast\circ
(\Phi_\tau)_\ast=(\Phi_{k\tau})_\ast\circ(\psi^k)_\ast
\end{equation}
because of (\ref{e:4.77}). By (\ref{e:4.66}) and the isomorphism
\begin{equation}\label{e:4.88}
(\Phi_{k\tau})_\ast: C_\ast({\cal L}_{k\tau}, \psi^k({\cal O}); \K)
 \cong C_\ast(F_{k\tau}, \psi^k({\cal O});\K)
\end{equation}
we only need to prove:

\begin{lemma}\label{lem:4.10}
The Gromoll-Meyer pairs  in (\ref{e:4.80}) can be chosen so that the
homomorphism in (\ref{e:4.86}) is an isomorphism provided that
\begin{equation}\label{e:4.89}
m^-_{k\tau}(\psi^k({\cal O}))= m^-_{\tau}({\cal O})\quad{\rm
and}\quad m^0_{k\tau}(\psi^k({\cal O}))= m^0_{\tau}({\cal O}).
\end{equation}
\end{lemma}

\noindent{\bf Proof.}\quad By (\ref{e:4.58}), (\ref{e:4.72}) and
(\ref{e:4.74}) we have
\begin{equation}\label{e:4.90}
\left.\begin{array}{ll}
&C_\ast(F_{\tau}, {\cal O};
\K)=C_\ast(\Upsilon_{\tau}+
\Xi_{\tau}, {\cal O}; \K),\\
&C_\ast(F_{k\tau}, \psi^k({\cal O}); \K)= C_\ast(\Upsilon_{k\tau}+
\Xi_{k\tau}, \psi^k({\cal O}); \K)
\end{array}\right\}
\end{equation}
We shall imitate the proof of Lemma~\ref{lem:4.6} to prove that the
homomorphism
\begin{equation}\label{e:4.91}
(\psi^k)_\ast: C_\ast(\Upsilon_{\tau}+ \Xi_{\tau}, {\cal O}; \K)\to
C_\ast(\Upsilon_{k\tau}+ \Xi_{k\tau}, \psi^k({\cal O}); \K)
\end{equation}
 is an isomorphism.

Let $(W_0({\cal O}), W_0^-({\cal O}))$ be a Gromoll-Meyer pair of
$\Upsilon_\tau$ at ${\cal O}\subset N({\cal O})^0(\epsilon)$. Since
(\ref{e:4.89}) implies that $\psi^k: N({\cal O})^0(\epsilon)\to
N(\psi^k({\cal O}))^0(\sqrt{k}\epsilon)$ is a bundle isomorphism.
Hence
$$\bigl(\psi^k(W_0({\cal O})), \psi^k(W_0^-({\cal O}))\bigr)$$
is a Gromoll-Meyer pair of $\Upsilon_{k\tau}$ at $\psi^k({\cal
O})\subset N(\psi^k({\cal O}))^0(\sqrt{k}\epsilon)$.
 For $j=1, k$ let us write
$N(\psi^j({\cal O}))^\bot= N(\psi^j({\cal O}))^+\oplus
N(\psi^j({\cal O}))^-$ and
\begin{eqnarray*}
&&N(\psi^j({\cal O}))^\bot(\sqrt{j}\epsilon)= N(\psi^j({\cal
O}))^+(\sqrt{j}\epsilon)\oplus N(\psi^j({\cal
O}))^-(\sqrt{j}\epsilon),\\
&&\Xi_{j\tau}(y, v^\bot)=\Xi_{j\tau}^+(y, v^+) + \Xi_{j\tau}(y,
v^-),\quad v^\bot=v^+ + v^-.
\end{eqnarray*}
By (\ref{e:4.78}), for  $p\in N({\cal O})^\pm(\epsilon)$, there hold
\begin{equation}\label{e:4.92}
    \Xi^\pm_{k\tau}(\psi^k(p))=k\Xi^\pm_{\tau}(p).
    \end{equation}

Let $(W_{11}({\cal O}), W_{11}^-({\cal O}))$ be a Gromoll-Meyer pair
of $\Xi^-_\tau$ at ${\cal O}\subset N({\cal O})^-(\epsilon)$. Then
\begin{equation}\label{e:4.93}
(\psi^k(W_{11}({\cal O})), \psi^k(W_{11}^-({\cal O}))
\end{equation}
is a Gromoll-Meyer pair of $\Xi^-_{k\tau}$  at $\psi^k({\cal
O})\subset N(\psi^k({\cal O}))^-(\sqrt{k}\epsilon)$ because
(\ref{e:4.89}) implies that $\psi^k:N({\cal O})^-(\epsilon)\to
N(\psi^k({\cal O}))^-(\sqrt{k}\epsilon)$ is a bundle isomorphism.
For $0<\delta\ll\epsilon$, set
\begin{eqnarray*}
&&W_{12}:=\{(y, v)\in N({\cal O})^+(\epsilon)\,|\,
\|v\|_\tau\le\delta\,\},\\
&& W^-_{12}:=\{(y, v)\in
N({\cal O})^+(\epsilon)\,|\, \|v\|_\tau=\delta\,\},\\
&&V_{12}:=\{(y, v)\in N(\psi^k({\cal O}))^+(\epsilon)\,|\,
\|v\|_{k\tau}\le
\sqrt{k}\delta\,\},\\
&& V^-_{12}:=\{(y, v)\in N(\psi^k({\cal O}))^+(\epsilon)\,|\,
\|v\|_{k\tau}=\sqrt{k}\delta\,\}.
\end{eqnarray*}
 Then $(W_{12}, W^-_{12})$ (resp. $(V_{12}, V^-_{12})$) is a Gromoll-Meyer pair of
$\Xi^+_{\tau}$ (resp. $\Xi^+_{k\tau}$) at ${\cal O}\subset N({\cal
O})^+(\epsilon)$ (resp.  $\psi^k({\cal O})\subset N(\psi^k({\cal
O}))^+(\sqrt{k}\epsilon)$), and that
\begin{equation}\label{e:4.94}
\left(\psi^k(W_{12}), \psi^k(W^-_{12})\right)\subset (V_{12},
V^-_{12}).
\end{equation}
By Lemma 5.1 on the page 51 of \cite{Ch}, we may take
\begin{eqnarray*}
&&W_1({\cal O}):=W_{11}({\cal O})\oplus W_{12},\\
&&W_1^-({\cal O}):=(W_{11}({\cal O})\oplus W_{12}^-)\cup
(W_{11}^-({\cal
O})\oplus W_{12}),\\
&&V:=\psi^k(W_{11}({\cal O}))\oplus V_{12}, \\
&&V^-:=(\psi^k(W_{11}({\cal O}))\oplus V_{12}^-)\cup
(\psi^k(W_{11}^-({\cal O}))\oplus V_{12})
\end{eqnarray*}
and get a Gromoll-Meyer pair of $\Upsilon_\tau+ \Xi_{\tau}$  at
${\cal O}\subset N({\cal O})(\epsilon)$, $(W({\cal O}), W({\cal
O})^-)$, where
\begin{eqnarray}
&&W({\cal O}):=W_0({\cal O})\oplus W_{11}({\cal O})\oplus W_{12},\label{e:4.95}\\
&&W^-({\cal O}):=\bigl(W_0({\cal O})\oplus [(W_{11}({\cal O})\oplus
W_{12}^-)\cup
(W_{11}^-({\cal O})\oplus W_{12})]\bigr)\nonumber\\
&&\hspace{20mm} \cup \bigl(W_0^-({\cal O})\oplus W_{11}({\cal
O})\oplus W_{12}\bigr).\label{e:4.96}
\end{eqnarray}
Therefore
\begin{equation}\label{e:4.97}
C_\ast(\Upsilon_{\tau}+ \Xi_{\tau}, 0; \K)=H_\ast(W({\cal O}),
W^-({\cal O});\K).
\end{equation}

Similarly,  we have a Gromoll-Meyer pair of $\Upsilon_{k\tau} +
\Xi_{k\tau}$ at $\psi^k({\cal O})\subset N(\psi^k({\cal
O}))(\sqrt{k}\epsilon)$, $(W(\psi^k({\cal O})), W(\psi^k({\cal
O}))^-)$, where
\begin{eqnarray*}
&&W(\psi^k({\cal O})):=\psi^k(W_0({\cal O}))\oplus \psi^k(W_{11}({\cal O}))\oplus V_{12},\\
&&W^-(\psi^k({\cal O})):=\bigl(\psi^k(W_0({\cal O}))\oplus
[(\psi^k(W_{11}({\cal O}))\oplus V_{12}^-)\cup
(\psi^k(W_{11}^-({\cal O}))\oplus V_{12})]\bigr)\nonumber\\
&&\hspace{25mm} \cup \bigl(\psi^k(W_0^-({\cal O}))\oplus
\psi^k(W_{11}({\cal O}))\oplus V_{12}\bigr).
\end{eqnarray*}
It follows that
\begin{equation}\label{e:4.98}
C_\ast(\Upsilon_{k\tau}+ \Xi_{k\tau}, \psi^k({\cal O});
\K)=H_\ast(W(\psi^k({\cal O})), W^-(\psi^k({\cal O}));\K).
\end{equation}
Note that $\psi^k(W({\cal O}))=\psi^k(W_0({\cal O}))\oplus
\psi^k(W_{11}({\cal O}))\oplus \psi^k(W_{12})$ and
\begin{eqnarray}\label{e:4.99}
&&\psi^k(W^-({\cal O}))=\bigl(\psi^k(W_0({\cal O}))\oplus
(\psi^k(W_{11}({\cal O}))\oplus \psi^k(W_{12}^-))\nonumber\\
&&\hspace{30mm}\cup (\psi^k(W_{11}^-({\cal O}))\oplus
\psi^k(W_{12}))\bigr)\nonumber\\
 &&\hspace{30mm}\cup \bigl(\psi^k(W_0^-({\cal O}))\oplus
\psi^k(W_{11}({\cal O}))\oplus\psi^k(W_{12})\bigr).
\end{eqnarray}
Since $\psi^k: N^+({\cal O})\to N^+(\psi^k({\cal O}))$ is a
continuous bundle injection, by (\ref{e:4.94}) and the constructions
of $(V_{12}, V_{12}^-)$ and $(W_{12}, W_{12}^-)$ above
(\ref{e:4.94}) it is readily checked that $\bigl(\psi^k(W_{12}),
\psi^k(W^-_{12})\bigr)$ is a deformation retract of $(V_{12},
V^-_{12})$. It follows that
$$
\bigl(\psi^k(W({\cal O})), \psi^k(W^-({\cal O}))\bigr)\subset
\bigl(W(\psi^k({\cal O})), W^-(\psi^k({\cal O}))\bigr)
$$
is a deformation retract of $\bigl(W(\psi^k({\cal O})),
W^-(\psi^k({\cal O}))\bigr)$. Hence
$$
(\psi^k)_\ast: H_\ast(W({\cal O}), W^-({\cal O});\K\bigr)\to
H_\ast\bigl(W(\psi^k({\cal O})), W^-(\psi^k({\cal O}));\K\bigr)
$$
is an isomorphism. Therefore, by (\ref{e:4.97}) and (\ref{e:4.98}),
the homomorphism $(\psi^k)_\ast$ in (\ref{e:4.91}) is an
isomorphism. Lemma~\ref{lem:4.10} is proved. $\Box$.

When $\gamma$ is constant, i.e. ${\cal O}=S_\tau\cdot\gamma$ is an
isolated critical point, this case has been proved in
Theorem~\ref{th:4.4}. Combing this with Lemma~\ref{lem:4.10}, and
(\ref{e:4.66}) and (\ref{e:4.88}) we get

\begin{theorem}\label{th:4.11}
  For an isolated critical submanifold ${\cal
O}={\rm S}_{\tau}\cdot \gamma$ of ${\cal L}_{\tau}$ in
$H_{\tau}(\alpha)$, suppose that for some $k\in\N$ the critical
submanifold $\psi^k({\cal O})={\rm S}_{k\tau}\cdot \gamma^k$ of
${\cal L}_{k\tau}$ in $H_{k\tau}(\alpha^k)$ is also isolated, and
that (\ref{e:4.89}) is satisfied, i.e. $m^-_{k\tau}(S_{k\tau}\cdot
\gamma^k)=m^-_{\tau}(S_{\tau}\cdot \gamma)$ and
 $m^0_{k\tau}(S_{k\tau}\cdot \gamma^k)=m^0_{\tau}(S_{\tau}\cdot \gamma)$.
  Then for $c={\cal L}_\tau|_{\cal O}$ and small
   $\epsilon>0$ there exist
  topological Gromoll-Meyer pairs  of ${\cal L}_{\tau}$ at ${\cal O}\subset H_\tau(\alpha)$ and
 of ${\cal L}_{k\tau}$ at $\psi^k({\cal O})\subset H_{k\tau}(\alpha^k)$
\begin{eqnarray*}
&&(\widehat W({\cal O}), \widehat W({\cal
  O})^-)\subset\bigl(({\cal L}_\tau)^{-1}[c-\epsilon,
c+ \epsilon], ({\cal L}_\tau)^{-1}(c-\epsilon)\bigr)
 \quad\hbox{and}\\
&&\bigl(\widehat W(\psi^k({\cal O})), \widehat W(\psi^k({\cal
  O}))^-\bigr)\subset\bigl(({\cal
L}_{k\tau})^{-1}[kc-k\epsilon, kc+ k\epsilon], ({\cal
L}_{k\tau})^{-1}(kc-k\epsilon)\bigr),
\end{eqnarray*}
such that
$$
\bigl(\psi^k(\widehat W({\cal O})), \psi^k(\widehat W({\cal
  O})^-)\bigr)\subset \bigl(\widehat W(\psi^k({\cal O})), \widehat W(\psi^k({\cal
  O}))^-\bigr)
$$
and that the iteration map $\psi^k: H_{\tau}(\alpha)\to
H_{k\tau}(\alpha^k)$ induces an isomorphism:
\begin{eqnarray*}
 &&\psi^k_*:  C_*({\cal L}_{\tau}, {\cal O}; \K):=H_\ast\bigl(\widehat W({\cal O}), \widehat W({\cal
  O})^-; \K\bigr)\\
 &&\qquad \longrightarrow
         C_*({\cal L}_{k\tau}, \psi^k({\cal O}); \K):
         =H_\ast\bigl(\widehat W(\psi^k({\cal O})), \widehat W(\psi^k({\cal
  O}))^-; \K\bigr).
\end{eqnarray*}
\end{theorem}

\begin{lemma}\label{lem:4.12}
Suppose that $C_q({\cal L}_{\tau}, {\cal O};\K)\ne 0$ for ${\cal
O}=S_{\tau}\cdot \gamma$. Then
\begin{equation}\label{e:4.100}
q-2n\le q-1- m^0_{\tau}({\cal O})\le m^-_{\tau}({\cal O})\le q-1
\end{equation}
if ${\cal O}$ is not a single point critical orbit, i.e. $\gamma$ is
not constant,  and
\begin{equation}\label{e:4.101}
q-2n\le q- m^0_{\tau}({\cal O})\le m^-_{\tau}({\cal O})\le q
\end{equation}
otherwise.
\end{lemma}

\noindent{\bf Proof}.\quad  If $\gamma$ is not a constant solution,
 it follows from (\ref{e:4.66}) and (\ref{e:4.67}) that
\begin{equation}\label{e:4.102}
   C_{q-1-m^-_{\tau}({\cal O})}
\bigl(F_{\tau}|_{N({\cal O})^0_{\gamma}(\epsilon)},
\gamma;\K\bigr)\cong C_q({\cal L}_{\tau},\, {\cal O};\K)\ne 0.
\end{equation}
Since $\gamma$ is an isolated critical point of $F_{\tau}|_{N({\cal
O})^0_{\gamma}(\epsilon)}$ in $N({\cal O})^0_{\gamma}(\epsilon)$ and
$N({\cal O})^0_{\gamma}(\epsilon)$ has dimension $m^0_{\tau}({\cal
O})$, we get
\begin{equation}\label{e:4.103}
0\le q-1-m^-_{\tau}({\cal O})\le \dim N({\cal
O})^0_{\gamma}(\epsilon)=m^0_{\tau}({\cal O}).
\end{equation}
 By (\ref{e:4.53}), $m^0_{\tau}({\cal
O})=m^0_{\tau}(\gamma)-1\le 2n-1$. (\ref{e:4.100})  easily follows
from this and (\ref{e:4.103}).

If $\gamma$ is a constant solution, i.e. ${\cal O}=\{\gamma\}$,
using the isomorphisms above (\ref{e:4.11}) and (\ref{e:4.24}) we
derive
$$
C_q(\tilde\alpha_{\tau}+ \tilde\beta_{\tau}, 0; \K)\cong
C_q(\tilde{\cal L}_{\tau}, \tilde\gamma; \K)\ne 0,\quad{\rm
where}\quad \tilde\gamma=(\phi_{\tau})^{-1}(\gamma).
$$
On the other hand,  (\ref{e:3.11}) and the shifting theorem
(\cite{GM1} and \cite[pp.50]{Ch}) imply
$$
C_q(\alpha_{\tau}+ \beta_{\tau}, 0; \K) \cong
C_{q-m^-_{\tau}(\gamma)}(\tilde\alpha_{\tau}, 0; \K).
$$
Since $\tilde\alpha_{\tau}$ is defined on a manifold of dimension
$m^0_{\tau}(\gamma)\le 2n$, (\ref{e:4.101}) follow immediately.
$\Box$\vspace{2mm}

\begin{lemma}\label{lem:4.13}
Suppose that $C_q({\cal L}_{\tau}, {\cal O};\K)\ne 0$ for ${\cal
O}=S_{\tau}\cdot \gamma$. If either ${\cal O}$ is not a single point
critical orbit and $q>1$, or ${\cal O}$ is a single point critical
orbit and $q>0$, then each point in ${\cal O}$ is non-minimal saddle
point.
\end{lemma}

\noindent{\bf Proof}.\quad  When ${\cal O}$ is a single point
critical orbit and $q>0$, the conclusion follows from \cite[Ex.1,
pp.33]{Ch}. Now assume that ${\cal O}$ is not a single point
critical orbit and $q>1$. For any $y\in{\cal O}$, by (\ref{e:4.66})
and the formula above (\ref{e:4.67}) we have
\begin{eqnarray*}
  0\ne C_q(F_{\tau}, {\cal O}; \K)
&\cong& \oplus^q_{j=0}\left[C_{q-j}
  \left(F_{\tau}\Bigm|_{N({\cal O})_y(\epsilon)}, y; \K\right) \otimes
     H_j(S_{\tau}; \K)\right]  \\
    &\cong& C_{q-1}
      \left(F_{\tau}\Bigm|_{N({\cal O})_y(\epsilon)},y; \K\right).
 \end{eqnarray*}
Since $y$ is an isolated critical point of $F_{\tau}\Bigm|_{N({\cal
O})_y(\epsilon)}$ and $q-1>0$, we derive from \cite[Ex.1, pp.33]{Ch}
that $y$ is a non-minimal saddle point of $F_{\tau}\Bigm|_{N({\cal
O})_y(\epsilon)}$. This implies that $y$ is a non-minimal saddle
point of ${\cal L}_\tau$ on the submanifold $\Psi_\tau(N({\cal
O})_y(\epsilon))\subset H_\tau(\alpha)$ (and therefore on
$H_\tau(\alpha)$). $\Box$\vspace{2mm}

\section{Proof of Theorem~\ref{th:1.1}}
\setcounter{equation}{0}

\noindent{\bf 5.1. Proof of (i)}.\quad For any $\tau\in\N$, let
$H_\tau(\alpha^k)$ denote the Hilbert manifold of $W^{1,2}$-loops
$\gamma:\R/\tau\Z\to M$ representing $\alpha^k$. Since
$H_r(C(\R/\tau\Z, M; \alpha^k);\K)= H_r(C(\R/\Z, M; \alpha^k);\K)$
and the inclusion $H_\tau(\alpha^k)\hookrightarrow C(\R/\tau\Z, M;
\alpha^k)$ is a homotopy equivalence,
\begin{equation}\label{e:5.1}
{\rm rank}H_r(H_\tau(\alpha^k); \K)\ne 0\quad\forall \tau, k\in\N.
\end{equation}
By \cite{Be} the functional ${\cal L}_\tau$ on the Hilbert manifold
$H_\tau(\alpha^k)$ is $C^2$-smooth, bounded below, satisfies the
Palais-Smale condition, and all critical points of it have finite
Morse indexes and nullities.  In particular,  \textsf{the critical
set ${\cal K}({\cal L}_{\tau}, \alpha^k)$ of ${\cal L}_\tau$ on
$H_\tau(\alpha^k)$ is nonempty} because ${\cal L}_{\tau}$ can attain
the minimal value on $H_{\tau}(\alpha^k)$. Clearly, for any $\tau,
k\in\N$ we may assume that each  critical point of
 ${\cal L}_\tau$ on $H_\tau(\alpha^k)$ is isolated.
  By contradiction we make:

 \noindent{\bf Assumption $F(\alpha)$}: (i) {\it For any given
 integer $k>0$, the system (\ref{e:1.6}) only possesses  finitely many
 distinct,  $k$-periodic  solutions representing $\alpha^k$},
 (ii) {\it there exists an integer $k_0>1$ such that for each integer $k>k_0$, any
 $k$-periodic
 solution $\tilde\gamma$ of the system (\ref{e:1.6}) representing $\alpha^k$
must be an iteration of some $l$-periodic solution $\gamma$ of the
system (\ref{e:1.6}) representing $\alpha^l$ with $l\le k_0$ and
$k=ls$ for some $s\in\N$.}

Under this assumption  we have  integer periodic solutions
$\hat\gamma_i$ of the system (\ref{e:1.6}) of period $\tau_i\le k_0$
and representing $\alpha^{\tau_i}$, $i=1,\cdots,p$, such that for
each integer $k>k_0$ any  integer $k$-periodic solution $\gamma$ of
the system (\ref{e:1.6}) representing $\alpha^k$  must be an
iteration of some $\hat\gamma_i$, i.e. $\gamma=\hat\gamma_i^l$ for
some $l\in\N$ with $l\tau_i=k$. Set $\tau:=k_0!$ (the factorial of
$k_0$) and $\gamma_i=\hat\gamma_i^{\tau/\tau_i}$, $i=1,\cdots, p$.
Then each $\gamma_i$ is a  $\tau$-periodic solution of the system
(\ref{e:1.6}) representing $\alpha^\tau$. We conclude

\begin{claim}\label{cl:5.1}
 For any $k\in\N$, it holds that
\begin{equation}\label{e:5.2}
{\cal K}({\cal L}_{k\tau},
\alpha^{k\tau})=\bigl\{\,\gamma^k_j\,\bigm|\,1\le j\le p\,
    \bigr\}.
\end{equation}
 \end{claim}

\noindent{\bf Proof}.\quad Let $\gamma\in {\cal K}({\cal L}_{k\tau},
\alpha^{k\tau})$. Since $k\tau>k_0$, by (ii) in Assumption
$F(\alpha)$ we have $\gamma=\hat\gamma_i^l$ for some $l\in\N$ with
$l\tau_i=k\tau$. Hence
$\gamma=\hat\gamma_i^l=(\hat\gamma_i)^{k\tau/\tau_i}
=(\hat\gamma_i^{\tau/\tau_i})^k=\gamma_i^k$. $\Box$\vspace{2mm}

Since $M$ is not assumed to be orientable, it is possible that the
pullback bundle $\gamma_j^\ast TM\to \R/\tau\Z$ is not trivial.
However,  each $2$-fold iteration $\gamma_j^2$, $(\gamma_j^2)^\ast
TM\to \R/2\tau\Z$ is always trivial. Note that (\ref{e:5.2}) implies
\begin{equation}\label{e:5.3}
{\cal K}({\cal L}_{2k\tau},
\alpha^{2k\tau})=\bigl\{\,(\gamma^2_j)^k=\gamma_j^{2k}\,\bigm|\,1\le
j\le p\, \bigr\}.
\end{equation}
Hence replacing $\{\gamma_1\,\cdots\, \gamma_p\}$  by
$\{\gamma^2_1\,\cdots\, \gamma^2_p\}$ we may \textbf{assume}:
\begin{equation}\label{e:5.4}
\gamma_j^\ast TM\to \R/\tau\Z,\quad j=1,\cdots, p,\quad  \textsf{are
all trivial}.
\end{equation}

\begin{lemma}\label{lem:5.2}
For each $k\in\N$ there exists  $\gamma'_k\in{\cal K}({\cal
L}_{k\tau}, \alpha^{k\tau})$
 such that
$$
C_r({\cal L}_{k\tau}, \gamma'_k; \K)\ne 0\quad{\rm and}\quad r-2n\le
r-m^0_{k\tau}(\gamma'_k)\le m^-_{k\tau}(\gamma'_k)\le r.
$$
\end{lemma}

\noindent{\bf Proof}.\quad Let $c_1<\cdots <c_l$ be all critical
values of ${\cal L}_\tau$, $l\le p$. Then $kc_1<\cdots <kc_l$ are
all critical values of ${\cal L}_{k\tau}$, $k=1, 2, \cdots$. In
particular, $\inf{\cal L}_{k\tau}=kc_1$ because ${\cal L}_{k\tau}$
is bounded below and satisfies the Palais-Smale condition.

 By (\ref{e:5.1}), ${\rm rank}
H_r(H_{k\tau}(\alpha^{k\tau}); \K)\ge m$ for some $m\in\N$. Recall
that a subset of an abelian group is defined to be {\bf linearly
independent} if it satisfies the usual condition with integer
coefficients, cf. \cite[pp. 87]{Ma}. Take linearly independent
elements of $H_r(H_{k\tau}(\alpha^{k\tau}); \K)$, $\beta_1, \cdots,
\beta_m$, and singular cycles $Z_1,\cdots, Z_m$ of
$H_{k\tau}(\alpha^{k\tau})$ which represent them. Let $S$ be a
compact set containing the supports of $Z_1,\cdots, Z_m$. Then
$S\subset ({\cal L}_{k\tau})_b:=\{{\cal L}_{k\tau}\le b\}$ for a
sufficiently large regular value $b>kc_l$. Note that $Z_1,\cdots,
Z_m$ are also singular cycles of $({\cal L}_{k\tau})_b$, and that
non-trivial $\K$-linear combination of them cannot be homologous to
zero in $({\cal L}_{k\tau})_b$ (otherwise the same combination is
homologous to zero in $H_{k\tau}(\alpha^{k\tau})$.) Hence we get
$$
{\rm rank}H_r(({\cal L}_{k\tau})_b;\K)\ge m>0.
$$
 Take the regular values of
${\cal L}_{k\tau}$, $a_0<a_1<\cdots <a_l=b$ such that $kc_i\in
(a_{i-1}, a_i)$, $i=1,\cdots, l$. By Theorem 4.2 of \cite[pp.
23]{Ch},
\begin{equation}\label{e:5.5}
H_r(({\cal L}_{k\tau})_{a_i}, ({\cal L}_{k\tau})_{a_{i-1}}; \K)\cong
\bigoplus_{{\cal L}_{k\tau}(z)=kc_i,\, d{\cal L}_{k\tau}(z)=0
}C_r({\cal L}_{k\tau}, z; \K).
\end{equation}
Since each critical point has finite Morse index, it follows from
the generalized Morse lemma that each group $C_r({\cal L}_{k\tau},
z; \K)$ has finite rank, and therefore that
$$
{\rm rank}H_r(({\cal L}_{k\tau})_{a_i}, ({\cal L}_{k\tau})_{
a_{i-1}}; \K)<+\infty,\;i=1,\cdots, l.
$$
By the arguments on the page 38 of \cite{Ch} and the fact (b) on the
page 87 of \cite{Ma}, for a triple $Z\subset Y\subset X$ of
topological spaces it holds that
$$
{\rm rank}H_q(X, Z;\K)\le {\rm rank}H_q(X, Y;\K)+ {\rm rank}H_q(X,
Y;\K)
$$
if these three numbers are finite. It follows that
\begin{eqnarray*}
0<m&\le &{\rm rank}H_r(({\cal L}_{k\tau})_b;\K)\\
&=& {\rm rank}H_r(({\cal
L}_{k\tau})_{a_l}, ({\cal L}_{k\tau})_{a_{0}}; \K)\\
&\le& \sum^m_{i=1}{\rm rank}H_r(({\cal L}_{k\tau})_{a_i}, ({\cal
L}_{k\tau})_{a_{i-1}}; \K)<+\infty.
\end{eqnarray*}
Hence ${\rm rank}H_r(({\cal L}_{k\tau})_{a_i}, ({\cal L}_{k\tau})_{
a_{i-1}}; \K)\ge 1$ for some $i$. By (\ref{e:5.5}) we get a
$\gamma'_k\in{\cal K}({\cal L}_{k\tau}, \alpha^{k\tau})$ such that
${\rm rank}C_r({\cal L}_{k\tau}, \gamma'_k; \K)\ne 0$ and thus
$C_r({\cal L}_{k\tau}, \gamma'_k; \K)\ne 0$. Noting (\ref{e:5.4}),
we can use
 the isomorphism
above (\ref{e:4.11}) to derive
$$
C_r(\tilde{\cal L}_{k\tau}, \tilde\gamma'_k; \K)\ne 0,\quad{\rm
where}\quad \tilde\gamma'_k=(\phi_{k\tau})^{-1}(\gamma_k').
$$
Replaceing $\tilde\gamma^k$ in (\ref{e:4.12}) by $\tilde\gamma_k'$,
and using the isomorphism above (\ref{e:4.24}), (\ref{e:3.11}) and
the shifting theorem (\cite{GM1} and \cite[pp.50]{Ch}) we get
$$
C_{r-m^-_{k\tau}(\gamma'_k)}(\tilde\alpha_{k\tau}, 0; \K)\cong
C_r(\alpha_{k\tau}+ \beta_{k\tau}, 0; \K)\cong C_r(\tilde{\cal
L}_{k\tau}, \tilde\gamma'_k; \K)\ne 0.
$$
Since $\tilde\alpha_{k\tau}$ is defined on a manifold of dimension
$m^0_{k\tau}(\gamma'_k)\le 2n$, the desired inequalities follow
immediately.  $\Box$\vspace{2mm}

\begin{lemma}\label{lem:5.3}
\textsf{Without Assumption $F(\alpha)$}, let $\gamma$ be an isolated
critical point  of ${\cal L}_\tau$ in $H_\tau(\alpha^\tau)$ such
that $\gamma^\ast TM\to S_{\tau}$ is trivial. For every integer
$q\ge n+1$, let $k(q, \gamma)=1$ if $\hat m^-_\tau(\gamma)=0$, and
$k(q, \gamma)= \frac{q+n}{\hat m^-_\tau(\gamma)}$ if $\hat
m^-_\tau(\gamma)\ne 0$. Assume that $\gamma^k$ is also an isolated
critical point of ${\cal L}_{k\tau}$ for some integer $k> k(q,
\gamma)$. Then
\begin{equation}\label{e:5.6}
C_q({\cal L}_{k\tau},  \gamma^k; \K)=0.
\end{equation}
\end{lemma}

\noindent{\bf Proof.}\quad Let $\phi_{k\tau}: W^{1,2}(S_\tau,
B^n_\rho(0))\to H_{k\tau}(\alpha^{k\tau})$ be a coordinate chart on
$H_{k\tau}(\alpha^{k\tau})$ around $\gamma^k$ as in (\ref{e:3.8}).
Set $\tilde\gamma=(\phi_{\tau})^{-1}(\gamma)$. Then
$\tilde\gamma^k=(\phi_{k\tau})^{-1}(\gamma^k)$ and
$m^-_{\tau}(\tilde\gamma)=m^-_{\tau}(\gamma)$,
$m^0_{\tau}(\tilde\gamma)=m^0_{\tau}(\gamma)$ and
$m^-_{k\tau}(\tilde\gamma^k)=m^-_{k\tau}(\gamma^k)$ and
$m^0_{k\tau}(\tilde\gamma^k)=m^0_{k\tau}(\gamma^k)$. As in the proof
of Lemma~\ref{lem:5.2}, by the isomorphisms above (\ref{e:4.11}) and
(\ref{e:4.24}) we have
\begin{eqnarray*}
C_q({\cal L}_{k\tau}, \gamma^k; \K)&\cong &C_r(\tilde{\cal
L}_{k\tau},
\tilde\gamma^k; \K)\\
&\cong&C_q(\tilde\alpha_{k\tau}+ \tilde\beta_{k\tau}, 0; \K)\\
&\cong&C_{q-m^-_{k\tau}(\gamma^k)}(\tilde\alpha_{k\tau}, 0; \K).
\end{eqnarray*}
Here $\tilde\alpha_{k\tau}$ is defined on a manifold of dimension
$m^0_{k\tau}(\gamma^k)\le 2n$.

 If $\hat m^-_\tau(\gamma)=0$, by
(\ref{e:3.2}) (or (\ref{e:3.18})) we have $0\le
m^-_{k\tau}(\gamma^k)\le n-m^0_{k\tau}(\gamma^k)$. Hence
$$
q- m^-_{k\tau}(\gamma^k)\ge q-(n-m^0_{k\tau}(\gamma^k))\ge 1+
m^0_{k\tau}(\gamma^k).
$$
This gives $C_{q-m^-_{k\tau}(\gamma^k)}(\tilde\alpha_{k\tau}, 0;
\K)=0$.

If $\hat m^-_\tau(\gamma)>0$,  by (\ref{e:3.2}) (or (\ref{e:3.18}))
we have $k\hat m^-_\tau(\gamma)-n\le m^-_{k\tau}(\gamma^k)$ and thus
$$
q- m^-_{k\tau}(\gamma^k)\le q- (k\hat m^-_\tau(\gamma)-n)= q+
n-k\hat m^-_\tau(\gamma)<0
$$
 if $k>\frac{q+ n}{\hat m^-_{\tau}(\gamma)}$. This also leads to
$C_{q-m^-_{k\tau}(\gamma^k)}(\alpha_{k\tau}, 0; \K)=0$.
$\Box$\vspace{2mm}

So we immediately get the following generalization of Lemma 4.2 in
\cite{Lo2}.

\begin{corollary}\label{cor:5.4}
Under Assumption $F(\alpha)$, for every integer $q\ge n+1$  there
exists a constant $k_0(q)>0$ such that for every integer $k\ge
k_0(q)$  there holds
$$
C_q({\cal L}_{k\tau},  y; \K)=0\quad \forall y\in{\cal K}({\cal
L}_{k\tau}, \alpha^{k\tau}).
$$
Here $k_0(q)=1$ if $\hat m^-_{\tau}(\gamma_j)=0$ for all $1\le j\le
p$, and
$$
k_0(q)=1+ {\rm max}\bigg\{\left[\frac{q+ n}{\hat
m^-_{\tau}(\gamma_j)}\right]\,\Bigm|\,\hat m^-_{\tau}(\gamma_j)\ne
0,\, 1\le j\le p\biggr\}
$$ otherwise. ($[s]$ denotes the largest integer less than or equal to $s$).
\end{corollary}

Indeed, by (\ref{e:5.2}) we may assume $y=\gamma^k_j$ for some $1\le
j\le p$. Then Lemma~\ref{lem:5.3} yields the desired conclusion.
\vspace{2mm}

Clearly, if $r\ge n+1$ then Lemma~\ref{lem:5.2} and
Lemma~\ref{lem:5.3} immediately give a contradiction.
Theorem~\ref{th:1.1}(i) is proved in this case.

 \verb"In the following we"
\verb"consider the case" $r=n$.

Under Assumption $F(\alpha)$ we apply Lemma \ref{lem:5.2} to all
$k\in \{2^m\,|\,m\in\{0\}\cup{\N}\}$ to get an infinite subsequence
$Q$ of $\{2^m\,|\,m\in\{0\}\cup{\N}\}$, some $l\in\N$ and an
$\gamma\in\{\gamma_1,\cdots, \gamma_p\}$ such that $C_n({\cal
L}_{kl\tau}, \gamma^{kl}; \Z_2)\ne 0$,
 $m^-_{kl\tau}(\gamma^{kl})=
m^-_{l\tau}(\gamma^{l})$ and
$m^0_{kl\tau}(\gamma^{kl})=m^0_{l\tau}(\gamma^{l})$ for any $k\in
Q$. In order to save notations we {\bf always assume $l=1$} in the
following. That is, we have $\gamma^k\in{\cal K}({\cal L}_{k\tau},
\alpha^{k\tau})$ with
\begin{equation}\label{e:5.7}
\left.\begin{array}{ll}
 &C_n({\cal L}_{k\tau}, \gamma^{k}; \K)\ne 0,\\
 &m^-_{k\tau}(\gamma^{k})=
m^-_{\tau}(\gamma),\quad m^0_{k\tau}(\gamma^{k})=m^0_{\tau}(\gamma)
\end{array}\right\}
\end{equation}
for any $k\in Q$. By Corollary~\ref{cor:5.4} there exists $k_0>0$
such that for any $\gamma\in\{\gamma_1,\cdots, \gamma_p\}$,
\begin{equation}\label{e:5.8}
C_{n+1}({\cal L}_{k\tau},  \gamma^k; \K)=0\quad\forall k\in
Q(k_0):=\{k\in Q\,|\, k\ge k_0\}.
\end{equation}

{\bf To avoid} the finite energy homology introduced and used in
\cite{Lo2} we need to improve the proof and conclusions of Theorem
4.3 in \cite{Lo2}. Let $c={\cal L}_\tau(\gamma)$. Take $\epsilon>0$
sufficiently small so that for each $k\in\N$ the interval
$[k(c-3\epsilon), k(c+ 3\epsilon)]$ contains an unique critical
value $kc$ of ${\cal L}_{k\tau}$ on $H_{k\tau}(\alpha^{k\tau})$,
i.e.
$$
{\cal L}_{k\tau}\bigr({\cal K}({\cal L}_{k\tau},
\alpha^{k\tau})\bigl)\cap [k(c-3\epsilon), k(c+ 3\epsilon)]=\{kc\}.
$$
By Theorem~\ref{th:4.4}, for each integer $k\in Q$ we may choose
topological Gromoll-Meyer pairs  of ${\cal L}_\tau$ at $\gamma$ and
${\cal L}_{k\tau}$ at $\gamma^k$, $\bigl( W(\gamma),
W(\gamma)^-\bigr)$ and $\bigl( W(\gamma^k), W(\gamma^k)^-\bigr)$,
such that
\begin{eqnarray}
\bigl( W(\gamma), W(\gamma)^-\bigr)\subset \left(({\cal
L}_\tau)^{-1}([c-2\epsilon, c+ 2\epsilon]), ({\cal
L}_\tau)^{-1}(c-2\epsilon)\right),\label{e:5.9}\\
\bigl( W(\gamma^k), W(\gamma^k)^-\bigr)\subset  \bigl(({\cal
L}_{k\tau})^{-1}([kc-2k\epsilon, kc+ 2k\epsilon]), ({\cal
L}_{k\tau})^{-1}(kc-2k\epsilon)\bigr),\label{e:5.10}\\
\bigl(\psi^k(W(\gamma)),  \psi^k(W(\gamma)^-)\bigr)\subset
\bigl(W(\gamma^k),  W(\gamma^k)^-\bigr)\label{e:5.11}
\end{eqnarray}
and that  the iteration map $\psi^k: H_\tau(\alpha)\to
H_{k\tau}(\alpha^{k})$ induces isomorphisms
\begin{eqnarray*}
&&(\psi^k)_\ast: C_\ast({\cal L}_\tau, \gamma;\K)=
H_\ast\bigl(W(\gamma),  W(\gamma)^-; \K\bigr)\\
&&\hspace{40mm} \to C_\ast({\cal L}_{k\tau},
\gamma^k;\K)=H_\ast\bigl(W(\gamma^k), W(\gamma^k)^-;\K\bigr).
\end{eqnarray*}
 For $j=1, k$, denote by the inclusions
\begin{eqnarray*}
&&{\bf h}_1^j:\bigl( W(\gamma^j), W(\gamma^j)^-\bigr)\hookrightarrow
\bigl(({\cal L}_{j\tau})_{j(c+ 2\epsilon)},
({\cal L}_{j\tau})_{j(c-2\epsilon)}\bigr),\\
&&{\bf h}_2^j:\bigl(({\cal L}_{j\tau})_{j(c+ 2\epsilon)}, ({\cal
L}_{j\tau})_{j(c-2\epsilon)}\bigr) \hookrightarrow \bigl(({\cal
L}_{j\tau})_{j(c+ 2\epsilon)}, ({\cal
L}_{j\tau})^\circ_{j(c-\epsilon)}\bigr),\\
&&{\bf h}_3^j:  \bigl(({\cal L}_{j\tau})_{j(c+ 2\epsilon)}, ({\cal
L}_{j\tau})^\circ_{j(c-\epsilon)}\bigr)\hookrightarrow
\bigl(H_{j\tau}, ({\cal L}_{j\tau})^\circ_{j(c-\epsilon)}\bigr).
\end{eqnarray*}
Hereafter $B^\circ$ denote the interior of $B$ without special
statements.  The arguments above \cite[Th.4.3]{Lo2} show that
\begin{eqnarray*}
&&({\bf h}_2^j\circ{\bf h}_1^j)_\ast:H_\ast\bigl( W(\gamma^j),
W(\gamma^j)^-;\K\bigr)\to H_\ast\bigl(({\cal L}_{j\tau})_{j(c+
2\epsilon)}, ({\cal
L}_{j\tau})^\circ_{j(c-\epsilon)};\K\bigr),\\
&&({\bf h}_3^j)_\ast: H_\ast\bigl(({\cal L}_{j\tau})_{j(c+
2\epsilon)}, ({\cal L}_{j\tau})^\circ_{j(c-\epsilon)};\K\bigr)\to
H_\ast(\bigl(H_{j\tau}, ({\cal
L}_{j\tau})^\circ_{j(c-\epsilon)};\K\bigr)
\end{eqnarray*}
are monomorphisms on homology modules. For $j=1, k$, we have also
inclusions
\begin{eqnarray*}
&&{\bf I}_j:\bigl( W(\gamma^j), W(\gamma^j)^-\bigr)\hookrightarrow
\bigl(({\cal L}_{j\tau})^{-1}([jc-2j\epsilon, jc+ 2j\epsilon]),
({\cal L}_{j\tau})^{-1}(jc-2j\epsilon)\bigr),\\
&&{\bf J}_j:\bigl(({\cal L}_{j\tau})^{-1}([jc-2j\epsilon, jc+
2j\epsilon]), ({\cal
L}_{j\tau})^{-1}(jc-2j\epsilon)\bigr)\hookrightarrow
\bigl(H_{j\tau}, ({\cal L}_{j\tau})^\circ_{jc-j\epsilon}\bigr).
\end{eqnarray*}
It is clear that
\begin{equation}\label{e:5.12}
{\bf J}_j\circ{\bf I}_j={\bf h}^j_3\circ{\bf h}^j_2\circ{\bf
h}^j_1,\quad j=1,k.
\end{equation}
By (\ref{e:5.11}),  we have also
$$
\psi^k\circ{\bf I}_1={\bf I}_k\circ\psi^k
$$
as maps from $\bigl( W(\gamma), W(\gamma)^-\bigr)$ to $\bigl(({\cal
L}_{k\tau})^{-1}([kc-2k\epsilon, kc+ 2k\epsilon]), ({\cal
L}_{k\tau})^{-1}(kc-2k\epsilon)\bigr)$.
 So we get the following result, which is a slightly strengthened
 version of \cite[Th. 4.3]{Lo2} in the case $M=T^n$.

\begin{proposition}\label{prop:5.5}
Under Assumption $F(\alpha)$, there exist a periodic solution
$\gamma$ of (\ref{e:1.6}) of integer period $\tau$ and representing
$\alpha$, a large integer $k_0>0$, an infinite integer set $Q$
containing $1$, and a small $\epsilon>0$ having properties: For any
$k\in Q(k_0):=\{k\in Q\,|\, k\ge k_0\}$ there exist topological
Gromoll-Meyer pairs $\bigl( W(\gamma), W(\gamma)^-\bigr)$ and
$\bigl( W(\gamma^k), W(\gamma^k)^-\bigr)$ satisfying
(\ref{e:5.9})-(\ref{e:5.11}) such that for the inclusion
$$
j_{k\tau}={\bf h}^k_3\circ {\bf h}^k_2\circ {\bf
h}^k_1:\bigl(W(\gamma^k), W(\gamma^k)^-\bigr)\to
\bigl(H_{k\tau}(\alpha^{k\tau}), ({\cal
L}_{k\tau})^\circ_{k(c-\epsilon)}\bigr)
$$
the following diagram holds:
\begin{eqnarray}\label{e:5.13}
 &&0\ne C_n({\cal L}_{\tau},
\gamma;\K)\stackrel{(\psi^k)_*}{\longrightarrow} C_n({\cal
L}_{k\tau},
\gamma^k;\K)\nonumber\\
&&\hspace{30mm}\stackrel{(j_{k\tau})_*}{\longrightarrow}
H_n\bigl(H_{k\tau}(\alpha^{k\tau}), ({\cal
L}_{k\tau})^\circ_{k(c-\epsilon)};\K\bigr)\equiv{\cal H}_k,
 \end{eqnarray}
 where  $c={\cal L}_\tau(\gamma)$, $(\psi^k)_*$ is an isomorphism, and $(j_{k\tau})_*$ is a
monomorphism among the singular homology modules. In particular, if
$\omega$ is a generator of $C_n({\cal L}_\tau,
\gamma;\K)=H_n\left(W(\gamma), W(\gamma)^-;\K\right)$, then
\begin{eqnarray}
&&(j_{k\tau})_*\circ(\psi^k)_*(\omega)\ne 0\quad{\rm in}\;
      {\cal H}_k,\label{e:5.14}\\
&&(j_{k\tau})_\ast\circ(\psi^k)_\ast(\omega)=({\bf
J}_k)_\ast\circ({\bf I}_k)_\ast\circ(\psi^k)_\ast(\omega)\nonumber\\
&&\hspace{29mm}=({\bf J}_k)_\ast\circ(\psi^k)_\ast\circ ({\bf
I}_1)_\ast(\omega)\quad{\rm in}\;
      {\cal H}_k.\label{e:5.15}
\end{eqnarray}
\end{proposition}

\textsf{It is (\ref{e:5.15}) that helps us avoiding to use the
finite energy homology.}

The notion of a $C^1$-smooth singular simplex in Hilbert manifolds
can be defined as on page 252 of \cite{Ma}.

\begin{proposition}\label{prop:5.6}
For $\tau\in\N$, $c\in\R$, $\epsilon>0$, $q\geq 0$, and a
$C^1$-smooth $q$-simplex
$$\eta:(\Delta_q,\partial\Delta_q)\to
\left(H_{\tau}(\alpha^\tau),({\cal
L}_{\tau})_{c-\epsilon}^{\circ}\right),
$$
there exists
an integer $k(\eta)>0$ such that for every integer $k\geq k(\eta)$,
the $q$-simplex
$$
\eta^k\equiv
\psi^k(\eta):(\Delta_q,\partial\Delta_q)
            \to \left(H_{k\tau}(\alpha^{k\tau}),({\cal L}_{k\tau})_{k(c-\epsilon)}^{\circ}\right)
            $$
is  homotopic to a  singular $q$-simplex
\begin{equation}\label{e:5.16}
\eta_k:(\Delta_q,\partial\Delta_q)\to
 \left(({\cal L}_{k\tau})_{k(c-\epsilon)}^{\circ},({\cal L}_{k\tau})_{k(c-\epsilon
 )}^{\circ}\right)
 \end{equation}
  with $\eta^k = \eta_k$ on $\partial\Delta_q$
 and the homotopy fixes $\eta^k|_{\partial\Delta_q}$.
\end{proposition}

This is an analogue of \cite[Th.1]{BK}, firstly proved by Y. Long
\cite[Prop. 5.1]{Lo2} in the case $M=T^n$. Proposition 5.1 in
\cite{Lo2} actually gave stronger conclusions under weaker
assumptions: If the $q$-simplex $\eta$ above is only a finite energy
one ($C^1$-smooth simplex must be of finite energy), then the
simplex $\eta^k$ is finite energy homotopic to a finite energy
$q$-simplex $\eta_k$. Hence Proposition~\ref{prop:5.6} can be
derived with the same reason as in \cite[Prop. 5.1]{Lo2} as long as
we generalize an inequality as done in Lemma A.4 of Appendix. But we
also give necessary details for the reader's
convenience.\vspace{2mm}

\noindent{\bf Proof of Proposition~\ref{prop:5.6}}.\quad Recall that
for paths $\sigma:[a_1, a_2]\to M$ and $\delta:[b_1, b_2]\to M$ with
$\sigma(a_2)=\delta(b_1)$ one often define new paths
$\sigma^{-1}:[a_1, a_2]\to M$ by $ \sigma^{-1}(t):=\sigma(a_2+
a_1-t)$
 and $\sigma\ast\delta:[a_1,
a_2+ b_2-b_1]\to M$ by $\sigma\ast\delta|_{[a_1, a_2]}=\sigma$ and
$$
\sigma\ast\delta(t):=\delta(t-a_2+ b_1)\quad{\rm for}\quad t\in
[a_2, a_2+ b_2-b_1].
$$

{\bf Given a $C^1$-path $\rho:[a,b]\to H_\tau(\alpha^\tau)$ and an
integer $k\ge 3$ we want to construct a path $\rho_k:[a,b]\to
H_{k\tau}(\alpha^{k\tau})$ such that
$$
\rho_k(a)=\psi^k(\rho(a))\quad{\rm and}\quad
\rho_k(b)=\psi^k(\rho(b)).
$$}
Define the initial point curve $\beta_\rho$ of $\rho$ by
$$
[a, b]\to M,\; s\mapsto\beta_\rho(s)=\rho(s)(0).
$$
It is $C^1$-smooth.  Following \cite[pp. 460]{Lo2} and \cite[pp.
381]{BK}, for $0\le s\le (b-a)/k$ and $1\le j\le
k-2$ define\\
$\bullet$ $\tilde\rho_k(a+s)=\rho(a)^{k-1}\ast\left(\beta_\rho|_{[a,
a+ ks]}\right)\ast\rho(a+ ks)\ast\left(\beta_\rho|_{[a,
a+ ks]}\right)^{-1}$,\\
$\bullet$ $\tilde\rho_k(a+ j(b-a)/k +
s)=\rho(a)^{k-j-1}\ast\left(\beta_\rho|_{[a, a+
ks]}\right)\ast\rho(a+ ks)\ast\left(\beta_\rho|_{[a,
a+ ks]}\right)\ast\rho(b)^j\ast(\beta_\rho)^{-1}$,\\
$\bullet$ $\tilde\rho_k(b-(b-a)/k + s)=\rho(a+
ks)\ast\left(\beta_\rho|_{[a, a+
ks]}\right)\ast\rho(b)^{k-1}\ast\left(\beta_\rho|_{[a, a+
ks]}\right)^{-1}$.

\noindent{These are piecewise $C^1$-smooth loops in $M$ representing
$\alpha^k$, and
$$
\tilde\rho(a)=\rho(a)^{k-1}\quad{\rm and}\quad
\tilde\rho(b)=\rho(b)\ast\beta_\rho\ast\rho(b)^{k-1}\ast\beta_\rho^{-1}.
$$
For each $u\in [a, b]$}, reparametrising the loop $\tilde\rho_k(u)$
on $\R/k\tau$ as in \cite[pp.461]{Lo2} we get a piecewise
$C^1$-smooth loop $\rho_k(u)\in H_{k\tau}(\alpha^{k\tau})$ and
therefore a piecewise $C^1$-smooth path $\rho_k:[a,b]\to
H_{k\tau}(\alpha^{k\tau})$ with
$\rho_k(a)=\psi^k(\rho(a))=\rho(a)^k$ and
$\rho_k(b)=\psi^k(\rho(b))=\rho(b)^k$.

Replacing all the terms of powers of $\rho(a)$ and $\rho(b)$ by the
constant point paths in the definition of $\tilde\rho_k$ above, we
get a piecewise $C^1$-smooth path $\beta_{\rho,k}:[a, b]\to
H_\tau(\alpha)$. For $s\in [a, b]$ and $j=[k(s-a)/(b-a)]$, by the
arguments of \cite[pp. 461]{Lo2},
\begin{eqnarray}\label{e:5.17}
{\cal L}_{k\tau}(\rho_k(s))&=& (k-j-1){\cal L}_\tau(\rho(a))+
j{\cal L}_\tau(\rho(b))+ {\cal L}_\tau(\beta_{\rho,k}(s))\nonumber\\
&\le &(k-1)M_0(\rho)+ M_1(\rho)+ 2M_2(\rho),
\end{eqnarray}
where $M_0(\rho)=\max\{{\cal L}_\tau(\rho(a)), {\cal
L}_\tau(\rho(b))\}$, $M_1(\rho)=\max_{a\le s\le b}|{\cal
L}_\tau(\rho(s))|$ and
\begin{equation}\label{e:5.18}
M_2(\rho)=\int^b_a\bigl|L(s, \beta_\rho(s),
\dot\beta_\rho(s))\bigr|ds.
\end{equation}
Note that ($L_3$) implies
\begin{equation}\label{e:5.19}
|L(t, q, v)|\le C(1+ \|v\|^2)\quad\forall (t,q, v)\in \R\times TM
\end{equation}
for some constant $C>0$. Therefore it follows from Lemma A.4 that
\begin{eqnarray*}
M_2(\rho)&=&\int^b_a\bigl|L(s, \beta_\rho(s),
\dot\beta_\rho(s))\bigr|ds\\
&\le& (b-a)C + C\int^b_a\bigl|\dot\beta_\rho(s)\bigr|^2ds\le (b-a)C+
\frac{1+\tau}{2\tau}Cc(\rho).
\end{eqnarray*}
This and (\ref{e:5.17}) yield
\begin{equation}\label{e:5.20}
\lim_{k\to +\infty}\sup\max_{a\le s\le b}\frac{1}{k}{\cal
L}_{k\tau}(\rho_k(s))\le M_0(\rho).
\end{equation}

Next replacing \cite[Lem.2.3]{Lo2} by Lemma A.4, and almost
repeating the reminder arguments of the proof of
\cite[Prop.5.1]{Lo2}, we can complete the proof of
Proposition~\ref{prop:5.6}. $\Box$\vspace{3mm}

\begin{lemma}\label{lem:5.7}\!\!\!{\rm (\cite[Lem.1]{BK})}\quad
Let $(X, A)$ be a pair of topological spaces and $\beta$ a singular
relative $p$-cycle of $(X,A)$. Let $\Sigma$ denote the set of
singular simplices of $\beta$ together with all their faces. Suppose
to every $\sigma\in\Sigma$, $\sigma:\Delta^q\to X$, $0\le q\le p$,
there is assigned a map $P(\sigma):\Delta^q\times [0, 1]\to X$ such
that
\begin{description}
\item[(i)] $P(\sigma)(z,0)=\sigma(z)$ for $z\in\Delta^q$,
\item[(ii)] $P(\sigma)(z,t)=\sigma(z)$ if $\sigma(\Delta^q)\subset
A$,
\item[(iii)] $P(\sigma)(\Delta^q\times\{1\})\subset A$,
\item[(iv)] $P(\sigma)\circ (e^i_q\times id)=P(\sigma\circ e^i_q)$
for $0\le i\le q$.
\end{description}
Then the homology class $[\beta]\in H_p(X,A)$ vanishes.
\end{lemma}

For the class $\omega$ in (\ref{e:5.15}), by the definition of ${\bf
I}_1$ above (\ref{e:5.12}) we have
\begin{equation}\label{e:5.21}
({\bf I}_1)_\ast(\omega)\in
H_n\bigl(({\cal L}_{\tau})^{-1}([c-2\epsilon, c+ 2\epsilon]),
({\cal L}_{\tau})^{-1}(c-2\epsilon);\K\bigr).
\end{equation}
Since both $({\cal L}_\tau)^{-1}([c-2\epsilon, c+ 2\epsilon])$ and
$({\cal L}_\tau)^{-1}(c-2\epsilon)$ are at least $C^2$-smooth
Hilbert manifolds, we can choose a $C^1$-smooth cycle representative
$\sigma$ of the class $({\bf I}_1)_\ast(\omega)$. Denote by
$\Sigma(\sigma)$  the set of all simplexes together with all their
faces contained in $\sigma$. By \cite[Ex.1, pp.33]{Ch} each
$\gamma^k$ in (\ref{e:5.7}) is a non-minimal saddle point of ${\cal
L}_{k\tau}$ on $H_{k\tau}(\alpha^{k\tau})$. As in the proof of
\cite[Prop. 5.2]{Lo2} we can use Proposition~\ref{prop:5.6} and
Lemma A.4 to get the corresponding result \textsf{without using the
finite energy homology}.

\begin{proposition}\label{prop:5.8}
  There exists a sufficiently large integer
$k(\sigma)\geq k_0$ such that for every integer $k\in Q(k(\sigma))$
and for every $\mu\in\Sigma(\sigma)$ with $\mu:\Delta_r\to
H_{\tau}(\alpha^\tau)$ and $0\leq r\leq n$, there exists a homotopy
$P(\psi^k(\mu)):\Delta_r\times [0,1]\to H_{k\tau}(\alpha^{k\tau})$
such that the properties (i) to (iv) in Lemma~\ref{lem:5.7} hold for
$(X,A)=\bigl(H_{k\tau}(\alpha^{k\tau}),({\cal
L}_{k\tau})_{k(c-\epsilon)}^{\circ}\bigr)$.
\end{proposition}

It follows that the homology class $({\bf
J}_k)_\ast\circ(\psi^k)_\ast\circ ({\bf I}_1)_\ast(\omega)\in
      {\cal H}_k$ vanishes. By (\ref{e:5.15}),
$(j_{k\tau})_\ast\circ(\psi^k)_\ast(\omega)=0$ in
      ${\cal H}_k$.
This  contradicts to (\ref{e:5.14}).
 Therefore Assumption $F(\alpha)$
can not hold.   Theorem~\ref{th:1.1}(i) is proved. \vspace{2mm}

\noindent{\bf 5.2. Proof of (ii)}.\quad
 Since the
inclusion $E_\tau\hookrightarrow C(\R/\tau\Z, M)$ is a homotopy
equivalence, and therefore ${\rm rank}H_r(E_{\tau}; \K)\ne 0$ for
all $\tau\in\N$. Consider the functional ${\cal L}_{k\tau}$ on
$E_{k\tau}$. It has still a nonempty critical point set. Replace
Assumption $F(\alpha)$ by

 \noindent{\bf Assumption $F$}: (i)
{\it For any given
 integer $k>0$, the system (\ref{e:1.6}) only possesses  finitely many
 distinct,  $k$-periodic  solutions},
 (ii) {\it there exists an integer $k_0>1$ such that for each integer $k>k_0$, any
 $k$-periodic solution $\tilde\gamma$ of the system (\ref{e:1.6})
must be an iteration of some $l$-periodic solution $\gamma$ of the
system (\ref{e:1.6}) with $l\le k_0$ and $k=ls$ for some $s\in\N$.}

Then slightly modifying the proof of (i) above one can complete the
proof.
 $\Box$

\section{Proof of Theorem~\ref{th:1.4}}
\setcounter{equation}{0}

The proof is similar to that of Theorem~\ref{th:1.1}. We only give
the main points. Identifying $\R/\tau\Z=[-\frac{\tau}{2},
\frac{\tau}{2}]/\{-\frac{\tau}{2}, \frac{\tau}{2}\}$, let
$$
C(\R/\tau\Z, M)_e:=\{x\in C(\R/\tau\Z, M)\,|\, x(-t)=x(t)\;
-\tau/2\le t\le\tau/2\}.
$$
We have a contraction from $C(\R/\tau\Z, M)_e$ to the subset of
constant loops in $C(\R/\tau\Z, M)_e$ which is identified with $M$:
$$
[0, 1]\times C(\R/\tau\Z, M)_e\to C(\R/\tau\Z, M)_e,\;(s, x)\mapsto
x_s,
$$
where $x_s(t)=x(st)$ for $-\tau/2\le t\le\tau/2$. Since the
inclusion $C(\R/\tau\Z, M)_e\hookrightarrow EH_\tau$ is also a
homotopy equivalence, we get
\begin{equation}\label{e:6.1}
H_n(EH_\tau; \Z_2)=H_n(C(\R/\tau\Z, M)_e; \Z_2)=H_n(M;\Z_2)\ne 0
\end{equation}
for any $\tau>0$. Note that ${\cal L}^E_\tau$ can always attain the
minimal value on $EH_\tau$ and therefore has a nonempty critical set
${\cal K}({\cal L}^E_{\tau})$. Under the conditions (L1)-(L4) we
replace the Assumption F($\alpha$) in \S 5 by\vspace{2mm}

 \noindent{\bf Assumption FE}: (i) {\it For any given
 integer $k>0$, the system (\ref{e:1.6}) possesses only finitely many distinct
 \verb"reversible" $k\tau$-periodic  solutions}, (ii) {\it there exists an
integer $k_0>1$ such that for each integer $k>k_0$, any
\verb"reversible" $k\tau$-periodic
  solution $\tilde\gamma$ of the system (\ref{e:1.6}) is an iteration of
some reversible $l\tau$-periodic solution $\gamma$ of the system
(\ref{e:1.6})
 with $l\le k_0$ and $k=ls$ for some
$s\in\N$.}\vspace{2mm}

  Under this assumption, as the arguments below Assumption F($\alpha$)
  we may get an integer $\tau\in\N$ and
finitely many reversible $\tau$-periodic  solutions of the system
(\ref{e:1.6}), $\gamma_1\,\cdots\, \gamma_p$, such that for any
$k\in\N$ every reversible $k\tau$-periodic  solution of the system
(\ref{e:1.6}) has form $\gamma^k_j$ for some $1\le j\le p$. Namely,
\begin{equation}\label{e:6.2}
{\cal K}({\cal L}^E_{k\tau})=\bigl\{\,\gamma^k_j\,\bigm|\,1\le j\le
p\,\bigr\}.
\end{equation}
Using the same proof as one of Lemma~\ref{lem:5.2} we may obtain:

\begin{lemma}\label{lem:6.1}
Under Assumption FE,  for each $k\in\N$ there exists  a critical
point $\gamma'_k$ of ${\cal L}^E_{k\tau}$ such that
\begin{equation}\label{e:6.3}
C_n({\cal L}^E_{k\tau}, \gamma'_k; \Z_2)\ne 0\quad{\rm and}\quad
-n\le n-m^0_{1, k\tau}(\gamma'_k)\le m^-_{1, k\tau}(\gamma'_k)\le n.
\end{equation}
\end{lemma}

Let $k_0=1$ if $\hat m^-_{1,\tau}(\gamma_j)=0$ for all $1\le j\le
p$, and
$$
k_0=1+{\rm max}\left\{\left[\frac{3n+2}{2\hat
m^-_{1,\tau}(\gamma_j)}\right]\,\Bigm|\,\hat
m^-_{1,\tau}(\gamma_j)\ne 0,\, 1\le j\le p\right\}
$$ otherwise. Corresponding with Corollary~\ref{cor:5.4} we have the
following generalization of \cite[Lem.4.4]{LuW2}.

\begin{lemma}\label{lem:6.2}
Under Assumption FE, for any integer number $k\ge k_0$, every
isolated critical point $z$ of ${\cal L}^E_{k\tau}$  has the
trivial $(n+1)$-th critical module, i.e.
$$ C_{n+1}({\cal L}^E_{k\tau}, z; \K)=0.$$
\end{lemma}

\noindent{\bf Proof.} Using the chart $\phi^E_{k\tau}$ in
(\ref{e:3.29}) let $\tilde z=(\phi^E_{k\tau})^{-1}(z)$. We only need
to prove
\begin{equation}\label{e:6.4}
 C_{n+1}(\tilde{\cal L}^E_{k\tau}, \tilde z; \K)=0\quad\forall k\ge k_0.
 \end{equation}
 Let $z=\gamma^k_j$ and thus $\tilde
z=\tilde\gamma^k_j$ with
$\tilde\gamma_j=(\phi^E_\tau)^{-1}(\gamma_j)$. By (\ref{e:4.41}), it
follows from Shifting theorem (\cite[p.50, Th. 5.4]{Ch}) and the
K\"{u}nneth formula that
$$\begin{array}{rcl}
C_{n+1}(\tilde{\cal L}^E_{k\tau}, \tilde z; \K)&\cong&C_{n+1}
(\alpha^E_{k\tau}+ \beta^E_{k\tau},0; \K)\\
&\cong&C_{n+1-m^-_{1,k\tau}(\tilde\gamma^k_j)}
(\alpha^E_{k\tau},0; G)\otimes C_{m^-_{1,m\tau}(\tilde\gamma^k_j)}(\beta^E_{k\tau},0; \K)\\
&\cong&C_{n+1-m^-_{1,k\tau}(\tilde\gamma^k_j)} (\alpha^E_{k\tau},0; \K)\otimes \K\\
&\cong&C_{n+1-m^-_{1,k\tau}(\tilde\gamma^k_j)} (\alpha^E_{k\tau},0;
\K)
\end{array}
$$
 because  $0$ is a nondegenerate critical
point of quadratic function $\beta^E_{k\tau}$. If (\ref{e:6.4})
does not hold, we get that $0\le n+1-m^-_{1,
k\tau}(\tilde\gamma^k_j)\le m^0_{1, k\tau}(\tilde\gamma^k_j)$
because $\tilde\gamma_{k\tau}$ is defined on a manifold of
dimension $m^0_{1,k\tau}(\tilde\gamma^k_j)$. Note that
\begin{eqnarray*}
&&m^-_{1, k\tau}(\tilde\gamma^k_j)=m^-_{k\tau}(\tilde{\cal
L}^E_{k\tau}, \tilde\gamma^k_j)= m^-_{1, k\tau}(\gamma^k_j),\\
&&m^0_{1, k\tau}(\tilde\gamma^k_j)=m^0_{k\tau}(\tilde{\cal
L}^E_{k\tau}, \tilde\gamma^k_j)= m^0_{1, k\tau}(\gamma^k_j).
\end{eqnarray*}
We have
\begin{eqnarray}
&&m^-_{1,k\tau}(\tilde{\cal L}^E_{k\tau}, \tilde\gamma^k_j)\le
n+1\le m^-_{1,k\tau}(\tilde{\cal L}^E_{k\tau}, \tilde\gamma^k_j)+
m^0_{1,k\tau}(\tilde{\cal L}^E_{k\tau},
\tilde\gamma^k_j)\label{e:6.5}\\
\hbox{\verb"or"} &&m^-_{1,k\tau}(\gamma^k_j)\le n+1\le
m^-_{1,k\tau}(\gamma^k_j)+ m^0_{1,k\tau}(\gamma^k_j)\label{e:6.6}
\end{eqnarray}
for any $k\in\N$.   By (\ref{e:2.24})
$$
m^-_{k\tau}(\tilde{\cal L}^E_{k\tau}, \tilde\gamma^k_j)-2n\le
m^-_{2, k\tau}(\tilde{\cal L}^E_{k\tau}, \tilde\gamma^k_j)\le
m^-_{k\tau}(\tilde{\cal L}^E_{k\tau}, \tilde\gamma^k_j)\quad\forall
k\in\N.
$$
Hence it follows from this, (\ref{e:3.26}) and (\ref{e:6.5}) that
\begin{eqnarray*}
  2k\hat m^-_\tau(\tilde{\cal L}^E_{\tau},\tilde\gamma_j)-n &\le &
m^-_{2, k\tau}(\tilde{\cal L}_{k\tau}, \tilde\gamma_j^k)+
m^-_{k\tau}(\tilde{\cal L}^E_{k\tau},
 \tilde\gamma_j^k)\\
 &\le& 2m^-_{k\tau}(\tilde{\cal L}^E_{k\tau},
 \tilde\gamma_j^k)\le 2n+ 2.
 \end{eqnarray*}
Therefore, when $\hat m^-_{1, \tau}(\gamma_j)=\hat
m^-_\tau(\tilde{\cal L}^E_{\tau}, \tilde\gamma_j)>0$, $k\le
\big[\frac{3n+2}{2\hat m_{1,\tau}(\gamma_j)}\big]$, which
contradicts to $k\ge k_0$.

When $\hat m^-_{1, \tau}(\gamma_j)=\hat m^-_\tau(\tilde{\cal
L}^E_{\tau}, \tilde\gamma_j)=0$, (\ref{e:3.22}) and (\ref{e:6.6})
also give a contradiction. The desired (\ref{e:6.4})  is proved.
$\Box$

Now as the arguments below Corollary~\ref{cor:5.4},  under
Assumption FE we may use Lemma~\ref{lem:6.1}  to get an infinite
subsequence $Q$ of $\{2^m\,|\,m\in\{0\}\cup{\N}\}$ and an
$\gamma\in\{\gamma_1,\cdots, \gamma_p\}$ such that
\begin{equation}\label{e:6.7}
\left.\begin{array}{ll}
 &C_n({\cal L}^E_{k\tau}, \gamma^{k}; \Z_2)\ne 0,\\
 &m^-_{1,k\tau}(\gamma^{k})=
m^-_{1,\tau}(\gamma),\quad
m^0_{1,k\tau}(\gamma^{k})=m^0_{1,\tau}(\gamma)
\end{array}\right\}
\end{equation}
for any $k\in Q$. By Lemma~\ref{lem:6.2}, for any
$x\in\{\gamma_1,\cdots, \gamma_p\}$ we have also
\begin{equation}\label{e:6.8}
C_{n+1}({\cal L}^E_{k\tau},  x^k; \K)=0\quad\forall k\in
Q(k_0):=\{k\in Q\,|\, k\ge k_0\}.
\end{equation}
 Then from
Proposition~\ref{prop:5.5} to the end of \S5.1 we only need to
make suitable replacements for some notations such as
$H_{j\tau}(\alpha^{j\tau})$, ${\cal L}_{j\tau}$  by $EH_{j\tau}$,
${\cal L}^E_{j\tau}$ for $j=1,k$, and so on,  and can complete the
proof of Theorem~\ref{th:1.4}.

\section{Proof of Theorem~\ref{th:1.6}}
\setcounter{equation}{0}

\noindent{\bf 7.1. Proof of (i)}.\quad Denote by ${\cal KO}({\cal
L}_{\tau}, \alpha^k)$ the \textsf{set of critical orbits}  of ${\cal
L}_\tau$ on $H_\tau(\alpha^k)$. It is always nonempty because ${\cal
L}_{\tau}$ can attain the minimal value on $H_{\tau}(\alpha^k)$.
Clearly,  we may assume that each  critical orbit of
 ${\cal L}_\tau$ on $H_\tau(\alpha^k)$ is isolated for any $k\in\N$.
  As in \S 5.1, by contradiction we
assume:\vspace{2mm}

 \noindent{\bf Assumption $FT(\alpha)$}: (i) {\it For any given
 integer $k>0$, the system (\ref{e:1.6}) only possesses  finitely many
 distinct,  $k\tau$-periodic  solution orbit towers based on $k\tau$-periodic
 solutions  of (\ref{e:1.6}) representing $\alpha^k$},
 (ii) {\it there exists an integer $k_0>1$ such that for each integer $k>k_0$, any
 $k\tau$-periodic
 solution $\tilde\gamma$ of the system (\ref{e:1.6}) representing $\alpha^k$
must be an iteration of some $l\tau$-periodic solution $\gamma$ of
the system (\ref{e:1.6}) representing $\alpha^l$ with $l\le k_0$ and
$k=lq$ for some $q\in\N$.}

Under this assumption, there only exist finitely many
 periodic solution orbit towers $\{s\cdot\hat\gamma_1^k\}^{s\in\R}_{k\in\N},\,\cdots,\,
 \{s\cdot\hat\gamma_p^k\}^{s\in\R}_{k\in\N}$
 of the system  (\ref{e:1.6}) such that\\
 $\bullet$ $\hat\gamma_i$ has period $k_i\tau\le k_0\tau$ and represents
 $\alpha^{k_i}$ for some $k_i\in\N$, $i=1,\cdots, p$;\\
 $\bullet$  for each
integer $k>k_0$ any  $k\tau$-periodic solution $\gamma$ of the
system (\ref{e:1.6}) representing $\alpha^k$  must be an iteration
of some $s\cdot\hat\gamma_i$, i.e.
$\gamma=(s\cdot\hat\gamma_i)^l=s\cdot \hat\gamma_i^l$  for some
$s\in\R$ and $l\in\N$ with $lk_i=k$.

Set $m:=k_0!$ (the factorial of $k_0$) and
$\gamma_i=\hat\gamma_i^{m/k_i}$, $i=1,\cdots, p$.  Then each
$\gamma_i$ is a $m\tau$-periodic solution of the system
(\ref{e:1.6}) representing $\alpha^m$. We conclude

\begin{claim}\label{cl:7.1}
 For any $k\in\N$, it holds that
$$
{\cal KO}({\cal L}_{km\tau},
\alpha^{km})=\bigl\{\,S_{km\tau}\cdot\gamma^k_j\,\bigm|\,1\le j\le
p\, \bigr\}.
$$
 \end{claim}

\noindent{\bf Proof}.\quad Let $\gamma\in {\cal K}({\cal
L}_{km\tau}, \alpha^{km})$. Since $km>k_0$, then
$\gamma=(s\cdot\hat\gamma_i)^l$ for some $s\in\R$ and $l\in\N$ with
$lk_i=km$. Hence $\gamma=s\cdot\hat\gamma_i^l=s\cdot
(\hat\gamma_i)^{km/k_i}=s\cdot
(\hat\gamma_i^{m/k_i})^k=s\cdot\gamma_i^k$. $\Box$

Hence replacing $\tau$ by $m\tau$ we may assume $m=1$ below, i.e.
\begin{equation}\label{e:7.1}
{\cal KO}({\cal L}_{k\tau},
\alpha^{k})=\bigl\{\,S_{k\tau}\cdot\gamma^k_j\,\bigm|\,1\le j\le p\,
\bigr\}\quad\forall k\in\N.
\end{equation}
 As in \S5.1 we can also \textbf{assume}: $\gamma_j^\ast TM\to
\R/\tau\Z$, $j=1,\cdots, p$, \textsf{are all trivial}.

\begin{lemma}\label{lem:7.2}
For each $k\in\N$ there exists  ${\cal O}_k\in{\cal KO}({\cal
L}_{k\tau}, \alpha^{k})$
 such that
$$
C_r({\cal L}_{k\tau}, {\cal O}_k; \K)\ne 0.
$$
Moreover, $r-2n\le r-1- m^0_{k\tau}({\cal O}_k)\le m^-_{k\tau}({\cal
O}_k)\le r-1$ if ${\cal O}_k$ is not a single point critical orbit,
and $r-2n\le r- m^0_{k\tau}({\cal O}_k)\le m^-_{k\tau}({\cal
O}_k)\le r$ otherwise.
\end{lemma}

\noindent{\bf Proof}.\quad By Lemma~\ref{lem:4.12} we only need to
prove the first claim. The proof is similar to that of
Lemma~\ref{lem:5.2}. Let $kc_1<\cdots <kc_l$ be all critical values
of ${\cal L}_{k\tau}$, $l\le p$, and $\inf{\cal L}_{k\tau}=kc_1$,
$k=1,2,\cdots$. As in the proof of Lemma~\ref{lem:5.2} we have a
large regular value $b$ of ${\cal L}_{k\tau}$ such that ${\rm rank}
H_r(({\cal L}_{k\tau})_b;\K)>0$.
 Take the regular values of
${\cal L}_{k\tau}$, $a_0<a_1<\cdots <a_l=b$ such that $kc_i\in
(a_{i-1}, a_i)$, $i=1,\cdots, l$. Noting (\ref{e:7.1}), by Theorem
2.1 of \cite{Wa} or the proof of Lemma 4 of \cite[pp. 502]{GM2}, we
get
$$
H_r(({\cal L}_{k\tau})_{a_i}, ({\cal L}_{k\tau})_{a_{i-1}}; \K)\cong
\bigoplus_{{\cal L}_{k\tau}(\gamma^k_j)=kc_i}C_r({\cal L}_{k\tau},
S_{k\tau}\cdot\gamma^k_j; \K).
$$
Since each critical point has finite Morse index, (\ref{e:4.67})
implies that each critical group $C_r({\cal L}_{k\tau},
S_{k\tau}\cdot\gamma^k_j; \K)$ has finite rank. Almost repeating the
proof of Lemma~\ref{lem:5.2} we get some $S_{k\tau}\cdot\gamma^k_j$
in ${\cal KO}({\cal L}_{k\tau}, \alpha^{k})$ such that ${\rm
rank}C_r({\cal L}_{k\tau}, S_{k\tau}\cdot\gamma^k_j; \K)>0$ and thus
${\rm rank}C_r({\cal L}_{k\tau}, S_{k\tau}\cdot\gamma^k_j; \K)\ne
0$.
  $\Box$\vspace{2mm}

 Corresponding to
Corollary~\ref{cor:5.4} we have

\begin{lemma}\label{lem:7.3}
Under Assumption FT($\alpha$), for every integer $q\ge n+1$ there
exists a constant $k_0(q)>0$ such that
$$
C_q({\cal L}_{k\tau}, {\cal O}_k; \K)=0
$$
for every integer $k\ge k_0(q)$ and ${\cal O}_k\in{\cal KO}({\cal
L}_{k\tau}, \alpha^{k})$. Here $k_0(q)=1$ if $\hat
m^-_{r}(\gamma_j)=0$ for all $1\le j\le p$, and
$$
k_0(q)=1+ {\rm max}\bigg\{\left[\frac{q + n}{\hat
m^-_{r}(\gamma_j)}\right]\,\Bigm|\,\hat m^-_{r}(\gamma_j)\ne 0,\,
1\le j\le p\biggr\}
$$ otherwise.
\end{lemma}

\noindent{\bf Proof.}\quad Let ${\cal
O}_k=S_{k\tau}\cdot\gamma_j^k$. If $\gamma_j$ is constant, by the
proof of Lemma~\ref{lem:5.3} we have
$$
C_q({\cal L}_{k\tau}, {\cal O}_k; \K)=C_q({\cal L}_{k\tau},
\gamma_j^k; \K)=0
$$
for any $k>k(q, \gamma_j)$, where $k(q, \gamma_j)=1$ if $\hat
m^-_\tau(\gamma_j)=0$, and $k(q, \gamma_j)= \frac{q+n}{\hat
m^-_\tau(\gamma_j)}$ if $\hat m^-_\tau(\gamma_j)\ne 0$.

 Suppose that $\gamma_j$ is not a constant solution. If
$C_q({\cal L}_{k\tau}, {\cal O}_k; \K)\ne 0$, Lemma~\ref{lem:4.12}
yields
\begin{equation}\label{e:7.2}
m^-_{k\tau}(S_{k\tau}\cdot\gamma_j^k)\le q-1\le
m^-_{k\tau}(S_{k\tau}\cdot\gamma_j^k)+
m^0_{k\tau}(S_{k\tau}\cdot\gamma_j^k).
\end{equation}
By (\ref{e:4.53}) this becomes
\begin{equation}\label{e:7.3}
m^-_{k\tau}(\gamma_j^k)\le q-1\le m^-_{k\tau}(\gamma_j^k)+
m^0_{k\tau}(\gamma_j^k)-1.
\end{equation}
If $\hat m_{\tau}(\gamma_j)>0$, it follows from (\ref{e:7.3}) and
(\ref{e:3.2}) that
$$
k\hat m^-_{\tau}(\gamma_j)-n\le m^-_{k\tau}(\gamma^k_j)\le q-1
$$
and therefore $k\le \frac{q + n-1}{\hat m^-_{r}(\gamma_j)}$. This
contradicts to $k\ge k_0(q)$. If $\hat m_{\tau}(\gamma_j)=0$, by
(\ref{e:3.2}),
$$
0\le m^-_{k\tau}(\gamma^k_j)\le n- m^0_{k\tau}(\gamma^k_j)
\quad\forall k\in\N.
$$
It follows that
$$
m^-_{k\tau}(S_{k\tau}\cdot\gamma_j^k)+
m^0_{k\tau}(S_{k\tau}\cdot\gamma_j^k)=m^-_{k\tau}(\gamma_j^k)+
m^0_{k\tau}(\gamma_j^k)-1\le n-1.
$$
Since $q\ge n+1$, (\ref{e:7.2}) implies that
$m^-_{k\tau}(S_{k\tau}\cdot\gamma_j^k)+
m^0_{k\tau}(S_{k\tau}\cdot\gamma_j^k)\ge n$. This also gives a
contradiction. Lemma~\ref{lem:7.3} is proved. $\Box$\vspace{2mm}

Clearly,  Lemma~\ref{lem:7.2} and Lemma~\ref{lem:7.3} imply
Theorem~\ref{th:1.1}(i) in the case  $r\ge n+1$.

 \verb"In the following we"
\verb"consider the case" $r=n$.

Under Assumption $FT(\alpha)$ we apply Lemma \ref{lem:7.2} to all
$k\in \{2^m\,|\,m\in\{0\}\cup{\N}\}$ to get an infinite subsequence
$Q$ of $\{2^m\,|\,m\in\{0\}\cup{\N}\}$, some $l\in\N$ and an
$\gamma\in\{\gamma_1,\cdots, \gamma_p\}$ such that $C_n({\cal
L}_{kl\tau}, S_{kl\tau}\cdot\gamma^{kl}; \K)\ne 0$,
 $m^-_{kl\tau}(S_{kl\tau}\cdot\gamma^{kl})=
m^-_{l\tau}(S_{l\tau}\cdot\gamma^{l})$ and
$m^0_{kl\tau}(S_{kl\tau}\cdot\gamma^{kl})=m^0_{l\tau}(S_{l\tau}\cdot\gamma^{l})$
for any $k\in Q$.  As before we {\bf always assume $l=1$} in the
following. Then we have
\begin{equation}\label{e:7.4}
\left.\begin{array}{ll}
 &C_n({\cal L}_{k\tau}, S_{k\tau}\cdot\gamma^{k}; \K)\ne 0\quad{\rm and}\\
 &m^-_{k\tau}(S_{k\tau}\cdot\gamma^{k})=
m^-_{\tau}(S_{\tau}\cdot\gamma),\quad
m^0_{k\tau}(S_{k\tau}\cdot\gamma^{k})=m^0_{\tau}(S_{\tau}\cdot\gamma)
\end{array}\right\}
\end{equation}
for any $k\in Q$.   By Lemma~\ref{lem:7.3} there exists $k_0>0$ such
that for any $\gamma\in\{\gamma_1,\cdots, \gamma_p\}$,
\begin{equation}\label{e:7.5}
C_{n+1}({\cal L}_{k\tau},  S_{k\tau}\cdot\gamma^k; \K)=0\quad\forall
k\in Q(k_0):=\{k\in Q\,|\, k\ge k_0\}.
\end{equation}

Denote by ${\cal O}=S_\tau\cdot\gamma$, and by $c={\cal
L}_\tau(\gamma)={\cal L}_\tau({\cal O})$. Under Assumption
FT($\alpha$), as in \S5.1 let us take $\nu>0$ sufficiently small
so that for each $k\in\N$ the interval $[k(c-3\nu), k(c+ 3\nu)]$
contains an unique critical value $kc$ of ${\cal L}_{k\tau}$ on
$H_{k\tau}(\alpha^{k})$, i.e.
$$
{\cal L}_{k\tau}\bigr({\cal KO}({\cal L}_{k\tau},
\alpha^{k})\bigl)\cap [k(c-3\nu), k(c+ 3\nu)]=\{kc\}.
$$
For any $k\in Q$, by Theorem~\ref{th:4.11}, we may choose a
topological Gromoll-Meyer pair of ${\cal L}_\tau$ at ${\cal
O}\subset H_\tau(\alpha)$, $(\widehat W({\cal O}), \widehat W({\cal
  O})^-)$ satisfying
\begin{equation}\label{e:7.6}
 \bigl(\widehat W({\cal O}), \widehat W({\cal
  O})^-\bigr)\subset \left(({\cal
L}_\tau)^{-1}([c-2\nu, c+ 2\nu]), ({\cal
L}_\tau)^{-1}(c-2\nu)\right),
\end{equation}
and a topological Gromoll-Meyer pair  of ${\cal L}_{k\tau}$ at
$\psi^k({\cal O})\subset H_{k\tau}(\alpha^k)$,
$$
\bigl(\widehat
W(\psi^k({\cal O})), \widehat W(\psi^k({\cal
  O}))^-\bigr)
  $$
  such that
\begin{eqnarray}
&& \bigl(\psi^k(\widehat W({\cal O})), \psi^k(\widehat W({\cal
  O})^-)\bigr)\subset \bigl(\widehat W(\psi^k({\cal O})), \widehat W(\psi^k({\cal
  O}))^-\bigr)\quad{\rm and}\label{e:7.7}\\
&&\bigl(\widehat W(\psi^k({\cal O})), \widehat W(\psi^k({\cal
  O}))^-\bigr)\subset \nonumber\\
 && \hspace{20mm} \bigl(({\cal
L}_{k\tau})^{-1}([kc-2k\nu, kc+ 2k\nu]), ({\cal
L}_{k\tau})^{-1}(kc-2k\nu)\bigr)\label{e:7.8}
\end{eqnarray}
and that the iteration map $\psi^k: H_{\tau}(\alpha)\to
H_{k\tau}(\alpha^k)$ induces an isomorphism:
\begin{eqnarray*}
 &&\psi^k_*:  C_*({\cal L}_{\tau}, {\cal O}; \K):=H_\ast\bigl(\widehat W({\cal O}), \widehat W({\cal
  O})^-; \K\bigr)\\
 &&\qquad \longrightarrow
         C_*({\cal L}_{k\tau}, \psi^k({\cal O}); \K):
         =H_\ast\bigl(\widehat W(\psi^k({\cal O})), \widehat W(\psi^k({\cal
  O}))^-; \K\bigr).
\end{eqnarray*}
Identifying $\psi({\cal O})={\cal O}$,  for $j=1, k$, denote by
the inclusions
\begin{eqnarray*}
&&{\bf h}_1^j:\bigl(\widehat W(\psi^j({\cal O})), \widehat
W(\psi^j({\cal
  O}))^-\bigr)\hookrightarrow
\bigl(({\cal L}_{j\tau})_{j(c+ 2\nu)},
({\cal L}_{j\tau})_{j(c-2\nu)}\bigr),\\
&&{\bf h}_2^j:\bigl(({\cal L}_{j\tau})_{j(c+ 2\nu)}, ({\cal
L}_{j\tau})_{j(c-2\nu)}\bigr) \hookrightarrow \bigl(({\cal
L}_{j\tau})_{j(c+ 2\nu)}, ({\cal
L}_{j\tau})^\circ_{j(c-\nu)}\bigr),\\
&&{\bf h}_3^j:  \bigl(({\cal L}_{j\tau})_{j(c+ 2\nu)}, ({\cal
L}_{j\tau})^\circ_{j(c-\nu)}\bigr)\hookrightarrow \bigl(H_{j\tau},
({\cal L}_{j\tau})^\circ_{j(c-\nu)}\bigr).
\end{eqnarray*}
As in \S5.1 we have monomorphisms on homology modules,
\begin{eqnarray*}
&&({\bf h}_2^j\circ{\bf h}_1^j)_\ast:H_\ast\bigl(\widehat
W(\psi^j({\cal O})), \widehat W(\psi^j({\cal
  O}))^-;\K\bigr)\to H_\ast\bigl(({\cal L}_{j\tau})_{j(c+
2\nu)}, ({\cal
L}_{j\tau})^\circ_{j(c-\nu)};\K\bigr),\\
&&({\bf h}_3^j)_\ast: H_\ast\bigl(({\cal L}_{j\tau})_{j(c+ 2\nu)},
({\cal L}_{j\tau})^\circ_{j(c-\nu)};\K\bigr)\to
H_\ast\bigl(H_{j\tau}, ({\cal L}_{j\tau})^\circ_{j(c-\nu)};\K\bigr).
\end{eqnarray*}
Moreover, the inclusions
\begin{eqnarray*}
&&\!\!\!\!\!\!\!{\bf I}_j:\bigl(\widehat W(\psi^j({\cal O})),
\widehat W(\psi^j({\cal
  O}))^-\bigr)\hookrightarrow
\bigl(({\cal L}_{j\tau})^{-1}([jc-2j\nu, jc+ 2j\nu]),
({\cal L}_{j\tau})^{-1}(jc-2j\nu)\bigr),\\
&&\!\!\!\!\!\!\! {\bf J}_j:\bigl(({\cal
L}_{j\tau})^{-1}([jc-2j\nu, jc+ 2j\nu]), ({\cal
L}_{j\tau})^{-1}(jc-2j\nu)\bigr)\hookrightarrow \bigl(H_{j\tau},
({\cal L}_{j\tau})^\circ_{jc-j\nu}\bigr)
\end{eqnarray*}
satisfy
\begin{equation}\label{e:7.9}
{\bf J}_j\circ{\bf I}_j={\bf h}^j_3\circ{\bf h}^j_2\circ{\bf
h}^j_1,\quad j=1,k.
\end{equation}
By (\ref{e:7.7}), we have also
\begin{equation}\label{e:7.10}
\psi^k\circ{\bf I}_1={\bf I}_k\circ\psi^k
\end{equation}
as maps from $\bigl(\widehat W({\cal O}), \widehat W({\cal
  O})^-\bigr)$ to
$\bigl(({\cal L}_{k\tau})^{-1}([kc-2k\nu, kc+ 2k\nu]), ({\cal
L}_{k\tau})^{-1}(kc-2k\nu)\bigr)$.
 These yield the following corresponding result with Proposition~\ref{prop:5.5}.

\begin{proposition}\label{prop:7.4}
Under  Assumption $FT(\alpha)$, there exist a $\tau$-periodic
solution $\gamma$ of (\ref{e:1.6}) representing $\alpha$, a large
integer $k_0>0$, an infinite integer set $Q$ containing $1$, and a
small $\epsilon>0$ having properties: For the orbit ${\cal
O}=S_\tau\cdot\gamma$ and any $k\in Q(k_0):=\{k\in Q\,|\, k\ge
k_0\}$ there exist topological Gromoll-Meyer pairs $\bigl(\widehat
W({\cal O}), \widehat W({\cal
  O})^-\bigr)$ and $\bigl(\widehat W(\psi^k({\cal O})), \widehat W(\psi^k({\cal
  O}))^-\bigr)$ satisfying (\ref{e:7.6})-(\ref{e:7.8}) such
that for the inclusion
$$
j_{k\tau}={\bf h}^k_3\circ {\bf h}^k_2\circ {\bf h}^k_1:
\bigl(\widehat W(\psi^k({\cal O})), \widehat W(\psi^k({\cal
  O}))^-\bigr)\to
\bigl(H_{k\tau}(\alpha^{k}), ({\cal
L}_{k\tau})^\circ_{k(c-\nu)}\bigr)
$$
 the following diagram holds:
\begin{eqnarray}\label{e:7.11}
 &&0\ne C_n({\cal L}_{\tau},
{\cal O};\K)\stackrel{\psi^k_*}{\longrightarrow} C_n({\cal
L}_{k\tau},
\psi^k({\cal O});\K)\nonumber\\
&&\hspace{30mm}\stackrel{(j_{k\tau})_*}{\longrightarrow}
H_n\bigl(H_{k\tau}(\alpha^{k}), ({\cal
L}_{k\tau})^\circ_{k(c-\nu)};\K\bigr)\equiv{\cal H}_k,
 \end{eqnarray}
 where  $c={\cal L}_\tau(\gamma)$, $\psi^k_*$ is an isomorphism, and $(j_{k\tau})_*$ is a
monomorphism among the singular homology modules. In particular,
if $\omega$ is a generator of $C_n({\cal L}_\tau, {\cal
O};\K)=H_n\bigl(\widehat W({\cal O}), \widehat W({\cal
  O})^-; \K\bigr)$, then
\begin{eqnarray}
&&(j_{k\tau})_*\circ(\psi^k)_*(\omega)\ne 0\quad{\rm in}\;
      {\cal H}_k,\label{e:7.12}\\
&&(j_{k\tau})_\ast\circ(\psi^k)_\ast(\omega)=({\bf
J}_k)_\ast\circ({\bf I}_k)_\ast\circ(\psi^k)_\ast(\omega)\nonumber\\
&&\hspace{29mm}=({\bf J}_k)_\ast\circ(\psi^k)_\ast\circ ({\bf
I}_1)_\ast(\omega)\quad{\rm in}\;
      {\cal H}_k.\label{e:7.13}
\end{eqnarray}
\end{proposition}

Now we can slightly modify the arguments from
Proposition~\ref{prop:5.6} to Proposition~\ref{prop:5.8} to complete
the proof of (i). The only place which should be noted is that for
$\psi^k({\cal O})$ in (\ref{e:7.11}) Lemma~\ref{lem:4.13} implies
each point $y\in \psi^k({\cal O})$ to be a non-minimum saddle point
of ${\cal L}_{k\tau}$ on $H_{k\tau}(\alpha^k)$ in the case $\dim
M=n>1$. \vspace{2mm}

\noindent{\bf 7.2. Proof of (ii)} can be completed  by the similar
arguments as in \S5.2.

\section{Questions and remarks}
\setcounter{equation}{0}

 For a $C^3$-smooth compact n-dimensional
manifold $M$ without boundary, and a $C^2$-smooth map  $H:\R\times
T^\ast M\to\R$ satisfying the conditions (H1)-(H5),  we have shown
in $1^\circ)$ of Theorem~\ref{th:1.12} that the Poincar\'{e} map
$\Psi^H$ has infinitely many distinct periodic points sitting in the
zero section $0_{T^\ast M}$ of $T^\ast M$. Notice that the condition
(H5) can be expressed as: $H(t, x)=H(-t, \tau_0(x))\;\forall
(t,x)\in\R\times M$, where  $\tau_0: T^\ast M\to T^\ast M,\;(q,
p)\mapsto (q, -p)$, is the standard anti-symplectic involution. So
it is natural to consider the following question: Let
$(P,\omega,\tau)$ be a {\bf real symplectic manifold} with  an
anti-symplectic involution $\tau$ on $(P,\omega)$, i.e.
$\tau^\ast\omega=-\omega$ and $\tau^2=id_P$. A smooth time dependent
Hamiltonian function $H:\R\times P\to\R, (t,x)\mapsto H(t,x)=H_t(x)$
is said to be {\bf $1$-periodic in time and symmetric} if it
satisfies
$$
H_t(x)=H_{t+1}(x)\quad{\rm and}\quad H(t, x)=H(-t, \tau(x))\;\forall
(t,x)\in\R\times P.
$$
In this case, the Hamiltonian vector fields $X_{H_t}$ satisfies
$X_{H_{t+1}}(x)=X_{H_{t}}(x)=-d\tau(\tau(x))X_{H_{-t}}(\tau(x))$ for
all $(t,x)\in\R\times P$. If  the global flow  of
\begin{equation}\label{e:8.1}
 \dot
x(t)=X_{H_{t}}(x(t))
\end{equation}
exists, denoted by $\Psi^H_t$,   then it is obvious that
 $$
\Psi^H_{t+1}=\Psi^H_t\circ\Psi^H_1\;\forall t\in\R,\quad
\Psi^H_1\circ\tau=\tau\circ(\Psi^H_1)^{-1}.
$$
So each $\tau$-invariant $k$-periodic point $x_0$, i.e.
$\tau(x_0)=x_0$,  of $\Psi^H=\Psi^H_1$ with $k\in\N$ yields a
$k$-periodic contractible solution $x(t)=\Psi^H_t(x_0)$ of
(\ref{e:8.1}) satisfying $x(-t)=\tau(x(t))$ for all $t\in\R$. Such a
solution is called $\tau$-{\bf reversible}. By [Vi1, p.4] the fixed
point set $L:={\rm Fix}(\tau)$ of $\tau$ is either empty or a
Lagrange submanifold. It is natural to ask the following more
general version of the Conley conjecture.\vspace{2mm}

\begin{question}\label{qu:8.1}
{\rm Suppose that $L$ is nonempty and compact, and that $(P,
\omega)$ satisfies some good condition (e.x.  geometrically bounded
for some $J\in \R{\cal J}(P,\omega):=\{J\in {\cal J}(P,\omega)\,|\,
J\circ d\tau=-d\tau\circ J\}$ and  Riemannian metric $\mu$ on $P$).
Has the system (\ref{e:8.1}) infinitely many distinct
 $\tau$-reversible contractible periodic  solutions of integer periods?
  Furthermore, if the flow $\Psi^H_t$ exists globally,
   has the Poincar\'{e} map $\Psi^H=\Psi^H_1$ infinitely many
   distinct periodic points sitting in $L$?}
\end{question}

Let ${\cal P}_0(H,\tau)$ denote the set of all contractible
$\tau$-reversible $1$-periodic solutions of (\ref{e:8.1}). Since the
Conley conjecture came  from the Arnold conjecture,
Question~\ref{qu:8.1} naturally suggests the following more general
versions of the Arnold conjectures.

\begin{question}\label{qu:8.2}
{\rm Under the assumptions of Question~\ref{qu:8.1}, $\sharp{\cal
P}_0(H,\tau)\ge {\rm Cuplength}_{\F}(L)$ for $\F=\Z, \Z_2$?
Moreover, if some nondegenerate assumptions for elements of ${\cal
P}_0(H,\tau)$ are satisfied, $\sharp{\cal P}_0(H,\tau)\ge \sum^{\dim
L}_{k=0}b_k(L,\F)$? }
\end{question}

This question is closely related to the Arnold-Givental conjecture,
cf. \cite{Lu1}. In order to study it we try to construct a real
Floer homology $FH_\ast(P, \omega, \tau, H)$ with ${\cal
P}_0(H,\tau)$ under some nondegenerate assumptions for elements of
${\cal P}_0(H,\tau)$, which is expected to be isomorphic to
$H_\ast(M)$. Moreover, if $L\in C^2(\R/Z\times TM)$ satisfies
(L1)-(L4) and the functional ${\cal L}(\gamma)=\int^1_0L(t,
\gamma(t), \dot\gamma(t))dt$ on $EH_1$ has only nondegenerate
critical points, then one can, as in \cite[\S 2.2]{AbSc}, construct
a Morse complex $CM_\ast({\cal L})$ whose homology is isomorphic to
$H_\ast(M)$ as well. As in \cite{Vi3, SaWe, AbSc}, it is also
natural to construct an isomorphism between $HF_\ast(T^\ast M,
\omega_{\rm can}, \tau_0, H)$ and $H(CM_\ast({\cal L}))$ and to
study different product operations in them.

  The author believes that the techniques developed in
this paper are useful for one to generalize the results of multiple
periodic solutions of some Lagrangian and Hamiltonian systems on the
Euclidean space to manifolds.

\section{Appendix}
\setcounter{equation}{0}

\noindent{\bf A.1. Proof of Proposition A.}\quad The first claim is
a direct consequence of the following (\ref{e:A.4}). As to the
second, since for each $t\in\R$ the functions $L_t=L(t, \cdot)$ and
$H_t=H(t, \cdot)$ are Fenchel transformations of each other, we only
need to prove that (H2)-(H3) can be satisfied under the assumptions
(L2)-(L3).
 For conveniences we omit the time variable $t$. In any
local coordinates $(q_1,\cdots, q_n)$, we write $(q, v)=(q_1,\cdots,
q_n, v_1,\cdots, v_n)$. By definition of $H$ we have
\begin{equation}\label{e:A.1}
 H\bigl(q, \frac{\partial L}{\partial
v}(q,v)\bigr)=-L(q, v)+ \sum^n_{j=1}\frac{\partial L}{\partial
v_j}(q,v)v_j.
\end{equation}
Differentiating both sides with respect to the variable $v_i$ we get
$$
\sum^n_{j=1}\frac{\partial H}{\partial p_j}\bigl(q, \frac{\partial
L}{\partial v}(q,v)\bigr)\frac{\partial^2 L}{\partial v_i\partial
v_j}(q,v)= \sum^n_{j=1}v_j\frac{\partial^2 L}{\partial v_i\partial
v_j}(q, v).
$$
Since the matrix $\Bigl[\frac{\partial^2 L}{\partial v_i\partial
v_j}(q, v)\Bigr]$ is invertible, it follows that
\begin{equation}\label{e:A.2}
 \frac{\partial H}{\partial p_j}\bigl(q, \frac{\partial L}{\partial v}(q,v)\bigr)=v_j.
\end{equation}
Let $p=\frac{\partial L}{\partial v}(q,v)$. Differentiating both
sides of (\ref{e:A.1}) with respect to the variable $q_i$ and using
(\ref{e:A.2}) we obtain
\begin{eqnarray*}
 & &\sum^n_{j=1}v_j\frac{\partial^2 L}{\partial
q_i\partial v_j}(q, v)- \frac{\partial L}{\partial q_i}(q, v)\\
 &=&\frac{\partial H}{\partial q_i}\bigl(q, \frac{\partial L}{\partial
v}(q,v)\bigr)+  \sum^n_{j=1}\frac{\partial H}{\partial
p_j}\bigl(q,\frac{\partial L}{\partial v}(q,v)\bigr)
\frac{\partial^2 L}{\partial
q_i\partial v_j}(q,v)\\
&=& \frac{\partial H}{\partial q_i}\bigl(q, \frac{\partial
L}{\partial v}(q,v)\bigr)+  \sum^n_{j=1}v_j \frac{\partial^2
L}{\partial q_i\partial v_j}(q,v)
\end{eqnarray*}
and hence
\begin{equation}\label{e:A.3}
 \frac{\partial H}{\partial q_i}\bigl(q, \frac{\partial L}{\partial v}(q, v)\bigr)=
 -\frac{\partial L}{\partial q_i}(q, v).
\end{equation}
Differentiating both sides of (\ref{e:A.2}) with respect to the
variable $v_i$ yields
\begin{eqnarray}\label{e:A.4}
\sum^n_{k=1}\frac{\partial^2 H}{\partial p_j\partial p_k}(q, p)
\frac{\partial^2 L}{\partial v_k\partial v_i}(q, v)=\delta_{ij},\;
i.e.\nonumber\\
\Bigl[\frac{\partial^2 H}{\partial p_i\partial p_j}(q,
p)\Bigr]=\Bigl[\frac{\partial^2 L}{\partial v_i\partial v_j}(q,
v)\Bigr]^{-1}.
\end{eqnarray}
Differentiating both sides of (\ref{e:A.2}) with respect to the
variable $q_i$, and both sides of (\ref{e:A.3}) with respect to the
variable $q_j$ respectively, we arrive at
\begin{eqnarray*}
&&\frac{\partial^2 H}{\partial p_j\partial q_i}\bigl(q,
\frac{\partial L}{\partial v}(q,v)\bigr)+
\sum^n_{k=1}\frac{\partial^2 H}{\partial p_j\partial p_k}\bigl(q,
\frac{\partial L}{\partial v}(q,v)\bigr)\frac{\partial^2 L}{\partial
v_k\partial q_i}(q,v)=0,\\
&&\frac{\partial^2 H}{\partial q_i\partial q_j}\bigl(q,
\frac{\partial L}{\partial v}(q,v)\bigr)+
\sum^n_{k=1}\frac{\partial^2 H}{\partial q_i\partial p_k}\bigl(q,
\frac{\partial L}{\partial v}(q,v)\bigr)\frac{\partial^2 L}{\partial
v_k\partial q_j}(q,v)= -\frac{\partial^2 L}{\partial q_i\partial
q_j}(q, v),
\end{eqnarray*}
or their equivalent expressions of matrixes,
\begin{eqnarray*}
&&\Bigl[\frac{\partial^2 H}{\partial p_i\partial q_j}\bigl(q,
\frac{\partial L}{\partial v}(q,v)\bigr)\Bigr]+
\Bigl[\frac{\partial^2 H}{\partial p_i\partial p_j}\bigl(q,
\frac{\partial L}{\partial
v}(q,v)\bigr)\Bigr]\,\Bigl[\frac{\partial^2 L}{\partial
v_i\partial q_j}(q,v)\Bigr]=0,\\
&&\Bigl[\frac{\partial^2 H}{\partial q_i\partial q_j}\bigl(q,
\frac{\partial L}{\partial v}(q,v)\bigr)\Bigr]+
\Bigl[\frac{\partial^2 H}{\partial p_i\partial q_j}\bigl(q,
\frac{\partial L}{\partial
v}(q,v)\bigr)\Bigr]^t\,\Bigl[\frac{\partial^2 L}{\partial
v_i\partial q_j}(q,v)\Bigr] = -\Bigl[\frac{\partial^2 L}{\partial
q_i\partial q_j}(q, v)\Bigr].
\end{eqnarray*}
It follows from these that
\begin{eqnarray}\label{e:A.5}
&&\quad\Bigl[\frac{\partial^2 L}{\partial q_i\partial q_j}(q,
v)\Bigr]\nonumber\\
\!\!\!\!\!\!\!\!&& = \Bigl[\frac{\partial^2 H}{\partial p_i\partial
q_j}\bigl(q, \frac{\partial L}{\partial
v}(q,v)\bigr)\Bigr]^t\,\Bigl[\frac{\partial^2 H}{\partial
p_i\partial p_j}\bigl(q, \frac{\partial L}{\partial
v}(q,v)\bigr)\Bigr]^{-1}\,\Bigl[\frac{\partial^2 H}{\partial
p_i\partial q_j}\bigl(q, \frac{\partial L}{\partial
v}(q,v)\bigr)\Bigr]\nonumber\\
 \!\!\!\!\!\!\!\!&&-\Bigl[\frac{\partial^2 H}{\partial q_i\partial
q_j}\bigl(q, \frac{\partial L}{\partial v}(q,v)\bigr)\Bigr].
\end{eqnarray}
Finally, differentiating both sides of (\ref{e:A.3}) with respect to
the variable $v_j$ we get
\begin{eqnarray}\label{e:A.6}
\frac{\partial^2 L}{\partial q_i\partial v_j}(q,
v)\!\!\!\!\!\!\!\!&&= -\sum^n_{k=1}\frac{\partial^2 H}{\partial
q_i\partial p_k}\bigr(q, \frac{\partial L}{\partial v}(q,v)\bigl)
\frac{\partial^2 L}{\partial v_k\partial v_j}(q, v),\;
i.e.\nonumber\\
\Bigl[\frac{\partial^2 L}{\partial q_i\partial v_j}(q,
v)\Bigr]\!\!\!\!\!\!\!\!&&= - \Bigl[\frac{\partial^2 H}{\partial
q_i\partial p_j}\bigl(q, \frac{\partial L}{\partial
v}(q,v)\bigr)\Bigr]\,
\Bigl[\frac{\partial^2 L}{\partial v_i\partial v_j}(q, v)\Bigr]\nonumber\\
\!\!\!\!\!\!\!\!&&= - \Bigl[\frac{\partial^2 H}{\partial q_i\partial
p_j}\bigl(q, \frac{\partial L}{\partial v}(q,v)\bigr)\Bigr]\,
\Bigl[\frac{\partial^2 H}{\partial p_i\partial p_j}\bigl(q,
\frac{\partial L}{\partial v}(q,v)\bigr)\Bigr]^{-1}.
 \end{eqnarray}
Here the final equality is due to (\ref{e:A.4}). Since
$p=\frac{\partial L}{\partial v}(q,v)$ and $v=\frac{\partial
H}{\partial p}(q, p)$,  the desired conclusions will follow from
(\ref{e:A.4})-(\ref{e:A.6}). Indeed, by (\ref{e:A.4}) it is easily
seen that (L2) is equivalent to
\begin{description}
\item[(H2')] $\sum_{ij}\frac{\partial^2H}{\partial p_i\partial p_j}(t,q,
p)u_iu_j\le \frac{1}{c}|{\bf u}|^2\quad\forall {\bf
u}=(u_1,\cdots,u_n)\in\R^n$.
\end{description}
Moreover, the three inequalities in (L3) have respectively the
following equivalent versions in terms of matrix norms:
\begin{eqnarray*}
 &&\Biggl|\biggl[ \frac{\partial^2 L}{\partial q_i\partial
q_j}(t,q,v)\biggr]\Biggr|\le C(1+ |v|^2),\quad \Biggl|
\biggl[\frac{\partial^2 L}{\partial q_i\partial
v_j}(t,q,v)\biggr]\Biggr|\le C(1+
|v|)\\
&&\quad\hbox{and}\qquad \Biggl| \biggl[\frac{\partial^2L}{\partial
v_i\partial v_j}(t,q,v)\biggr]\Biggr|\le C.
\end{eqnarray*}
Then (L3)  is equivalent to the following
\begin{description}
\item[(H3')] $\biggl|\Bigl[\frac{\partial^2 H}{\partial
p_i\partial q_j}\bigl(t, q, p\bigr)\Bigr]^t\,\Bigl[\frac{\partial^2
H}{\partial p_i\partial p_j}\bigl(t, q,
p\bigr)\Bigr]^{-1}\,\Bigl[\frac{\partial^2 H}{\partial p_i\partial
q_j}\bigl(t, q, p\bigr)\Bigr]-\Bigl[\frac{\partial^2 H}{\partial
q_i\partial q_j}\bigl(t, q, p\bigr)\Bigr]\biggr|\\
\le C\biggl(1+ \Bigl|\frac{\partial H}{\partial p}(t, q,
p)\Bigr|^2\biggr)$, \\
$\biggl|\Bigl[\frac{\partial^2 H}{\partial q_i\partial p_j}\bigl(t,
q, p\bigr)\Bigr]\, \Bigl[\frac{\partial^2 H}{\partial p_i\partial
p_j}\bigl(t, q, p\bigr)\Bigr]^{-1}\biggl| \le C\biggl(1+
\Bigl|\frac{\partial H}{\partial p}(t, q, p)\Bigr|\biggr)$,\quad and\vspace{4mm}\\
 $\biggl|\Bigl[ \frac{\partial^2 H}{\partial p_i\partial
p_j}(t,q, p)\Bigr]^{-1}\biggr|\le C$.\\
 Here  $\frac{\partial H}{\partial p}(t,
q, p)=\Bigr(\frac{\partial H}{\partial p_1}(t, q, p),\cdots,
\frac{\partial H}{\partial p_n}(t, q, p)\Bigl)$, and $|A|$ denotes
the standard norm of matrix $A\in\R^{n\times n}$, i.e.
$|A|=(\sum^n_{i=1}\sum^n_{j=1}a_{ij}^2)^{1/2}$ if $A=(a_{ij})$.
\end{description}

Note that $|A|=\sup_{|x|=1}|(Ax, x)_{\R^n}|$ for any symmetric
matrix $A\in\R^{n\times n}$, and $|A|=\sup_{|x|=1}(Ax, x)_{\R^n}$ if
$A$ is also positive definite,
 where $(\cdot,\cdot)_{\R^n}$ is
the standard inner product in $\R^n$. As usual, for two symmetric
positive matrixes  $A, B\in\R^{n\times n}$, by ``$A\le B$'' we mean
that $(Ax, x)_{\R^n}\le (Bx, x)_{\R^n}$ for any $x\in\R^n$. Then it
is easily proved that
\begin{equation}\label{e:A.7}
\biggl|\Bigl[ \frac{\partial^2H}{\partial p_i\partial p_j}(t,q,
p)\Bigr]^{-1}\biggr|\le C \Longleftrightarrow\Bigl[
\frac{\partial^2H}{\partial p_i\partial p_j}(t,q, p)\Bigr]\ge
\frac{1}{C}I_n.
\end{equation}
This and (H2') yield
$$
 \frac{1}{C}I_n\le\Bigl[ \frac{\partial^2H}{\partial p_i\partial
p_j}(t,q, p)\Bigr]\le \frac{1}{c}I_n.
$$

\noindent{\bf Lemma A.1.}\quad{\it For a matrix $B\in\R^{n\times n}$
and symmetric matrixes $A, B\in\R^{n\times n}$, suppose that there
exist constants
$0<c<C$ and $\alpha\ge 0$ such that\\
\noindent{\rm (i)} $\frac{1}{C}I_n\le A\le\frac{1}{c}I_n$,\\
\noindent{\rm (ii)} $|BA^{-1}|\le C(1+ \alpha)$,\\
\noindent{\rm (iii)} $|B^tA^{-1}B-E|\le C(1+ \alpha^2)$.\\
Then it holds that
\begin{equation}\label{e:A.8}
|B|\le \frac{C}{c}(1+ \alpha)\quad{\rm and}\quad |E|\le
\Bigl(\frac{2C^3}{c^2}+ C\Bigr)(1+ \alpha^2).
\end{equation}
Conversely, if (i) and (\ref{e:A.8}) are satisfied, then
\begin{equation}\label{e:A.9}
|BA^{-1}|\le \frac{C^2}{c}(1+ \alpha)\quad{\rm and}\quad
|B^tA^{-1}B-E|\le \Bigl(\frac{4C^3}{c^2}+ C\Bigr)(1+ \alpha^2).
\end{equation} }

\noindent{\bf Proof.}\quad By (i), $|A|\le\frac{1}{c}$ and
$|A^{-1}|\le C$. Hence
\begin{eqnarray*}
&&|B|=|BA^{-1}A|\le |BA^{-1}||A|\le  \frac{C}{c}(1+ \alpha),\\
&&|E|=|B^tA^{-1}B-E-B^tA^{-1}B|\le |B^tA^{-1}B-E|+ |B^tA^{-1}B|\\
&&\hspace{6mm}\le  C(1+ \alpha^2)+ |B|^2|A^{-1}|\le  C(1+ \alpha^2)+
\frac{C^3}{c^2}(1+ \alpha)^2\\
&&\hspace{6mm}\le C(1+ \alpha^2)+ \frac{2C^3}{c^2}(1+
\alpha^2)\\
&&\hspace{6mm}\le (C+ \frac{2C^3}{c^2})(1+ \alpha^2).
\end{eqnarray*}
(\ref{e:A.8}) is proved. The ``conversely'' part is easily proved as
well. $\Box$\vspace{2mm}

By this lemma we get immediately:

\noindent{\bf Proposition A.2.}\quad{\it In any local coordinates
$(q_1,\cdots, q_n)$, the conditions (L2)-(L3) are equivalent to the
fact that there exist constants $0<C_1<C_2$, depending on the local
coordinates,  such that
\begin{eqnarray*}
&&C_1I_n\le \Bigl[\frac{\partial^2 H}{\partial p_i\partial p_j}(t,q,
p)\Bigr] \le
C_2I_n,\label{e:A.10}\\
&&\biggl|\Bigl[\frac{\partial^2 H}{\partial q_i\partial p_j}\bigl(t,
q, p\bigr)\Bigr]\biggr| \le C_2\Bigl(1+ \Bigl|\frac{\partial
H}{\partial p}(t, q, p)\Bigr|\Bigr),\label{e:A.11}\\
&& \biggl|\Bigl[\frac{\partial^2 H}{\partial q_i\partial
q_j}\bigl(t, q, p\bigr)\Bigr]\biggr| \le C_2\Bigl(1+
\Bigl|\frac{\partial H}{\partial p}(t, q,
p)\Bigr|^2\Bigr).\label{e:A.12}
\end{eqnarray*} }

For each $(t, q)\in\R/\Z\times M$, since the function $T_q^\ast
M\to\R,\;p\mapsto H(t, q, p)$ is strictly convex, it has a unique
minimal point $\bar p=\bar p(t, q)$. In particular, $D_pH(t, q, \bar
p)=0$. Recall that the diffeomorphism $\frak{L}_H$ in (\ref{e:1.3})
is the inverse of $\frak{L}_L$ in (\ref{e:1.5}), and that $L(t,q, v)
= \langle p(t,q, v), v\rangle- H(t,q, p(t,q, v))$, where
$p=p(t,q,v)$ is a unique point determined by the equality
$v=D_pH(t,q, p)$.  It follows that
$$
\{(t, q, \bar p(t, q))\in \R/\Z\times T^\ast M\,|\, (t,
q)\in\R/\Z\times M\}=\frak{L}_H(\R/\Z\times 0_{TM})
$$
is a compact
subset. So in any local coordinates $(q_1,\cdots, q_n)$,  there
exists a constant $C_3>0$, depending on the local coordinates,  such
that the expression of  $\bar p=\bar p(t, q)$ in the local
coordinate $(q_1,\cdots, q_n)$, denoted by $\bar p=(\bar p_1,\cdots,
\bar p_n)$, satisfies
\begin{equation}\label{e:A.10}
|\bar p|=|(\bar p_1,\cdots, \bar p_n)|\le C_3.
\end{equation}
By the mean value theorem we have $0<\theta=\theta(t, q, p)<1$ such
that
\begin{eqnarray*}
\Bigl|\frac{\partial H}{\partial p}(t, q,
p)\Bigr|&=&\Bigl|\frac{\partial H}{\partial p}(t, q, p)-
\frac{\partial H}{\partial p}(t, q, \bar p)\Bigr|\\
&=&\biggl|\Bigl[\frac{\partial^2 H}{\partial q_i\partial
q_j}\bigl(t, q, \theta p+ (1-\theta)\bar p\bigr)\Bigr](p-\bar
p)^t\biggr|.
\end{eqnarray*}
Since the first inequality in Proposition A.2  implies
$$
C_1|{\bf u}|\le \biggl|\frac{\partial^2H}{\partial p_i\partial
p_j}(t,q, p){\bf u}\biggr|\le C_2|{\bf u}|\quad\forall {\bf
u}=(u_1,\cdots,u_n)^t\in\R^n,
$$
 using (\ref{e:A.10}) and the inequality $ab\le \frac{\varepsilon}{2}a^2+
\frac{1}{2\varepsilon}b^2\;\forall\varepsilon>0$ we easily get
\begin{eqnarray*}
C_1|p|-C_1C_3\le C_1|p-\bar p|\le\!\!\!\!\!\!
&&\!\!\!\Bigl|\frac{\partial H}{\partial
p}(t, q, p)\Bigr|\le C_2|p-\bar p|\le C_2|p|+ C_2C_3,\\
\frac{C^2_1}{2}|p|^2- 2C^2_1C^2_3\le
\!\!\!\!\!\!&&\!\!\!\Bigl|\frac{\partial H}{\partial p}(t, q,
p)\Bigr|^2\le  2C^2_2|p|^2+ 2C^2_2C^3_3.
\end{eqnarray*}
These two inequalities and Proposition A.2 lead to:
 In any local coordinates
$(q_1,\cdots, q_n)$, the conditions (L2)-(L3) are equivalent to the
fact that there exist constants $0<c<C$, depending on the local
coordinates,  such that
\begin{eqnarray*}
&&cI_n\le \Bigl[\frac{\partial^2 H}{\partial p_i\partial p_j}(t,q,
p)\Bigr] \le
CI_n\quad{\rm and}\\
&&\biggl|\Bigl[\frac{\partial^2 H}{\partial q_i\partial p_j}\bigl(t,
q, p\bigr)\Bigr]\biggr| \le C (1+ |p|),\qquad
\biggl|\Bigl[\frac{\partial^2 H}{\partial q_i\partial q_j}\bigl(t,
q, p\bigr)\Bigr]\biggr| \le C (1+ |p|^2).
\end{eqnarray*}
Proposition A is proved. $\Box$

 \noindent{\bf A.2. An inequality for $C^1$-simplex in $C^1$ Riemannian-Hilbert
manifolds.}   For every integer $q\ge 0$ we denote by $\triangle_q$
the standard closed $q$-dimensional simplex in $\R^q$ with vertices
$e_0=0, e_1,\cdots,e_q$, i.e. $\triangle_0=\{0\}$ and
$$
\triangle_q:=\{(t_1,\cdots, t_q)\in\R^n_{\ge 0}\,|\, t_1+ \cdots +
t_q\le 1\}
$$
with $q\ge 1$. For $1\le i\le q$ denote by
$F^i_q:\triangle_{q-1}\to\triangle_q$ the $i$-th face. Let
$e(s)=(s,\cdots,s)\in\R^q$  with $s\in [0,1]$, ${\hat
e}=e(1/(q+1))$, and $L$ be the straight line passing through $e(0)$
and $\hat e$ successively in $\R^q$, i.e. $L=\{s\hat e\,|\,
s\in\R\}$. Then we have an orthogonal subspace decomposition
$$
\R^q=V_{q-1}\times L,
$$
and each $w\in\triangle_q$  may be uniquely written as $w=(v,
s_0)\in [V_{q-1}\times L]\cap\triangle_q$. Denote by $l(v)$ the
intersection segment of $\triangle_q$ with the straight line passing
through $w$ and parallel to $L$, i.e. $l(v)=\{w+ s\hat e\in
\triangle_q\,|\,s\in\R\}=\{(v, s)\,|\, s_1\le s\le s_2\}$ for some
$s_1\le s_0$ and $s_2\ge s_0$. Clearly, each $l(v)$ has length no
more than $\sqrt{q}/2$. \vspace{2mm}

Let $({\cal M}, \langle, \rangle)$ be a $C^1$ Riemannian-Hilbert
manifold and $\|\cdot\|$ be the induced Finsler metric.  For
$\phi\in C(\triangle_q, {\cal M})$ and each $w=(v, s_0)\in
[V_{q-1}\times L]\cap\triangle_q$ we define
$$
\widetilde\phi_v:l(v)\to {\cal M}, s\mapsto\phi(v,s).
$$
 If  $\phi\in C^1(\triangle_q,
{\cal M})$, i.e. $\phi$ can be extended into a $C^1$-map from some
open neighborhood of $\triangle_q$ in $\R^q$ to ${\cal M}$, then
there exists a constant $c=c(\phi)>0$ such that
$$
\Bigl\|\frac{\partial}{\partial s}\phi(v, s)\Bigr\|^2\le
c(\phi),\quad\forall (v, s)\in\triangle_q.
$$
So for any $(v,s)\in\triangle_q$ we get
\begin{equation}\label{e:A.11}
\int_{l(v)}\Bigl\|\frac{d}{ds}\tilde\phi_v(s)\Bigr\|^2ds\le
c(\phi){\rm Length}(l(v))\le \frac{\sqrt{q}}{2}c(\phi).
\end{equation}
 Now
consider the case ${\cal M}=E_{\tau}=W^{1,2}(S_{\tau}, M)$ with the
Riemannian metric given by (\ref{e:1.15}). Using the local
coordinate chart in (\ref{e:3.8}) it is easy to prove\vspace{2mm}

\noindent{\bf Lemma A.3.}\quad{\it For each $t\in S_\tau$ the
evaluation map
$$
{\bf EV}_t: W^{1,2}(S_\tau, M)\to M,\;\gamma\mapsto\gamma(t)
$$
is continuous and maps $W^{1,2}$-curves in $E_\tau$ to
$W^{1,2}$-curves in $M$.} \vspace{2mm}

\noindent{\bf Proof.}\quad We only need to prove the case $M=\R^n$.
Let $[a,b]\to \gamma(s)$ be a $W^{1,2}$-curve in $W^{1,2}(S_\tau,
\R^n)$. Then $\xi(s):=\frac{d}{ds}\gamma(s)$ is a $W^{1,2}$-vector
field along $\gamma(s)$. Since $T_{\gamma(s)}W^{1,2}(S_\tau,
\R^n)=W^{1,2}(S_\tau, \R^n)$, $\xi(s)\in W^{1,2}(S_\tau, \R^n)$ and
$$
\lim_{\epsilon\to
0}\left\|\frac{\gamma(s+\epsilon)-\gamma(s)}{\epsilon}-\xi(s)\right\|_{W^{1,2}(S_\tau,
\R^n)}=0.
$$
Carefully checking the proof of Proposition 1.2.1(ii) in \cite[pp.
9]{Kl} one easily derives
\begin{equation}\label{e:A.12}
\|\eta\|_{C^0}\le\sqrt{\frac{1+\tau}{\tau}}\|\eta\|_{W^{1,2}}\quad\forall
\eta\in W^{1,2}(S_\tau,\R^n).
\end{equation}
Hence we get
$$
\lim_{\epsilon\to
0}\left\|\frac{\gamma(s+\epsilon)(t)-\gamma(s)(t)}{\epsilon}-\xi(s)(t)\right\|_{\R^n}=0
$$
uniformly in $t$. This means that $[a, b]\to M,\; s\to {\bf
EV}_t(\gamma(s))$, is differentiable and
\begin{equation}\label{e:A.13}
\frac{d}{ds}{\bf EV}_t(\gamma(s))=\xi(s)(t)\quad\hbox{at each}\quad
s\in [a,b].
\end{equation}
Fix a $\epsilon>0$ such that
$$
\left\|\frac{\gamma(s+\epsilon)-\gamma(s)}{\epsilon}-\xi(s)\right\|_{W^{1,2}(S_\tau,
\R^n)}\le\sqrt{\frac{1+\tau}{\tau}}.
$$
By (\ref{e:A.12})  we get
$$
\left\|\frac{\gamma(s+\epsilon)(t)-\gamma(s)(t)}{\epsilon}-\xi(s)(t)\right\|^2_{\R^n}\le
1\quad\forall t\in\R.
$$
It follows that for any $s\in [a, b]$,
\begin{eqnarray*}
\left\|\xi(s)(t)\right\|^2_{\R^n}&\le&
2\Biggl[\left\|\frac{\gamma(s+\epsilon)(t)-\gamma(s)(t)}{\epsilon}-\xi(s)(t)\right\|^2_{\R^n}\\
&&\qquad +\left\|\frac{\gamma(s+\epsilon)(t)-\gamma(s)(t)}{\epsilon}\right\|^2_{\R^n}\Biggr]\\
&\le& 2\left[1 +
\frac{1}{\epsilon^2}\left\|\gamma(s+\epsilon)(t)-\gamma(s)(t)\right\|^2_{\R^n}\right]\\
&\le& 2\left[1 +
\frac{1+\tau}{\tau\epsilon^2}\left\|\gamma(s+\epsilon)-\gamma(s)\right\|^2_{W^{1,2}(S_\tau,
\R^n)}\right].
\end{eqnarray*}
Here the final inequality is due to (\ref{e:A.12}).  Hence
$\int^b_a\left\|\xi(s)(t)\right\|^2_{\R^n} ds<+\infty$, and thus
$\int^b_a\left\|\frac{d}{ds}{\bf EV}_t(\gamma(s))\right\|^2_{\R^n}
ds<+\infty$ because of (\ref{e:A.13}). $\Box$\vspace{2mm}

 For a singular simplex $\sigma$ from $\triangle_q$ to $E_\tau$
and every $w=(v, s_0)\in\triangle_q$, define curves
\begin{equation}\label{e:A.14}
\widetilde\sigma_v^t:l(v)\to M,\, s\mapsto {\bf
EV}_t(\widetilde\sigma_v(s))=\widetilde\sigma_v(s)(t)
\end{equation}
for each $t\in S_\tau$. The curve $\widetilde\sigma_v^0$ is called
the {\bf initial point curve}. Suppose that  $\sigma\in
C^1(\triangle_q, E_\tau)$. Then $\tilde\sigma_v\in C^1(l(v),
E_\tau)$, and by (\ref{e:A.11}) there exists a positive constant
$c(\sigma)$ such that
\begin{equation}\label{e:A.15}
\int_{l(v)}\left\|\frac{d}{ds}\widetilde\sigma_v(s)
\right\|^2_{W^{1,2}(\widetilde\sigma_v(s)^\ast TM)}ds\le
\frac{\sqrt{q}}{2}c(\sigma)
\end{equation}
 for any $(v,s)\in\triangle_q$,
where  $\frac{d}{ds}\widetilde\sigma_v(s)\in
T_{\widetilde\sigma_v(s)
 }E_\tau=W^{1,2}(\widetilde\sigma_v(s)^\ast TM)$.
Specially, by Lemma A.3 we get each $\widetilde\sigma_v^t\in
W^{1,2}(l(v), M)$ for any $t$. As in the proof of Proposition
1.2.1(ii) in \cite[pp. 9]{Kl} one can easily derive that
$$
\|\xi\|_{C^0(\gamma^\ast
TM)}\le\sqrt{\frac{1+\tau}{\tau}}\|\xi\|_{W^{1,2}(\gamma^\ast TM)}
$$
for any $\gamma\in W^{1,2}(S_\tau, M)$ and $\xi\in
W^{1,2}(\gamma^\ast TM)$. Applying to $\gamma=\widetilde\sigma_v(s)$
and $\xi=\frac{d}{ds}\widetilde\sigma_v(s)$
 we get
\begin{equation}\label{e:A.16}
\left\|\frac{d}{ds}\widetilde\sigma_v(s)
\right\|^2_{C^0(\widetilde\sigma_v(s)^\ast TM)}\le
\frac{1+\tau}{\tau} \left\|\frac{d}{ds}\widetilde\sigma_v(s)
\right\|^2_{W^{1,2}(\widetilde\sigma_v(s)^\ast TM)}.
\end{equation}
Moreover, it follows from (\ref{e:A.13}) and (\ref{e:A.14}) that
$$
\left(\frac{d}{ds}\widetilde\sigma_v(s)\right)(t)
=\frac{d}{ds}\widetilde\sigma_v^t(s)=\frac{d}{ds}(\widetilde\sigma_v(s)(t))
\in T_{\widetilde\sigma_v(s)(t)}M
$$
for all $s\in [a, b]$ and $t\in S_\tau$.  Hence for any $t\in
S_\tau$, we can derive from (\ref{e:A.16}) that
\begin{eqnarray*}
\left\|\frac{d}{ds}\widetilde\sigma_v^t(s)\right\|^2_{T_{\widetilde\sigma_v(s)(t)}M}&=&
\left\|\left(\frac{d}{ds}\widetilde\sigma_v(s)\right)(t)\right\|^2_{T_{\widetilde\sigma_v(s)(t)}M}\\
&\le&\left(\max_{t\in
S_\tau}\left\|\left(\frac{d}{ds}\widetilde\sigma_v(s)\right)(t)\right\|_{T_{\widetilde\sigma_v(s)(t)}M}\right)^2\\
&=&\left\|\frac{d}{ds}\widetilde\sigma_v(s)
\right\|^2_{C^0(\widetilde\sigma_v(s)^\ast TM)}\\
&\le& \frac{1+\tau}{\tau} \left\|\frac{d}{ds}\widetilde\sigma_v(s)
\right\|^2_{W^{1,2}(\widetilde\sigma_v(s)^\ast TM)}.
\end{eqnarray*}
This and (\ref{e:A.15}) together give the following generalization
of
 \cite[Lem. 2.3]{Lo2}.\vspace{3mm}

  \noindent{\bf Lemma A.4.}\quad{\it If  $\sigma\in C^1(\triangle_q, E_\tau)$, for every
$w=(v, s_0)\in\triangle_q$,  it holds that
$$\int_{l(v)}\left\|\frac{d}{ds}\widetilde\sigma^0_v(s)\right\|^2_{T_{\widetilde\sigma^0_v(s)}M} ds\le
\frac{(1+\tau)\sqrt{q}}{2\tau}c(\sigma).
$$ }

\end{document}